\documentclass{article}

\usepackage{arxiv}

\usepackage[utf8]{inputenc} % allow utf-8 input
\usepackage[T1]{fontenc}    % use 8-bit T1 fonts
\usepackage[hypertexnames=true]{hyperref}       % hyperlinks
\usepackage{url}            % simple URL typesetting
\usepackage{booktabs}       % professional-quality tables
\usepackage{amsmath, amsfonts}       % blackboard math symbols
\usepackage{amsthm}
\usepackage{nicefrac}       % compact symbols for 1/2, etc.
\usepackage{microtype}      % microtypography
\usepackage{cleveref}       % smart cross-referencing
\usepackage{graphicx}
\usepackage{doi}

\usepackage{amssymb}
\usepackage{algorithmic}
\usepackage[ruled,linesnumbered,vlined]{algorithm2e}
\usepackage{array}
\usepackage[caption=false,font=normalsize,labelfont=sf,textfont=sf]{subfig}
\usepackage{textcomp}
\usepackage{stfloats}
\usepackage{verbatim}
\usepackage{cite}
\usepackage{multirow} % Required for multirows
\usepackage{scalerel}
\usepackage{tikz}
\usepackage[title,titletoc]{appendix}
\usepackage{multirow}
\usepackage{amssymb}
\usepackage{pifont}

\DeclareMathOperator*{\argmin}{\arg\min}

\newcommand{\cmark}{\ding{51}}%
\newcommand{\xmark}{\ding{55}}%

\theoremstyle{plain}

\theoremstyle{definition}

\newtheorem{definition}{Definition}
\theoremstyle{remark}

\title{A Bi-Objective Optimization Based Acquisition Strategy for Batch Bayesian Global Optimization}

\author{\hspace{1mm}Francesco Carciaghi\\
	Global Optimization Laboratory (GOL) \\
	Department of Information Engineering \\
	University of Florence \\
	Via di Santa Marta, 3, 50139, Florence, Italy \\
	\texttt{francesco.carciaghi@unifi.it} \\
	\And
	\hspace{1mm}Simone Magistri \\
	Global Optimization Laboratory (GOL) \\
	Department of Information Engineering \\
	University of Florence \\
	Via di Santa Marta, 3, 50139, Florence, Italy \\
	\texttt{simone.magistri@unifi.it} \\
	\And
	\hspace{1mm}Pierluigi Mansueto \\
	Global Optimization Laboratory (GOL) \\
	Department of Information Engineering \\
	University of Florence \\
	Via di Santa Marta, 3, 50139, Florence, Italy \\
	\texttt{pierluigi.mansueto@unifi.it} \\
	\And
	\hspace{1mm}Fabio Schoen \\
	Global Optimization Laboratory (GOL) \\
	Department of Information Engineering \\
	University of Florence \\
	Via di Santa Marta, 3, 50139, Florence, Italy \\
	\texttt{fabio.schoen@unifi.it} \\
}

\renewcommand{\headeright}{Francesco Carciaghi and Simone Magistri and Pierluigi Mansueto and Fabio Schoen}
\renewcommand{\undertitle}{}
\renewcommand{\shorttitle}{}

\hypersetup{
pdftitle={A Bi-Objective Optimization Based Acquisition Strategy for Batch Bayesian Global Optimization},
pdfsubject={},
pdfauthor={Francesco Carciaghi, Simone Magistri, Pierluigi Mansueto, Fabio Schoen},
pdfkeywords={Global Optimization, Batch Bayesian Methods, Bi-Objective Optimization, Acquisition Function},
}

\begin{document}
\maketitle

\begin{abstract}
	In this paper, we deal with batch Bayesian Optimization (Bayes-Opt) problems over a box and we propose a novel bi-objective optimization (BOO) acquisition strategy to sample points where to evaluate the objective function. The BOO problem involves the Gaussian Process posterior mean and variance functions, which, in most of the acquisition strategies from the literature, are generally used in combination, frequently through scalarization. However, such scalarization could compromise the Bayes-Opt process performance, as getting the desired trade-off between exploration and exploitation is not trivial in most cases. We instead aim to reconstruct the Pareto front of the BOO problem based on optimizing both the posterior mean as well as the variance, thus generating multiple trade-offs without any a priori knowledge. The reconstruction is performed through the Non-dominated Sorting Memetic Algorithm (NSMA), recently proposed in the literature and proved to be effective in solving hard MOO problems. Finally, we present two clustering approaches, each of them operating on a different space, to select potentially optimal points from the Pareto front. We compare our methodology with well-known acquisition strategies from the literature, showing its effectiveness on a wide set of experiments.
\end{abstract}

\keywords{Global Optimization \and Batch Bayesian Methods \and Bi-Objective	Optimization \and Acquisition Function}
\MSCs{62C10 \and 90C29 \and 90C30}

\section{Introduction}
\label{sec::introduction}

Bayesian Optimization (Bayes-Opt) is a set of highly relevant heuristic optimization approaches whose aim is to approximately find a global optimum of the following problem:
\begin{equation}
	\label{eq::bo-prob}
	\min_{x \in \Omega} f(x),
\end{equation}
where $\Omega$ is typically a simple feasible set, i.e., a box $\{x \in \mathbb{R}^n \mid x \in [l, u],\ l, u \in \mathbb{R}^n \text{ s.t. } l^i \le u^i\}$, and $f: \mathbb{R}^n \rightarrow \mathbb{R}$ is a continuous, expensive-to-evaluate objective function. We assume that $f(\cdot)$ is evaluated without noise. The expensiveness of function evaluations may be caused by different factors: an high consumption of resources (CPU/GPU or clock time), monetary costs, investments, etc.. Application of such context can be found in various fields: hyper-parameters tuning for machine learning algorithms \cite{NIPS2012_05311655}, engineering systems design \cite{forrester2008engineering, jones1998efficient}, environmental models calibration \cite{shoemaker07}, reinforcement learning \cite{brochu2010tutorial, lizotte2007automatic} and materials/drug design \cite{Frazier2016, negoescu11} are just a few examples. Thus, the development of strategies to generate a good approximation of a global optimum, within a very limited budget of function evaluations, is of great practical interest and Bayes-Opt stands out as one of the most frequently adopted approaches.

Early papers on the subject can be traced back to \cite{Kushner64}, whose author can be considered as one of the fathers of Bayes-Opt, or to \cite{Mockus75}, who greatly contributed to the diffusion of the Bayesian Optimization ideas. Years later, a quite specific Bayes-Opt approach was presented in \cite{jones1998efficient} in an extremely clean and elegant way, attracting much interest from the communities. In particular, the attention of the machine learning community explains the increased importance towards black-box optimization of expensive functions. In recent years, Bayes-Opt was employed both in problems like \eqref{eq::bo-prob} and in more difficult ones (see, e.g., \cite{Audet2022}, where instances with both expensive-to-evaluate functions and constraints are taken into account, and \cite{binois22, swersky2017improving, oh18}, where Bayes-Opt have been considered in high-dimensional scenarios). A nice survey of recent approaches can be found in \cite{Frazier18}, while a detailed introduction to the main aspects of the theory is proposed in \cite{Garnett2023}. Furthermore, many software tools have been recently developed, with \texttt{BoTorch} \cite{Botorch} being considered as one of the most efficient ones for modern Bayesian Optimization.

We recall that Bayes-Opt mainly consists in two phases: first, a stochastic model of the objective function, i.e., a Gaussian Process, is defined \textit{a priori}, and, then, it is updated after observing $f(\cdot)$ at specific feasible points. The resulting \textit{a posteriori} model contains the updated knowledge we can use in order to predict the objective function value at new not-yet-visited feasible points. We refer the interested reader to \cite{Rsamussen06} for a deep and exhaustive analysis on how to model a surrogate function through Gaussian Processes. Given the stochastic model, the problem now becomes that of exploiting information contained in the model to guide the search towards likely locations of a global optimum. Assume that, after the objective function $f(\cdot)$ has been sampled at $k$ feasible points $x_1, \ldots, x_k \in \mathbb{R}^n$, a Gaussian model has been fit; let us denote by $\mu_k(\cdot)$ and $\sigma_k(\cdot)$ the posterior mean and standard deviation of $f(\cdot)$. One of the most frequently employed criteria consists on optimizing a suitably defined \textit{acquisition} function, such as the \textit{Expected Improvement} (\texttt{EI}) \cite{Mockus75} or the \textit{Lower Confidence Bound} (\texttt{LCB}) \cite{Srinivas_2012}, seeking a compromise between points close to already observed low valued feasible solutions (low $\mu_k(\cdot)$ values) and points in high uncertainty zones (high $\sigma_k(\cdot)$ values) typically located as far as possible from past observations. Here, we deal with Batch Bayesian Optimization, a sub-field of Bayes-Opt where, at each iteration, $q > 1$ feasible solutions are chosen through an \textit{acquisition} criterion to be evaluated on $f(\cdot)$. The idea of using Bayes-Opt for multiple simultaneous function evaluations is not new in the literature: it was exploited, e.g., in \cite{wang_parallel_2020}, based on the exploitation of the \texttt{q-EI}, or in \cite{wilson2017reparameterization}, where the \texttt{q-LCB} formulation can be found. It might seem to be contradictory, in such an expensive setting, the evaluation of the objective function several times before updating the Gaussian Process. However, in many situations, parallel computing machines are available; moreover, the idea of simultaneously and independently computing the objective function is worth considering when we try to reduce the wall-clock time of the overall procedure. 

In this paper, we propose a novel \textit{acquisition} strategy, based on reconstructing the Pareto front of a Bi-Objective Optimization (BOO) problem based on $\mu_k(\cdot)$ and $\sigma_k(\cdot)$. Unlike the mentioned approaches, we do not rely on combinations of the two posterior functions; indeed, such a way of proceeding could compromise the Bayes-Opt process performance, being difficult in some cases to get the desired trade-off between exploitation and exploration. Reconstructing the Pareto front instead leads to the generation of multiple compromises without any a priori setting, thus allowing to choose, based on some criteria, the most suitable points among many. In general, Multi-Objective Optimization (MOO) received much attention in the last decades; many approaches, mainly extensions of classical scalar optimization techniques, have been developed: \cite{fliege2000steepest, drummond2004projected, fliege2009newton, prudente_quasi-newton_nodate, LapucciLBFGS2023, cocchi2020convergence, custodio2011direct, deb2002, NSMA} just represent a non-exhaustive list of this type of methodologies. In particular, some of the most recent ones are capable of managing multiple solutions at each iteration, finally providing an approximation, as accurate as possible, of the Pareto front. 

The idea of an acquisition strategy based on MOO has been considered in some recent works from the literature. In \cite{death21}, $\mu_k(\cdot)$ and $\sigma_k(\cdot)$ are used as objective functions in a BOO problem; the latter is then solved with the \textit{Non-dominated Sorting Genetic Algorithm II} (\texttt{NSGA-II}) \cite{deb2002}, which is one of the most representative evolutionary algorithms for solving hard MOO problems; after reconstructing the Pareto front of the considered BOO problem, the new points where to evaluate $f(\cdot)$ are randomly selected. The same algorithmic scheme characterizes the work in \cite{gupta_exploiting_nodate}, where a variant of \texttt{NSGA-II} \cite{Deb2011} is employed. In \cite{CHEN2023101293}, \texttt{NSGA-II} and some of its variants are used to solve a BOO problem where, at each iteration of the Bayes-Opt procedure, the objectives are dynamically chosen based on an ad-hoc defined metric from a pool of well-known acquisition functions, such as \texttt{q-EI} and \texttt{q-LCB}; then, the candidate solutions selection involves a clustering approach employed in the domain space of the Pareto front points. In \cite{Feng2015}, the BOO problem is composed by the two terms of the \texttt{EI}; the method is the \textit{Multi-objective Evolutionary Algorithm based on Decomposition} (\texttt{MOEA/D}) \cite{zhang07}, based on a decomposition strategy to turn the BOO problem into $N$ scalar ones, each of them characterized by a specific weighted sum of the objectives; the selection of points is performed through clustering of the weights vectors defining the $N$ scalar problems. Finally, in \cite{lyu_batch_2018} the authors propose an MOO problem composed by three well-known acquisition functions, i.e., \texttt{EI}, \texttt{LCB} and \textit{Probability of Improvement} (\texttt{PI}) \cite{Kushner64}; the employed algorithm is the \textit{Differential Evolution for Multiobjective Optimization} (\texttt{DEMO}) \cite{demo_paper}; candidate points are again randomly chosen from the resulting Pareto front approximation.

Our novel acquisition strategy differs from the mentioned ones for the following aspects.
\begin{itemize}
	\item At each iteration of the batch Bayes-Opt procedure, we aim to reconstruct, as accurately as possible, the Pareto front of the BOO problem
	\begin{equation}
		\label{eq::acq-f-biob}
		\min_{x \in \Omega} (\mu_k(x), -\sigma_k^2(x))^\top.
	\end{equation}
	These two functions are the basis for some of the most employed acquisition functions (e.g., \texttt{EI}, \texttt{LCB}): optimizing $\mu_k(\cdot)$ leads to exploit the current information on the objective function $f(\cdot)$, while $\sigma_k^2(\cdot)$ let us highlight regions where we know little or nothing about the objective function. This duality leads to a resulting Pareto front with a clear meaning, being composed by different tradeoffs between exploitation and exploration of the current information about $f(\cdot)$.
	%	without defining any a priori weighted scalarization of the two objective functions.
	\item We employ the recently proposed \textit{Non-dominated Sorting Memetic Algorithm} (\texttt{NSMA}) \cite{NSMA} to solve problem \eqref{eq::acq-f-biob}; \texttt{NSMA} combines the good exploration capabilities of \texttt{NSGA-II} with tools typical of gradient-based methodologies; as it is shown in \cite{NSMA}, \texttt{NSMA} was compared with well-known algorithms for MOO, including \texttt{NSGA-II}, widely used in the Batch Bayesian literature; the effectiveness and efficiency of \texttt{NSMA} was particularly tested in retrieving good Pareto front approximations of difficult MOO problems.
	\item We propose two clustering methodologies, each of them operating on a different space, to choose from the Pareto front the new points where to evaluate $f(\cdot)$; in particular, we analyze the potential effectiveness of such clustering approaches on high-dimensional scenarios ($n \ge 20$), where the curse of dimensionality usually prevents the employed acquisition functions to find good candidate solutions where to evaluate the objective function $f(\cdot)$.
\end{itemize}

The remainder of the manuscript is organized as follows. In Section \ref{sec::preliminaries}, we recall the main batch Bayes-Opt concepts, focusing on two main aspects: the Gaussian Process regression and the acquisition phase. In this section, we also review the main Pareto's theory concepts. In Section \ref{sec::bi-obj-acq-f-for-BO}, we describe our novel BOO based acquisition strategy, reporting a brief description of \texttt{NSMA}, adapted to the Bayes-Opt context, and the algorithmic scheme of the proposed clustering approaches. In Section \ref{sec::computational-experiments}, we analyze the results of thorough computational experiments, highlighting the good performance of our acquisition methodology w.r.t.\ well-known ones from the batch Bayes-Opt literature. Finally, in Section \ref{sec::concluding-remarks} we provide some concluding remarks.

\section{Preliminaries}
\label{sec::preliminaries}

In this section, we first recall the principal batch Bayes-Opt concepts, including the Gaussian Process regression and the acquisition phase. Then, we review some MOO key concepts from the Pareto's theory.

\subsection{Batch Bayesian Optimization Overview}
\label{subsec::overview_bo}

In Algorithm \ref{alg::Basic_BO}, we report the main steps of a batch Bayes-Opt general scheme.

\begin{algorithm}[!h]
	\caption{Batch Bayes-Opt} \label{alg::Basic_BO}
	Input: $f:\mathbb{R}^n\to \mathbb{R}$, $\Omega$ feasible set, $N_i \in \mathbb{N}^+$ number of initial points, $N_M > N_i \in \mathbb{N}$ maximum number of evaluated points, $q > 1 \in \mathbb{N}$ batch size (number of evaluated points at each iteration), $\mu_0(\cdot)$ mean function, $\Sigma_0(\cdot, \cdot)$ kernel.\\
	Let $x_{1:N_i}$ be $N_i$ initial points in $\Omega$ \label{line::first-points}\\
	Let $f(x_{1:N_i})$ be the values of $f(\cdot)$ at $x_{1:N_i}$ \label{line::first-evaluation}\\
	Let $k = N_i$\\
	\While{$k \le N_i + N_M$}{
		Let the \textit{prior distribution} on $f(x_{1:k})$ be 
		\begin{equation}
			\label{eq::prior}
			f(x_{1:k}) \sim \mathcal{N}(\mu_0(x_{1:k}), \Sigma_0(x_{1:k}, x_{1:k}))
		\end{equation}\label{line::prior}\\
		Let the \textit{posterior distribution} on $f(x)$ be
		\begin{equation}
			\label{eq::posterior}
			\begin{aligned}
				f(x) \mid &f(x_{1:k}) \sim \mathcal{N}\left(\mu_k(x), \sigma^2_k(x)\right)\\
				&\mu_k(x) = \Sigma_0(x, x_{1:k})\Sigma_0(x_{1:k}, x_{1:k})^{-1}(f(x_{1:k}) - \mu_0(x_{1:k})) + \mu_0(x)\\
				&\sigma^2_k(x) = \Sigma_0(x, x) - \Sigma_0(x, x_{1:k})\Sigma_0(x_{1:k}, x_{1:k})^{-1}\Sigma_0(x_{1:k}, x)
			\end{aligned}
		\end{equation}\label{line::posterior}\\
		%\State Let $q_b = \max\{\tilde{q} \in \mathbb{N} \mid \tilde{q} \le q \land k + \tilde{q} \le N_i + N_M\}$ \label{line::tilde_q}
		Let $x_{k+1},\ldots,x_{k + q}$ be the new points found by minimizing an acquisition function based on the posterior distribution \label{line::acq_f}\\
		Evaluate $f(\cdot)$ at $x_{k+1},\ldots,x_{k + q}$ \label{line::evaluation_of_f}\\
		Let $k = k + q$
	}	
	\Return $x^\star$ s.t. $f(x^\star) \le f(x),\ \forall x \in x_{1:N_i + N_M}$\label{line::return_best}
\end{algorithm}

The procedure starts from $N_i$ initial points in the feasible set $\Omega$ (Line \ref{line::first-points}). For these points, we assume to know the value of $f(\cdot)$ (Line \ref{line::first-evaluation}). Moreover, as anticipated in Section \ref{sec::introduction}, we suppose that each $f(\cdot)$ evaluation is performed in a noise-free setting.
As in \cite{Frazier18}, we make use of the following compact notations: given $k \ge 1$,
\begin{equation}
	\label{eq::compact-notations}
	\begin{gathered}
		x_{1:k} = [x_1,\ldots,x_k]^\top \in \mathbb{R}^{k \times n},\\
		f(x_{1:k}) = [f(x_1),\ldots,f(x_k)]^\top \in \mathbb{R}^k.
	\end{gathered}
\end{equation}

Now, in order to model the objective function $f(\cdot)$, the Gaussian Process regression (Lines \ref{line::prior}-\ref{line::posterior}) is carried out. Let us suppose that we have already evaluated $f(\cdot)$ at $k$ points. Given these evaluations, we can define the \textit{prior distribution} on $f(x_{1:k})$ as in Equation \eqref{eq::prior}, where $\mathcal{N}(\mu_0(x_{1:k}), \Sigma_0(x_{1:k}, x_{1:k}))$ denotes the multivariate normal distribution with $\mu_0(\cdot)$ and $\Sigma_0(\cdot, \cdot)$ being the \textit{mean} and \textit{covariance} (\textit{kernel}) functions. The kernel is chosen such that close points have a large positive correlation, based on the idea that similar function values take place in them more likely than in points distant from each other. Similarly to \eqref{eq::compact-notations}, we indicate with $\mu_0(x_{1:k})$ and $\Sigma_0(x_{1:k}, x_{1:k})$ the vector and matrix, respectively, derived from evaluating these two functions at $x_1,\ldots, x_k$, i.e., $\mu_0(x_{1:k}) = [\mu_0(x_1),\ldots,\mu_0(x_k)]^\top$ and
\begin{equation*}
	\Sigma_0(x_{1:k}, x_{1:k}) = 
	\begin{bmatrix} 
		\Sigma_0(x_1, x_1) & \cdots & \Sigma_0(x_1, x_k) \\ 
		\vdots & \ddots & \vdots \\ 
		\Sigma_0(x_k, x_1) & \cdots & \Sigma_0(x_k, x_k)
	\end{bmatrix}.
\end{equation*}
The matrix $\Sigma_0(x_{1:k}, x_{1:k})$ is positive semi-definite, regardless the considered points.

The most commonly used mean function is a constant value: $\mu_0(\cdot) = \mu$. As for the kernel, a popular choice is the \textit{Màtern} kernel \cite{Rsamussen06}; here, we report the formulation that can be found in \cite{Frazier18}: given $x, y \in \mathbb{R}^n$,
\begin{equation*}
	\Sigma_0(x, y) = \alpha_0\frac{2^{1-\nu}}{\Gamma(\nu)}(\sqrt{2\nu}\|x - y\|_{\alpha_{1:n}})^\nu K_\nu(\sqrt{2\nu}\|x - y\|_{\alpha_{1:n}}),
\end{equation*}
with $\Gamma(\cdot)$ being the \textit{gamma} function, $\|x - y\|_{\alpha_{1:n}} = \sqrt{\sum_{i=1}^{n} \alpha_i(x^i - y^i)^2}$ and $K_\nu(\cdot)$ being the \textit{Modified Bessel} function. The hyper-parameters contained in these formulas, i.e., $\mu, \alpha_{0:n}, \nu$, can be found by \textit{Maximum Likelihood Estimation} (\texttt{MLE}) given the observations $f(x_{1:k})$. For more information on this procedure or on other alternative formulations for $\mu_0(\cdot)$ and $\Sigma_0(\cdot, \cdot)$, the reader is referred to \cite{Frazier18}.

We can then compute the \textit{posterior} probability distribution, which describes potential values of the objective function $f(\cdot)$ at a new candidate point (Line \ref{line::posterior}). Using Bayes' rule, the conditional distribution of $f(\cdot)$ given $f(x_{1:k})$ can be computed as a normal distribution $\mathcal{N}(\mu_k(\cdot), \sigma_k^2(\cdot))$ \eqref{eq::posterior}, where $\mu_k(\cdot)$ and $\sigma^2_k(\cdot)$ are the \textit{posterior} mean and variance functions. The mean function $\mu_k(\cdot)$ works as a point estimate of $f(\cdot)$ and typically assumes small values in solutions close to already-observed low-valued points. On the other hand, $\sigma_k^2(\cdot)$ allows to quantify the uncertainty we have on each feasible set region: the farther we move from the already-observed points, the more the uncertainty on the values assumed by $f(\cdot)$ grows. The $\mu_k(\cdot)$ and $\sigma_k^2(\cdot)$ formulations require specific methodologies to be computed: for more details on the topic, we refer the reader to \cite{Rsamussen06, Frazier18}.

Given the posterior distribution, we then minimize an \textit{acquisition function} $a_k: \mathbb{R}^{q \times n} \rightarrow \mathbb{R}$ (Line \ref{line::acq_f}) in order to find $q > 1$ new candidate solutions $x_{k+1},\ldots, x_{k + q}$ where to evaluate the objective function $f(\cdot)$ (Line \ref{line::evaluation_of_f}):
\begin{equation}
	\label{eq::acq-f-prob}
	\min_{\tilde{x}_{k + 1},\ldots, \tilde{x}_{k + q} \in \Omega} a_k(\tilde{x}_{k + 1:k+q}).
\end{equation}
Unlike $f(\cdot)$, the acquisition functions are usually inexpensive to evaluate; moreover, it is feasible to compute their first and second derivatives, which is not possible with $f(\cdot)$. This fact allows the acquisition functions to be optimized using first-order or second-order optimization methods. A solver that has proved to perform well \cite{Frazier18} is the \textit{Limited-Memory Broyden-Fletcher-Goldfarb-Shanno} (\texttt{L-BFGS}) \cite{Liu1989} variant designed for box-constrained problems, i.e., \texttt{L-BFGS-B} \cite{zhu97}.

%Note that $q$ may be such that $k + q > N_i + N_M$. In these cases, in order to have the total number of evaluated points exactly equal to $N_M$ (the maximum number of allowed $f(\cdot)$ evaluations), we could find $q_b$ new candidate solutions, with $q_b$ defined such that $q_b < q$ and $k + q_b = N_i + N_M$. Here, we remind that the values of $f(\cdot)$ at $x_{1:N_i}$ are supposed to be known and, thus, these observations do not count in the total of $N_M$ evaluations.

When the number of maximum objective function evaluations $N_M$ is reached, the procedure is stopped and the best found solution $x^\star$ is returned (Line \ref{line::return_best}).

In the next section, we review some of the most employed acquisition functions in the Bayes-Opt literature, reporting their formulations and a brief description of each of them.

\subsubsection{Expected Improvement and Lower Confidence Bound}
\label{subsubsec::acq-f}

In this section, we review two of the most commonly used acquisition functions: the \textit{Expected Improvement} and the \textit{Lower Confidence Bound}.

Using the notation introduced in Section \ref{subsec::overview_bo}, the \textit{Expected Improvement} (\texttt{EI}) \cite{Mockus75, jones1998efficient} indicates the expected value of the improvement w.r.t.\ $f_k^\star = \min_{i \in \{1,\ldots,k\}}f(x_i)$, that is, the best found objective function value after the evaluation of $k$ points. Following the steps proposed in \cite{jones1998efficient, clark61}, we can evaluate \texttt{EI} as:
\begin{equation}
	\label{eq::ei}
	\begin{gathered}
		\begin{aligned}
			EI_k(x) &= \mathbb{E}_y\left[\max\{f_k^\star - y, 0\} \bigm\vert x_{1:k}, f(x_{1:k})\right] =\\
			&= (f_k^\star - \mu_k(x))\Phi\left(\frac{f_k^\star - \mu_k(x)}{\sigma_k(x)}\right) + \sigma_k(x)\phi\left(\frac{f_k^\star - \mu_k(x)}{\sigma_k(x)}\right),
		\end{aligned}\vspace{0.2cm}\\
		y \sim \mathcal{N}\left(\mu_k(x), \sigma_k^2(x)\right),
	\end{gathered}
\end{equation}
where $\Phi(\cdot)$ and $\phi(\cdot)$ are, respectively, the cumulative distribution function (CDF) and the probability density function (PDF) of a standard normal random variable.
Since the expected improvement should be reasonably maximized as much as possible, according to \eqref{eq::acq-f-prob}, the problem to solve would be $\min_{x\in\Omega}-EI_k(x)$.

The \textit{Lower Confidence Bound} (\texttt{LCB}) \cite{Srinivas_2012} prefers both points with low $\mu_k(\cdot)$, where we expect to have high reward, that is, good values of $f(\cdot)$, and points with large $\sigma_k(\cdot)$, where $f(\cdot)$ is uncertain. \texttt{LCB} can be evaluated as:

\begin{equation}
	\label{eq::lcb}
	\begin{gathered}
		LCB_k(x) = \mathbb{E}_y\left[\mu_k(x) - |y - \mu_k(x)| \bigm\vert x_{1:k}, f(x_{1:k})\right] = \mu_k(x) - \sqrt{\beta} \sigma_k(x),
		\vspace{0.2cm}\\
		y \sim \mathcal{N}\left(\mu_k(x), \beta\frac{\pi}{2}\sigma_k^2(x)\right),
	\end{gathered}
\end{equation}
where $\beta \in \mathbb{R}^n$ is a constant parameter.

Both \texttt{EI} and \texttt{LCB} seek a trade-off between points where an high reward is expected and points where we have an high uncertainty. In the \texttt{LCB} formula, the desired trade-off can be explicitly controlled by the constant parameter $\beta$: if $\beta \approx 0$, the solver used to optimize \texttt{LCB} will be driven towards points associated to low values of $\mu_k(\cdot)$; on the other hand, if $\beta \gg 0$, the points with high values of $\sigma_k(\cdot)$ will be preferred.

The two formulas \eqref{eq::ei}-\eqref{eq::lcb} are extended to be used for batch Bayes-Opt, where $q$ new feasible not-evaluated-yet solutions are required at each iteration. The formulas of the extensions, called \texttt{q-EI} \cite{ginsbourger:hal-00260579} and \texttt{q-LCB} \cite{wilson2017reparameterization} respectively, are reported in Equations \eqref{eq::qei}-\eqref{eq::qlcb}, with $\Sigma_k(\cdot, \cdot)$ indicating the \textit{posterior covariance} function given $k$ observations of $f(\cdot)$.

\begin{equation}
	\label{eq::qei}
	\begin{gathered}
		q\text{-}EI_k(x_{k+1:k+q}) = \mathbb{E}_y\left[\max\left\{f_k^\star - \max\limits_{i \in \{1,\ldots, q\}}y_i, 0\right\} \bigm\vert x_{1:k}, f(x_{1:k}) \right]\vspace{0.2cm}\\
		y \sim \mathcal{N}\left(\mu_k(x_{k+1:k+q}), \Sigma_k(x_{k+1:k+q}, x_{k+1:k+q})\right)
	\end{gathered}
\end{equation}
\vspace{0.4cm}
\begin{equation}
	\label{eq::qlcb}
	\begin{gathered}
		q\text{-}LCB_k(x_{k+1:k+q}) = \mathbb{E}_y\left[\max\limits_{i \in \{1,\ldots, q\}}\left(\mu_k(x_{k+i}) - |y_i - \mu_k(x_{k+i})|\right) \bigm\vert x_{1:k}, f(x_{1:k}) \right]\vspace{0.2cm}\\
		y \sim \mathcal{N}\left(\mu_k(x_{k+1:k+q}), \beta\frac{\pi}{2}\Sigma_k(x_{k+1:k+q}, x_{k+1:k+q})\right)
	\end{gathered}
\end{equation}

Both \texttt{q-EI} and \texttt{q-LCB} are more challenging than their original versions since the integrals arising from the expected value formulation are more intractable; as a consequence, the optimization process per se is not trivial. For this reason, a variety of methods have been proposed in the literature for the evaluation and the optimization of intractable acquisition functions. In this paper, when performing numerical experiments with \texttt{q-EI} and \texttt{q-LCB}, we rely on one of the most common techniques: Monte Carlo sampling based integration. For the sake of brevity, we do not go into the details of either this methodology or other possible strategies, referring the reader to \cite{wilson2017reparameterization} for more information on the two topics.

\subsection{Optimality Concepts in Multi-Objective Optimization}
\label{subsec::MOO}

In this section, we report some MOO key concepts from Pareto's theory. In the rest of the section, we consider the following MOO problem:
\begin{equation}
	\label{eq::mo-prob}
	\min_{x\in\Omega}\;\mathcal{F}(x)=(f_1(x),\ldots,f_m(x))^\top,
\end{equation}
where $\mathcal{F}:\mathbb{R}^n\to\mathbb{R}^m$ is a vector-valued, component-wise continuously differentiable function. 

%We denote by $J_F(\cdot) \in \mathbb{R}^{m \times n}$ the Jacobian matrix associated with $F(\cdot)$.

In order to introduce some MOO notions, we define a partial ordering of the points in the objectives space $\mathbb{R}^m$: considering two points $u, v \in \mathbb{R}^m$, we have that
\begin{equation*}
	\begin{gathered}
		u < v \iff u_i < v_i \hspace{0.9cm} \forall i \in \{1,\ldots, m\}, \\
		u \le v \iff u_i \le v_i \hspace{0.9cm} \forall i \in \{1,\ldots, m\}, \\
		u \lneqq v \iff u \le v, u \ne v.
	\end{gathered}
\end{equation*}
We say that a point $x \in \Omega$ \textit{dominates} $y \in \Omega$ w.r.t.\ $\mathcal{F}(\cdot)$ if $\mathcal{F}(x) \lneqq \mathcal{F}(y)$.

In MOO, we would like to obtain a point which minimizes all the objectives $f_1(\cdot),\ldots, f_m(\cdot)$ at once. However, such a solution is unlikely to exist. For this reason, we rely on optimality concepts different from the ones for single-objective optimization.

\begin{definition}
	A point $\bar{x}\in\Omega$ is \textit{Pareto optimal} for problem \eqref{eq::mo-prob} if a point $y\in\Omega$ such that $\mathcal{F}(y)\lneqq \mathcal{F}(\bar{x})$ does not exist.
	If the property is satisfied only in $\Omega\cap \mathcal{N}(\bar{x})$, with $\mathcal{N}(\bar{x})$ being a neighborhood of $\bar{x}$, then $\bar{x}$ is \textit{locally Pareto optimal}. 
\end{definition}

In practice, it is difficult to get solutions characterized by the Pareto optimality property. Then, a slightly more affordable condition is introduced.

\begin{definition}
	A point $\bar{x}\in\Omega$ is \textit{weakly Pareto optimal} for problem \eqref{eq::mo-prob} if a point $y\in\Omega$ such that $\mathcal{F}(y)< \mathcal{F}(\bar{x})$ does not exist. If the property is satisfied only in $\Omega\,\cap\, \mathcal{N}(\bar{x})$, with $\mathcal{N}(\bar{x})$ being a neighborhood of $\bar{x}$, then $\bar{x}$ is \textit{locally weakly Pareto optimal}.
\end{definition}
We refer to the set of Pareto optimal solutions as the \textit{Pareto set}, while by \textit{Pareto front} we indicate the image of the Pareto set through $\mathcal{F}(\cdot)$.

Finally, we introduce the concept of \textit{Pareto stationarity}. Under differentiability assumptions, the latter is a necessary condition for all types of Pareto optimality; also assuming the component-wise convexity of $\mathcal{F}(\cdot)$, Pareto stationarity is sufficient for Weak Pareto optimality. 

\begin{definition}
	\label{def::par-stat}
	A point $\bar{x}\in\Omega$ is \textit{Pareto-stationary} for problem \eqref{eq::mo-prob} if we have that 
	\begin{equation*}
		\max_{j\in\{1,\ldots,m\}}\nabla f_j(\bar{x})^\top d\ge 0, \qquad \forall d \in \mathcal{D}(\bar{x}),
	\end{equation*}
	where $\mathcal{D}(\bar{x})=\{v\in\mathbb{R}^n\mid \exists\bar{t}>0:\bar{x}+tv\in \Omega,\; \forall t\in[0,\bar{t}\,]\}$.
	The property can be also compactly re-written as $$\min\limits_{d \in \mathcal{D}(\bar{x})} \max\limits_{j\in\{1,\ldots,m\}} \nabla f_j(\bar{x})^\top d = 0.$$
\end{definition}

Unlike (Weak) Pareto optimality, Pareto stationarity represents a concrete goal for many MOO algorithms based on first-order information of $\mathcal{F}(\cdot)$ \cite{fliege2000steepest}.

\section{A Bi-Objective Optimization Based Acquisition Strategy}
\label{sec::bi-obj-acq-f-for-BO}

In Section \ref{subsubsec::acq-f}, we have reported a brief description of some of the most employed acquisition functions from the Bayes-Opt literature. In each of them, a desired trade-off between exploitation (low $\mu_k(\cdot)$ values) and exploration (high $\sigma_k(\cdot)$ values) is defined, either in an implicit way (\texttt{EI}, \texttt{q-EI}) or in an explicit one through a parameter (\texttt{LCB}, \texttt{q-LCB}). However, deciding \textit{a priori} such exploitation-exploration compromise is not trivial in some cases and the Bayes-Opt procedure performance may be affected by this additional difficulty.

Based on this reasoning, in this section, we want to present a novel acquisition methodology for batch Bayes-Opt. As anticipated in Section \ref{sec::introduction}, our proposal is based on the optimization of the BOO problem \eqref{eq::acq-f-biob} and, in particular, on the reconstruction of its Pareto front. Unlike the acquisition strategies of Section \ref{subsubsec::acq-f}, the exploitation-exploration trade-off is not decided \textit{a priori}; on the contrary, the Pareto front provides multiple solutions representing different trade-offs between $\mu_k(\cdot)$ and $\sigma_k(\cdot)$; having such a more complete, \textit{a posteriori} overview, we can then employ specific criteria to choose the solutions we expect to be optimal for problem \eqref{eq::bo-prob}.

We propose the following methodologies:
\begin{itemize}
	\item in order to reconstruct, as accurately as possible, the Pareto front of the BOO problem \eqref{eq::acq-f-biob}, we employ a recently proposed memetic algorithm for MOO problems, i.e., \texttt{NSMA} \cite{NSMA}; this approach was indeed proved to be efficient and effective on reconstructing Pareto front of hard MOO problems w.r.t.\ other state-of-the-art MOO algorithms, including \texttt{NSGA-II} which is widely employed for this scope in the Bayesian literature;
	\item we propose two clustering approaches, operating on different spaces, for the selection of $q$ Pareto front points that will be the next candidate solutions where to evaluate $f(\cdot)$; after the proposal, we also analyze the potential effectiveness of such approaches on high dimensional scenarios ($n \ge 20$), which is known to be more difficult in the Bayesian literature as the curse of dimensionality usually leads the acquisition functions failing in finding good candidate solutions.
\end{itemize}

\begin{figure*}[!h]
	\centering
	\includegraphics[width=\textwidth]{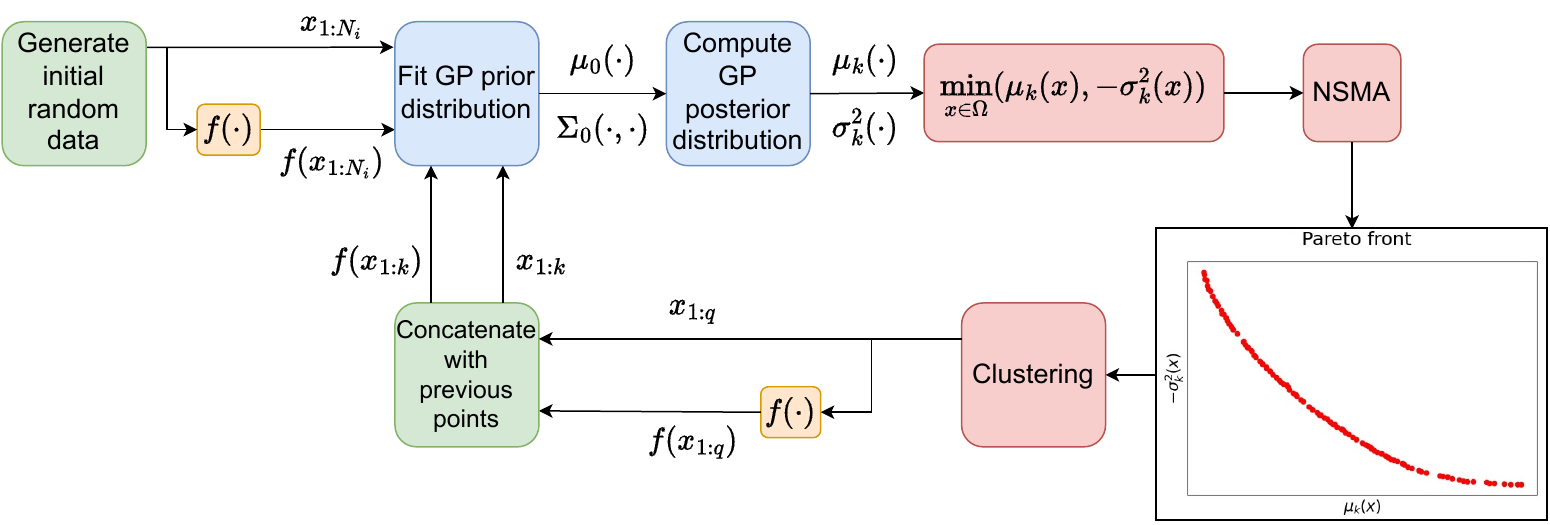}
	\caption{Flowchart of the batch Bayes-Opt procedure equipped with the BOO based acquisition strategy proposed in Section \ref{sec::bi-obj-acq-f-for-BO}.}
	\label{fig:methodology_flow}
\end{figure*}

In Figure \ref{fig:methodology_flow}, we report an overview of the entire batch Bayes-Opt procedure equipped with our BOO based acquisition strategy. The two phases (Pareto front reconstruction, clustering) are separately discussed in the following two sections.

\subsection{The Non-dominated Sorting Memetic Algorithm for Bayes-Opt}
\label{subsec::NSMA-for-BO}

In Algorithm \ref{alg::NSMA}, we report an algorithmic scheme of the \textit{Non-dominated Sorting Memetic Algorithm} (\texttt{NSMA}) \cite{NSMA}. Note that here we only report the main steps of the algorithm, adapted to the Bayes-Opt context; for a deeper understanding of the \texttt{NSMA} mechanisms, which is beyond the scope of this manuscript, we refer the interested reader to the cited paper. In the remainder of the section, we indicate with $\mathcal{F}(\cdot):\mathbb{R}^n\to \mathbb{R}^2$ the vector-valued objective function of problem \eqref{eq::acq-f-biob}, that is, $\mathcal{F}(\cdot) = (f_1(\cdot), f_2(\cdot))^\top = (\mu_k(\cdot), -\sigma_k^2(\cdot))^\top$, with $\mu_k(\cdot)$ and $\sigma_k^2(\cdot)$ being the posterior mean and variance functions (see Section \ref{subsec::overview_bo}, Lines \ref{line::mu_sigma}-\ref{line::F} of Algorithm \ref{alg::NSMA}). 

\begin{algorithm}[!h]
	\caption{Non-dominated Sorting Memetic Algorithm for Bayes-Opt} \label{alg::NSMA}
	Input: $\Omega$, $X^i \subset \Omega$, $N$ population size, $n_{opt} \in \mathbb{N}^+$.\\
	Let $\mu_k(\cdot), \sigma_k^2(\cdot)$ the posterior mean and variance functions (Algorithm \ref{alg::Basic_BO})\label{line::mu_sigma}\\
	Let $\mathcal{F}(\cdot) = (\mu_k(\cdot), -\sigma_k^2(\cdot))^\top$\label{line::F}\\
	$X^0 =$ \texttt{selection}($\mathcal{F}(\cdot)$, $\Omega$, $X^i$, $N$)\\
	$t = 0$\\
	\While{\textit{a stopping criterion is not satisfied}}{
		$P^t$ = \texttt{getParents}($\mathcal{F}(\cdot)$, $\Omega$, $X^t$) \label{line::parents}\\
		$O^t$ = \texttt{crossover}($\Omega$, $P^t$) \label{line::crossover}\\
		$\tilde{O}^t$ = \texttt{mutation}($\Omega$, $O^t$) \label{line::mutation}\\
		$\hat{X}^{t + 1} = X^t \cup \tilde{O}^t$ \label{line::new_offsprings}\\
		$X^{t + 1} =$ \texttt{selection}($\mathcal{F}(\cdot)$, $\Omega$, $\hat{X}^{t + 1}$, $N$)\\ \label{line::selection}
		\If{$t \mod n_{opt} = 0$} { \label{line::start_opt}
			$\hat{X}^{t + 1} =$ \texttt{optimizePopulation}($\mathcal{F}(\cdot)$, $\Omega$, $X^{t + 1}$) \label{line::optimize_population}\\
			$X^{t + 1} =$ \texttt{selection}($\mathcal{F}(\cdot)$, $\Omega$, $\hat{X}^{t + 1}$, $N$)
		}\label{line::end_opt}
		$t = t + 1$
	}
	$X^\star = X^t$ \label{line::end_1}\\
	\Return $X^\star$ \label{line::end_2}
\end{algorithm}

\texttt{NSMA} is a memetic approach that combines the great exploration capabilities of \texttt{NSGA-II} \cite{deb2002}, which represents the de facto popularity-wise standard for hard box-constrained MOO problems, and the tools typical of MOO gradient-descent methodologies, more efficient and effective than other classes of MOO algorithms on problems with reasonable regularity assumptions. 

Like \texttt{NSGA-II}, \texttt{NSMA} works with a fixed size population of $N$ solutions, starting from an initial one denoted by $X^i$; the latter is then evolved at each iteration through the application of specific random-based operators, such as \texttt{crossover}, \texttt{mutation} and \texttt{selection}. In particular, \texttt{crossover} (Line \ref{line::crossover}) and \texttt{mutation} (Line \ref{line::mutation}) aim to generate new, hopefully better, solutions from selected ones of the current population, which we call \textit{parents} (Line \ref{line::parents}); after inserting the new solutions, called \textit{offsprings}, into the population (Line \ref{line::new_offsprings}), the latter is reduced through the \texttt{selection} operator (Line \ref{line::selection}) in order to have exactly $N$ members. The selection is driven by two metrics, i.e., the \textit{Ranking} and the \textit{Crowding Distance}. While the first metric allows to sort the population based on the dominance relation (see Section \ref{subsec::MOO}), the second one indicates if a point is in a poorly populated area of the objectives space or not. In particular, points with an high crowding distance are extremes of the current Pareto front approximation or they are located far from the other solutions of the population; usually, these points are preferred for starting new search steps to get wider Pareto front approximations. More information about the mentioned operators and metrics can be found in \cite{deb2002}.

The main difference between \texttt{NSMA} and \texttt{NSGA-II} is represented by the \texttt{optimizePopulation} function (Line \ref{line::optimize_population}), which is called every $n_{opt}$ iterations. This function allows to refine the current population $X^{t + 1}$, performing some gradient-based search steps on selected points; the selection is again based on the \textit{Ranking} and the \textit{Crowding Distance} metrics. The search steps are carried out along \textit{steepest partial descent directions} \cite{cocchi2020convergence}: given a point $\bar{x} \in \Omega$, this type of direction can be found solving the following optimization problem
\begin{equation}
	\label{eq::par-dir}
	\theta^{\mathcal{I}}(\bar{x}) = \min_{\substack{d \in \mathcal{D}(\bar{x})\\\|d\|_\infty \le 1}}\max_{j \in \mathcal{I}} \nabla f_j(\bar{x})^\top d,
\end{equation}
where $\mathcal{D}(\bar{x})$ is defined as in Definition \ref{def::par-stat} and $\mathcal{I} \subseteq \{1,2\}$ denotes a subset of objectives indices. We denote as $v^\mathcal{I}(\bar{x})$ the (not necessarily unique) solution of problem \eqref{eq::par-dir} at $\bar{x}$. We trivially deduce that, if $\theta^\mathcal{I}(\bar{x}) = 0$, then $\bar{x}$ is Pareto Stationary (Definition \ref{def::par-stat}) w.r.t. the objective functions indicated by $\mathcal{I}$. On the other hand, if $\bar{x}$ is not Pareto stationary, then we can employ an \textit{Armijo-Type Line Search} to find a step along $v^\mathcal{I}(\bar{x})$, then leading to a new point that can be added in the population. Note that, depending on how $\mathcal{I}$ is defined, this mechanism can refine a point w.r.t.\ both $\mu_k(\cdot)$ and $\sigma_k^2(\cdot)$, only $\mu_k(\cdot)$ or only $\sigma_k^2(\cdot)$, thus contributing to explore the objectives space and to obtain wider Pareto front approximations. By design, the \texttt{optimizePopulation} function drives the selected points towards Pareto stationarity; however, as can be deduced from Section \ref{subsec::MOO}, the generated new solutions could be far from the actual Pareto front of the problem. In such cases, the \texttt{NSGA-II} operators could allow to escape from these non-optimal solutions, so that then the gradient-descent tools could be again effectively exploited from other solutions to refine the current Pareto front approximation.

At the end of the execution, \texttt{NSMA} returns the population obtained during the last iteration (Lines \ref{line::end_1}-\ref{line::end_2}). From this population, we will select $q$ new candidate solutions where to evaluate the objective function $f(\cdot)$ of the original problem \eqref{eq::bo-prob}. In the next section, we report two possible selection criteria.

\subsection{Clustering Approaches for Pareto Front Points Selection}
\label{subsec::clustering}

After reconstructing, as accurately as possible, the Pareto front of the BOO problem \eqref{eq::acq-f-biob} through \texttt{NSMA} (Section \ref{subsec::NSMA-for-BO}), we have to select from it $q$ candidate points, which will be then evaluated on the original problem \eqref{eq::bo-prob}. In this section, we propose two clustering-based selection strategies, whose schemes are reported in Algorithm \ref{alg::clustering_approaches}. 

\begin{algorithm}[!h]
	\caption{Clustering Approaches for Pareto Front Selection} \label{alg::clustering_approaches}
	Input: $\mathcal{F}:\mathbb{R}^n\to \mathbb{R}^2$, $\Omega$, $X^\star = \left\{x_1^\star,\ldots,x_{|X^\star|}^\star\right\} \subset \Omega$ set of solutions generated by \texttt{NSMA} (Algorithm \ref{alg::NSMA}), $q > 1 \in \mathbb{N}$ batch size.\\
	\If{\textit{Clustering in variables space}} {
		Perform \texttt{K-MEANS} on $X^\star$ to find $q$ clusters \label{line::k_means_X}\\
		Let $c_{1:q}$ be the $q$ resulting cluster centers \label{line::centroids_X}\\
		Let $x_{1:q} = c_{1:q}$ \label{line::new_X_X}
	}
	\Else { 
		Let $F^\star = \left\{f_1^\star,\ldots,f_{|X^\star|}^\star\right\} \subset \mathbb{R}^2$ be the image of $X^\star$ through $\mathcal{F}(\cdot)$ \label{line::eval_F}\\
		Perform \texttt{K-MEANS} on $F^\star$ to find $q$ clusters \label{line::k_means_F} \hspace{0.7cm}\textit{(Clustering in objectives space)}\\
		Let $c_{1:q}$ be the $q$ resulting cluster centers\\
		\For{$i \in \{1,\ldots, q\}$} {
			Let $z_i = \argmin\limits_{z \in\{1,\ldots, |X^\star|\}} \|f_z^\star - c_i\|$ \label{line::j_F}\\
			Let $x_i = x_{z_i}^\star$ \label{line::new_X_F}
		}\label{line::end_for_F}
	} 
	\Return $x_{1:q}$
\end{algorithm}

In the first clustering strategy (Lines \ref{line::k_means_X}-\ref{line::new_X_X}), we perform the \texttt{K-MEANS} algorithm \cite{hartigan79, lloyd82} to find $q$ clusters in the variables space (Line \ref{line::k_means_X}). \texttt{K-MEANS} is arguably the most popular local search algorithm for clustering problems; the objective consists in the minimization of the sum-of-squares of the Euclidean distances of the samples to their cluster centers \cite{Hansen1997}. Besides requiring a few parameters to set (number of clusters to find included), it proves to be particularly efficient and robust in many scenarios. The $q$ resulting cluster centers are taken as the new candidates for the original problem \eqref{eq::bo-prob} (Lines \ref{line::centroids_X}-\ref{line::new_X_X}).

On the other hand, the second clustering approach (Lines \ref{line::eval_F}-\ref{line::end_for_F}) involves a \texttt{K-MEANS} execution in the objectives space. For each cluster center $c_i$, with $i \in \{1,\ldots, q\}$, we then find the closest Pareto front point $f^\star_{z_i}$ (Line \ref{line::j_F}); the corresponding $x^\star_{z_i}$ will be one of the $q$ candidate points to be evaluated on \eqref{eq::bo-prob} (Line \ref{line::new_X_F}). 

Although clustering in the objectives space seems to be more intuitive than clustering in the variables space, the latter could be more effective in high-dimensional scenarios ($n > 20$), where a blind exploration of the feasible set boundaries is usually carried out employing the standard acquisition functions (Section \ref{subsubsec::acq-f}). This behavior on high-dimensional problems is well-known in the literature and is called the \textit{boundary issue} \cite{swersky2017improving, oh18}. In fact, in high dimensional spaces, it happens very frequently that along the boundaries of the feasible region many points are characterized by high $\sigma_k^2(\cdot)$ values (uncertainty). Thus, most methods based on a compromise between the two objectives, such as \texttt{q-EI} and \texttt{q-LCB}, tend to generate candidate points on the border of the feasible region. This fact can be considered as one of the main reasons Bayes-Opt methods have been applied in the past only to very low--dimensional problems, where this effect is less relevant. In this context, clustering in the objectives space does not represent a possible solution, since the considered BOO problem \eqref{eq::acq-f-biob} involves the same functions as \texttt{q-EI} and \texttt{q-LCB}. On the other hand, clustering in the variables spaces could be an \textit{heuristic} solution to the boundary issue: clustering points in the boundaries, far from each other, could end up with setting some cluster centers in the interior of the feasible set, thus allowing to collect information about $f(\cdot)$ in more internal regions; at the next iteration, this information would lead to more informative posterior mean and variance functions and, then, to a better exploitation of the Pareto front reconstruction. 

\begin{figure*}[!h]
	\centering
	\subfloat[Rastrigin, $n=2$]{\includegraphics[width=0.35\textwidth]{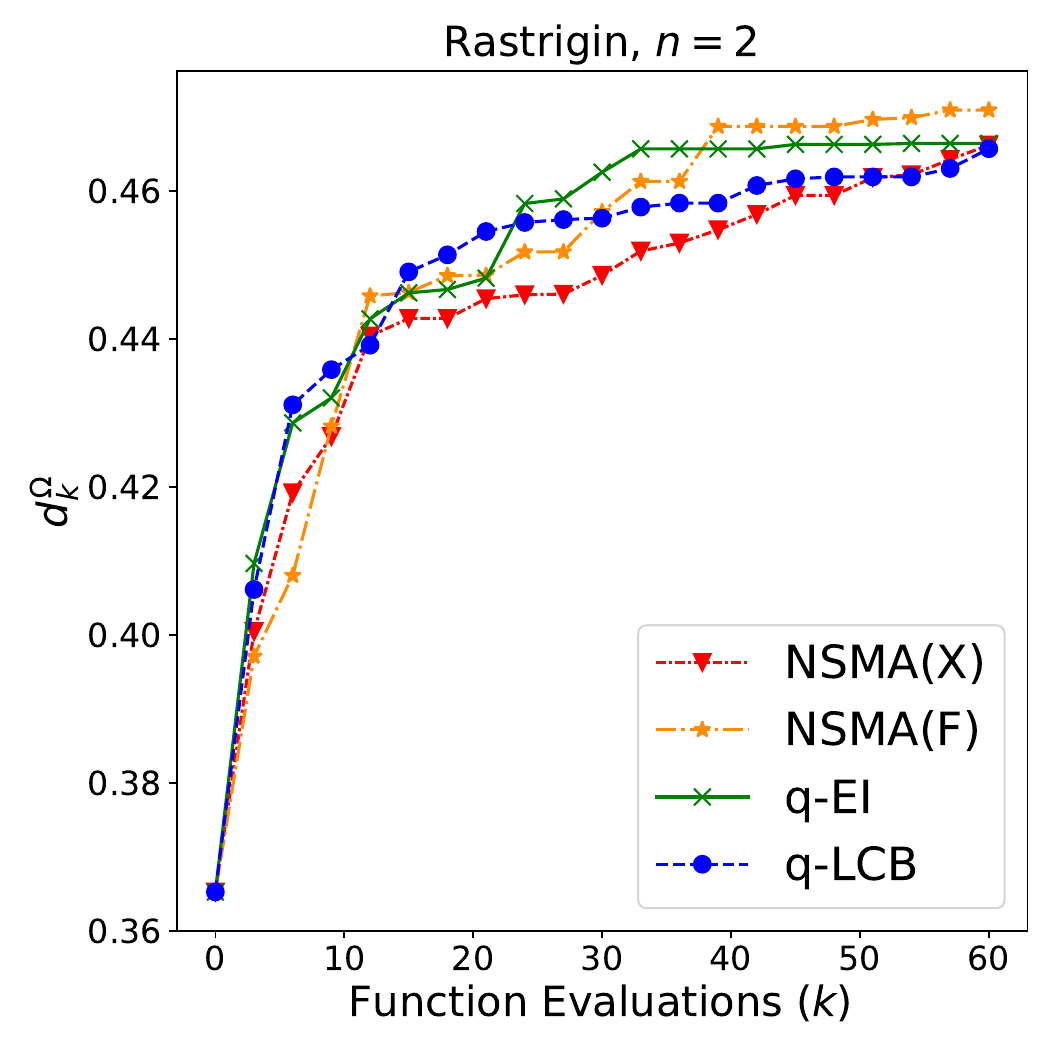}}
	\hfil
	\subfloat[Rastrigin, $n=10$]{\includegraphics[width=0.35\textwidth]{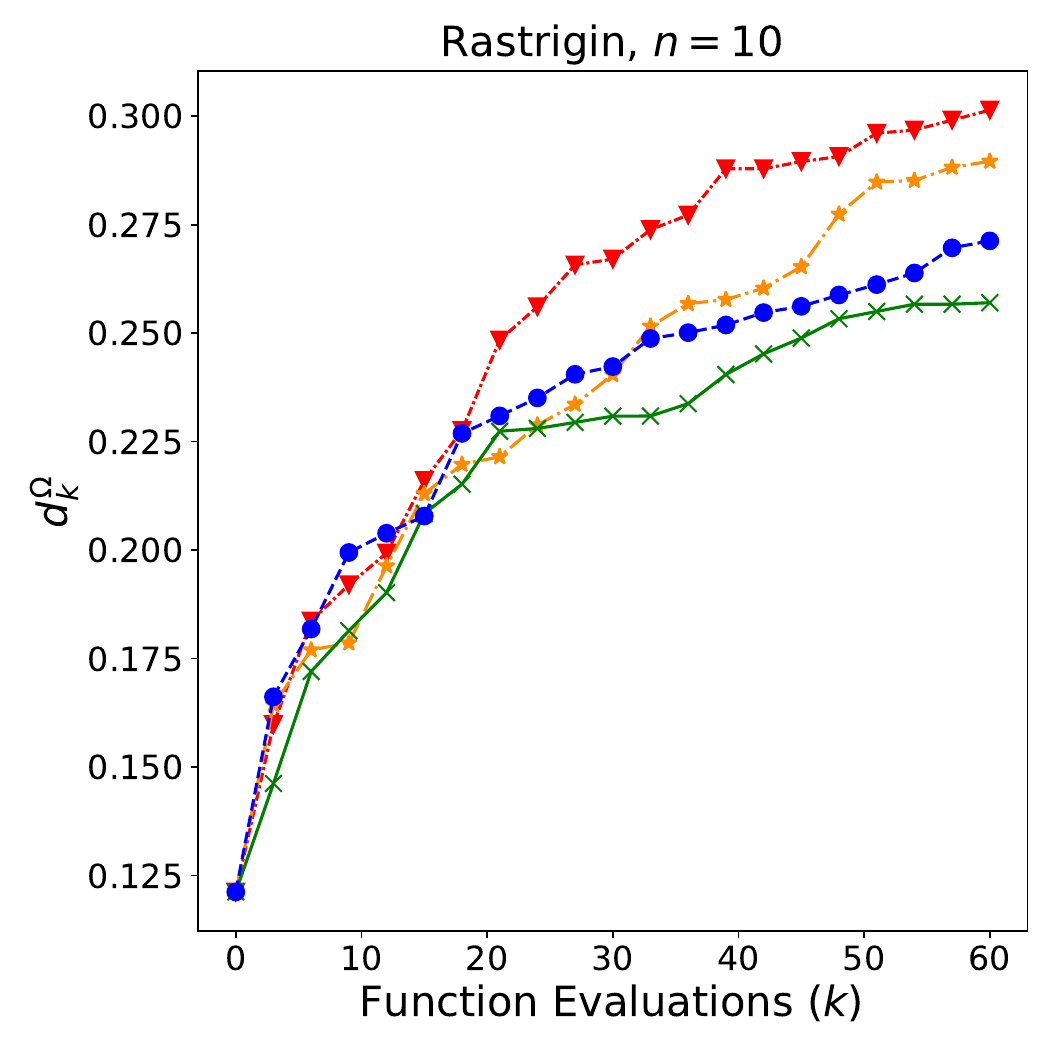}}
	\\
	\subfloat[Rastrigin, $n=50$]{\includegraphics[width=0.35\textwidth]{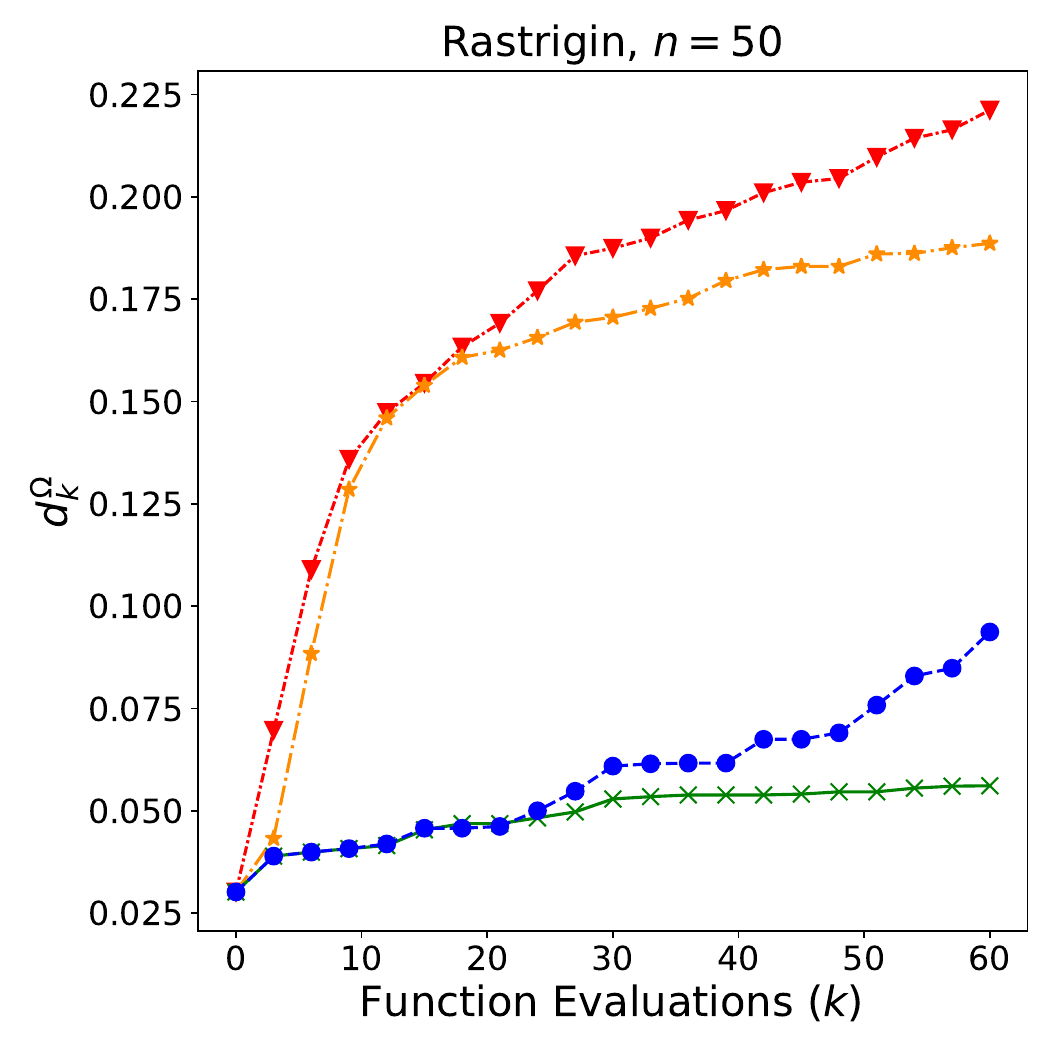}}
	\hfil
	\subfloat[Rastrigin, $n=100$]{\includegraphics[width=0.35\textwidth]{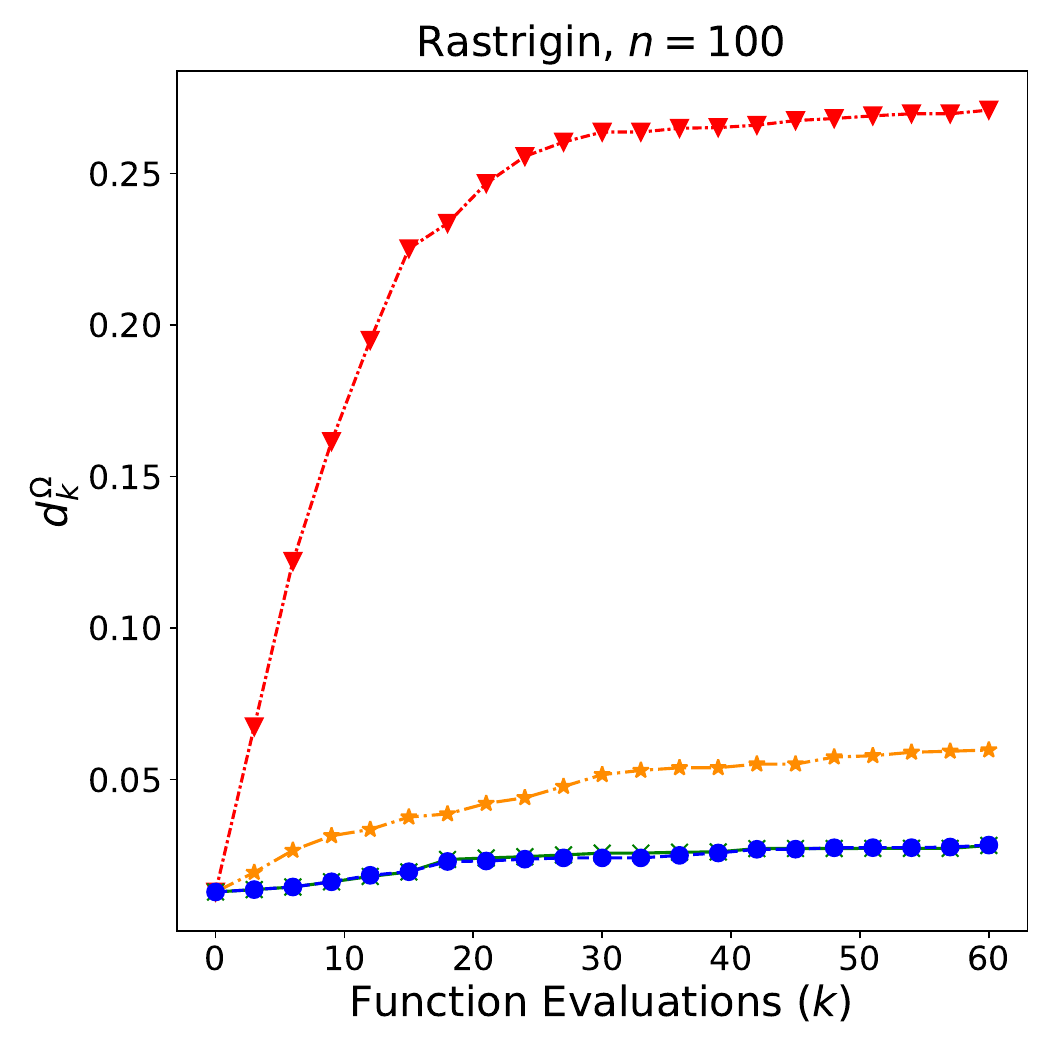}}
	\caption{Plot of $d^\Omega_k$ \eqref{eq::dist-from-bounds} w.r.t.\ the number of function evaluations $k$. The tests, whose setting can be found in Section \ref{subsec::experimental_settings}, were performed on the Rastrigin function \cite{torn1989global} at different dimensionalities.}
	\label{fig:boundary_issues}
\end{figure*}

In Figure \ref{fig:boundary_issues}, we report an empirical example supporting our arguments. Here, we considered four acquisition strategies: \texttt{NSMA(X)}, \texttt{NSMA(F)}, \texttt{q-EI}, \texttt{q-LCB} (for more details on the last two strategies, see Section \ref{subsubsec::acq-f}); the acronyms \texttt{(X)} and \texttt{(F)} indicates in which space \texttt{K-MEANS} was performed. For each observed point $x_k$ in the Bayes-Opt procedure (Algorithm \ref{alg::Basic_BO}), we calculated the component-wise minimum distance from the problem boundaries; this distance was compared with the one of the previous point $x_{k - 1}$ and the maximum was taken. The resulting metric is formalized in the following equation:
\begin{equation}
	\label{eq::dist-from-bounds}
	d^\Omega_k = 
	\begin{cases} 
		0, & \mbox{if }k = 0 \\ 
		\max\left\{d^\Omega_{k - 1}, \min\limits_{i \in \{1,\ldots, n\}}\min\left\{x^i_k - l^i, u^i - x^i_k\right\}\right\}, & \mbox{otherwise}
	\end{cases}.
\end{equation}

We can observe that, in the low-dimensional settings ($n \le 10$), all the acquisition strategies allowed to explore points distant from the boundaries from the very first iterations. The situation is different in the high-dimensional scenarios, where the effects of the boundary issue can be noticed on the performance of \texttt{q-EI} and \texttt{q-LCB}. Unlike the two mentioned approaches, \texttt{NSMA(F)} managed to escape from the boundaries in the $n=50$ case; however, in the hardest scenario ($n=100$), even this method struggled. The only methodology capable of finding points far from the boundaries regardless the problem dimensionality was \texttt{NSMA(X)}, with remarkable performance w.r.t.\ the competitors on the $n=100$ case. 

\section{Computational Experiments}
\label{sec::computational-experiments}

In this section, we report the results of thorough computational experiments, comparing the proposed acquisition strategy with state-of-the-art ones from the batch Bayes-Opt literature. The experimental code was written in \texttt{Python3}\footnote{The implementation code can be found at \href{https://github.com/FranciC19/biobj_acquistion_function_for_BO}{github.com/FranciC19/biobj\_acquistion\_function\_for\_BO}.} and was run on a computer with the following characteristics: Ubuntu 22.04 OS, Intel(R) Core(TM) i5-10600KF 6 cores 4.10GHz, 32 GB RAM. For the implementation of the Bayes-Opt procedure (Algorithm \ref{alg::Basic_BO}), including the Gaussian Process regression and the optimization of the standard acquisition functions, we relied on the \texttt{BoTorch}\footnote{The \texttt{BoTorch} code can be found at \href{https://github.com/pytorch/botorch}{github.com/pytorch/botorch}.} library \cite{Botorch}. Finally, for the \texttt{NSMA} implementation we used the code\footnote{The \texttt{NSMA} code can be found at \href{https://github.com/pierlumanzu/nsma}{github.com/pierlumanzu/nsma}.} distributed by its authors; in order to solve instances of problem \eqref{eq::par-dir}, we employed the Gurobi Optimizer (Version 10) \cite{gurobi}.

\subsection{Experimental Settings}
\label{subsec::experimental_settings}

In this section, we report the tested acquisition strategies, the functions used for the comparisons and, finally, the employed metrics.

\subsubsection{Acquisition Strategies}

Throughout the section, we indicate with \texttt{NSMA} the BOO based acquisition strategy proposed in Section \ref{sec::bi-obj-acq-f-for-BO} of this manuscript. We decided also to test our methodology with a different solver, that is, \texttt{NSGA-II} \cite{deb2002}. As mentioned in Section \ref{sec::introduction}, this algorithm has been already tested to solve a BOO problem similar to \eqref{eq::acq-f-biob}. Moreover, \texttt{NSGA-II} is arguably one of the most representative evolutionary algorithms for hard MOO problems. Thus, in this context, it represents the most direct competitor. For both the algorithms, we set the population size $N = 100$ and a maximum number of iterations equal to $20$; for \texttt{NSMA}, the parameter $n_{opt} = 5$; the other parameters values were set as recommended in \cite{NSMA, deb2002}.

Both algorithms were executed with the two clustering approaches proposed in Section \ref{subsec::clustering} for the Pareto front points selection. The acronyms \texttt{(X)} and \texttt{(F)} indicates if \texttt{K-MEANS} was performed in the variables space or in the objectives one.

Finally, we considered as competitors the acquisition strategies involving the optimization of \texttt{q-EI} and \texttt{q-LCB}, respectively. A brief description of these acquisition functions can be found in Section \ref{subsubsec::acq-f}. For the integrals estimation, we employed the Sobol Quasi Monte-Carlo Normal Sampler with $512$ samples. Moreover, as suggested in \cite{wilson2017reparameterization}, for \texttt{q-LCB} we set the parameter $\beta = \sqrt{3}$. Finally, we set a maximum number of iterations for the solver \texttt{L-BFGS-B} \cite{zhu97} equal to $100$. 

At each iteration of the Bayes-Opt procedure, the acquisition strategies were run with $100$ initial points, uniformly chosen in the feasible set $\Omega$.

Regarding the Gaussian Process, we used all the default settings and the recommendations proposed in the \texttt{BoTorch} library. Both in the Gaussian Process regression and in the evaluations of the posterior mean and variance functions, the solutions were first scaled in $[0, 1]$ considering the problem bounds; moreover, for numerical stability, the observed $f(\cdot)$ values were min-max scaled before the Gaussian Process regression.

\subsubsection{Benchmark Functions}
\label{subsubsec::benchmark_functions}

The functions\footnote{The formulations, as well as other information, of some of the tested functions can be found at \href{https://www.sfu.ca/~ssurjano/index.html}{sfu.ca/$\sim$ssurjano/index.html}.} employed for the comparisons are listed in Table \ref{tab::test_functions}. As mentioned in Sections \ref{sec::introduction}-\ref{subsec::overview_bo}, all the evaluations were supposed to be noise-free. Note that some of the tested functions are scalable w.r.t.\ the dimensionality $n$; this feature was crucial in our experimentation to compare the acquisition strategies in (hard) high-dimensional scenarios, where additional difficulties arise (e.g., the boundary issue; for more details on the topic, we refer the reader to Section \ref{subsec::clustering}).

\begin{table*}[!h]
	\caption{Benchmark functions used for the experiments. For each of them, we indicate: the source paper; the name; the tested values for the dimensionality $n$; the optimal value; if it has/has not the zero solution as a global minimum; the lower and upper bounds.}
	\label{tab::test_functions}
	\centering
	\footnotesize
	\begin{tabular}{cc||cccc}
		\midrule
		\multirow{2}{*}{\textbf{Source}} & \multirow{2}{*}{\textbf{Name}} & \multirow{2}{*}{$\mathbf{n}$} & \multirow{2}{*}{$\mathbf{f^\star}$} & \multirow{2}{*}{$\mathbf{x^\star = 0_n}$} & \multirow{2}{*}{\textbf{Bounds}} \\
		&&&&&\\
		\midrule
		\midrule
		\multirow{2}{*}{\cite{torn1989global}} & \multirow{2}{*}{Rastrigin} &  2, 4, 10, & \multirow{2}{*}{$0.0$} & \multirow{2}{*}{\cmark} & \multirow{2}{*}{$[-5.12, 5.12]^n$} \\ 
		& &  20, 50, 100 & & & \\ 
		\midrule
		\multirow{3}{*}{\cite{hartman}\vspace{0.2cm}} & \multirow{3}{*}{Hartmann\vspace{0.2cm}} & 3 & $-3.86278$ & \xmark & $[0.0, 1.0]^n$ \\
		\cmidrule{3-6}
		& & 6 & $-3.32237$ & \xmark & $[0.0, 1.0]^n$ \\
		\midrule
		\multirow{2}{*}{\cite{rosenbrock1960automatic}} & \multirow{2}{*}{Rosenbrock} & 2, 4, 10, & \multirow{2}{*}{$0.0$} & \multirow{2}{*}{\xmark} & \multirow{2}{*}{$[-5, 10]^n$} \\
		&&20, 50, 100&&&\\
		\midrule
		\multirow{2}{*}{\cite{back1993overview}} & \multirow{2}{*}{Ackley\_1} & 2, 4, 10, & \multirow{2}{*}{$0.0$} & \multirow{2}{*}{\cmark} & \multirow{2}{*}{$[-32.768, 32.768]^n$} \\
		&&20, 50, 100&&&\\
		\midrule
		\multirow{2}{*}{\cite{RAHNAMAYAN20071605}} & \multirow{2}{*}{Alpine\_1} & 2, 4, 10, & \multirow{2}{*}{$0.0$} & \multirow{2}{*}{\cmark} & \multirow{2}{*}{$[-10, 10]^n$} \\
		&&20, 50, 100&&&\\
		\midrule
		\multirow{2}{*}{\cite{branin1972widely}} & \multirow{2}{*}{Branin} & \multirow{2}{*}{2} & \multirow{2}{*}{$0.397887$} & \multirow{2}{*}{\xmark} & $x_1 \in [-5, 10]$ \\
		&&&&&$x_2 \in [0, 15]$ \\
		\midrule
		\multirow{2}{*}{\cite{schwefel1981numerical}} & \multirow{2}{*}{Schwefel} & 2, 4, 10, & \multirow{2}{*}{$0.0$} & \multirow{2}{*}{\xmark} & \multirow{2}{*}{$[-500, 500]^n$} \\
		&&20, 50, 100&&&\\
		\midrule
		\multirow{2}{*}{\cite{levy}} & \multirow{2}{*}{Levy\_8} & 2, 4, 10, & \multirow{2}{*}{$0.0$} & \multirow{2}{*}{\xmark} & \multirow{2}{*}{$[-10, 10]^n$} \\
		&&20, 50, 100&&&\\
		\midrule
		\multirow{4}{*}{\cite{molga2005test}} & \multirow{4}{*}{Michalewicz} & 2 & $-1.80130341$ & \xmark & $[0,\pi]^n$ \\
		\cmidrule{3-6}
		& & 5 & $-4.687658$ & \xmark & $[0,\pi]^n$ \\
		\cmidrule{3-6}
		& & 10 & $-9.66015$ & \xmark & $[0,\pi]^n$ \\
		\midrule
		\multirow{2}{*}{\cite{branin1972widely}} & \multirow{2}{*}{SixHumpCamel} & \multirow{2}{*}{2} & \multirow{2}{*}{$-1.0316$} & \multirow{2}{*}{\xmark} & $x_1 \in [-3, 3]$\\
		&&&&&$x_2 \in [-2, 2]$ \\
		\midrule
		\multirow{2}{*}{\cite{bukin1997new}} & \multirow{2}{*}{Bukin\_6} & \multirow{2}{*}{2} & \multirow{2}{*}{$0.0$} & \multirow{2}{*}{\xmark} & $x_1 \in [-15, -5]$\\
		&&&&&$x_2 \in [-3, 3]$ \\
		\midrule
		\multirow{2}{*}{\cite{styblinski1990experiments}} & \multirow{2}{*}{StyblinskiTang} & 2, 4, 10, & \multirow{2}{*}{$-39.166166n$} & \multirow{2}{*}{\xmark} & \multirow{2}{*}{$[-5, 5]^n$} \\
		&&20, 50, 100&&&\\
		\midrule
		\multirow{3}{*}{\cite{benchmark_survey}\vspace{0.2cm}} & HolderTable\_2 & 2 & $ -19.2085$ & \xmark & $[-10, 10]^n$ \\
		\cmidrule{2-6}
		& EggHolder & 2 & $-959.6407$ & \xmark & $[-512, 512]^n$ \\
		\midrule
		\cite{opacic1973heuristic} & Shekel ($m=10$) & 4 & $-10.5363$ & \xmark & $[0, 10]^n$ \\
		\midrule
	\end{tabular}
\end{table*}

For each function, we started the Bayes-Opt procedure with $N_i = 10$ points randomly chosen in $\Omega$; the same initial solutions were considered for all the acquisition strategies. We mainly tested the methodologies with the batch size parameter $q = 3$: ideally, at each iteration of the Bayes-Opt procedure, we would like to obtain as candidates for the original problem \eqref{eq::bo-prob} two points with good values for $\mu_k(\cdot)$ (exploitation) and for $\sigma_k^2(\cdot)$ (exploration), respectively, and one middle trade-off solution. For the sake of completeness, we also experimented other values for $q$; the results are commented at the end of the experimental section. In the Bayes-Opt procedure, we set the number of total objective function evaluations $N_M = 60$; thus, considering that $q=3$, for each function we generated $20$ point batches.

For each acquisition strategy and function, we run the execution $20$ times, each of which characterized by a different seed for the pseudo-random number generator, so as to reduce the experiments sensibility to the random operations. The results were then obtained by calculating the averages of the metrics values achieved in the twenty tests.

\subsubsection{Metrics}
\label{subsubsec::metrics}

We denote by $f^{best}_k$ the best observed objective function value after $k$ evaluations; moreover, we indicate with $f^{best}_L$ the best $f(\cdot)$ value found at the end of the Bayes-Opt procedure. According to Algorithm \ref{alg::Basic_BO}, the two metrics can be formalized in the following way: $f^{best}_k = \min_{i \in \{1,\ldots,k\}}f(x_i)$ and $f^{best}_L = f(x^\star)$. 

In order to observe the performance trend of the methodologies between one Bayes-Opt procedure iteration and another, the values obtained in terms of $f^{best}_k$ are plotted w.r.t. the number of objective function evaluations $k$. In the plots, vertical bars are also shown, representing the 95\% confidence interval \cite{Botorch} over the $f^{best}_k$ values obtained in the twenty runs. 

Another metric we employed is a normalized version of the \textit{regret} \cite{Garnett2023}: $$NR_k = \frac{f^{best}_k - f^\star}{f^{best}_0 - f^\star},$$ where $f^{best}_0$ indicates the best function value obtained by the $10$ initial points of the execution. We trivially deduce that: for all $k$, $NR_k \ge 0$; lower values for this metric mean better performance. In our experimentation, the $NR_k$ values obtained for all $k$ were summarized by calculating the \textit{AUC} (Area Under Curve); the resulting metric is called \textit{NR-AUC}; like $NR_k$, small values of \textit{NR-AUC} are associated with good performance. In brief, we would like to minimize the regret at the end of the optimization, but at the same time we would also like to converge in a few iterations to the final solution: this is particularly relevant when dealing with expensive-to-evaluate objective functions.

Lastly, we employed the variant introduced in \cite{CAPPANERA2024121953} of the performance profiles \cite{dolan2002benchmarking}. Given a solver and a certain metric, the new tool allows to plot the (cumulative) distribution function of the relative gap between the score achieved by that solver and the best score among those obtained by all the considered solvers. In other words, it shows the probability that the relative gap for a solver in a problem is within a factor $\tau \in \mathbb{R}$. Like the performance profiles, the tool proposed in \cite{CAPPANERA2024121953} allows to appreciate the relative performance and robustness of the considered algorithms. Moreover, distribution of relative gap is particularly useful when evaluating metrics involving objective function values, such as $f^{best}_k$ and \textit{NR-AUC}. For a more technical explanation of the tool, we refer the reader to \cite{CAPPANERA2024121953}.

\subsection{Preliminary Comparison of the Bi-Objective Optimization Based Acquisition Strategy with \texttt{q-LCB}}
\label{subsec::preliminary_assessment}

\begin{figure*}[!h]
	\centering
	\subfloat[Holdertable\_2, $n = 2$]{\includegraphics[width=0.32\textwidth]{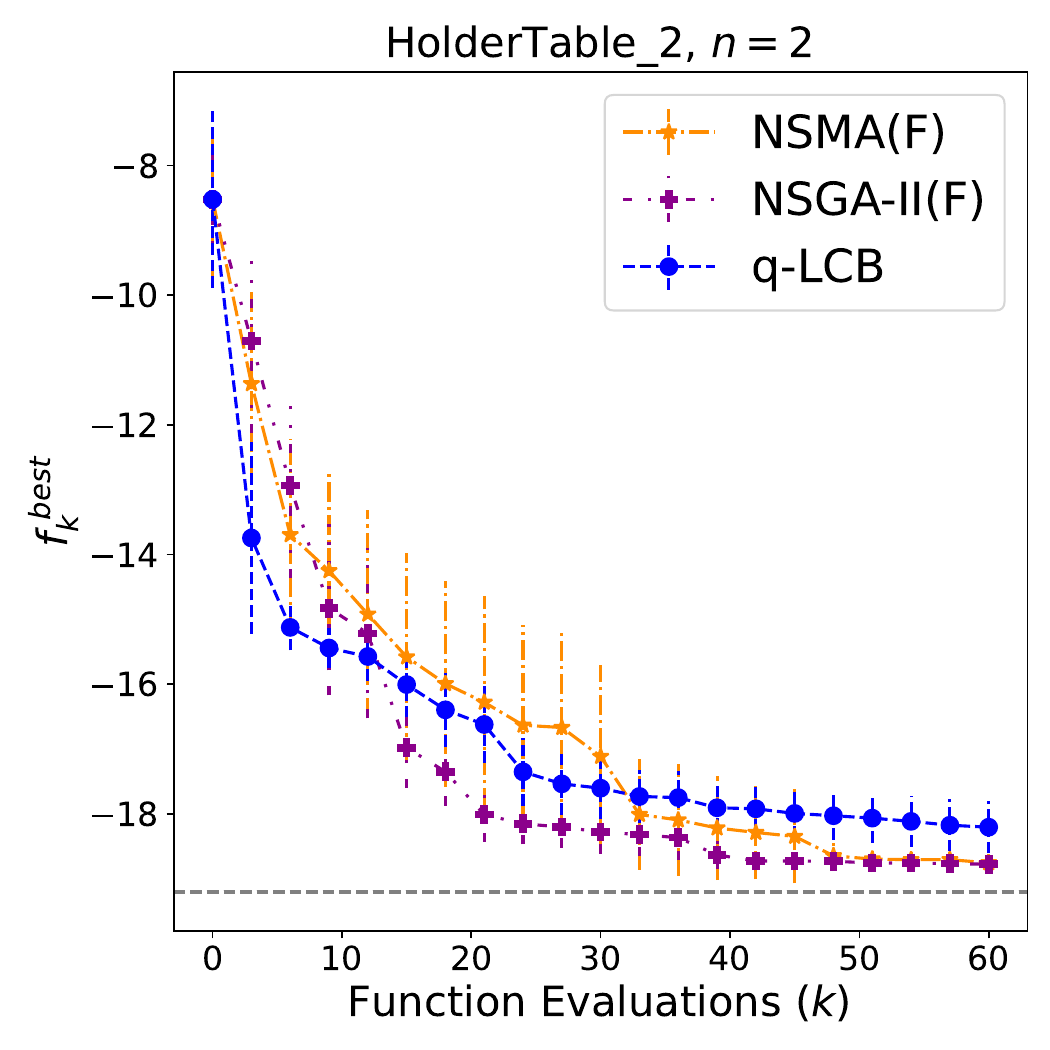}}
	\hfil
	\subfloat[Schwefel, $n = 2$]{\includegraphics[width=0.32\textwidth]{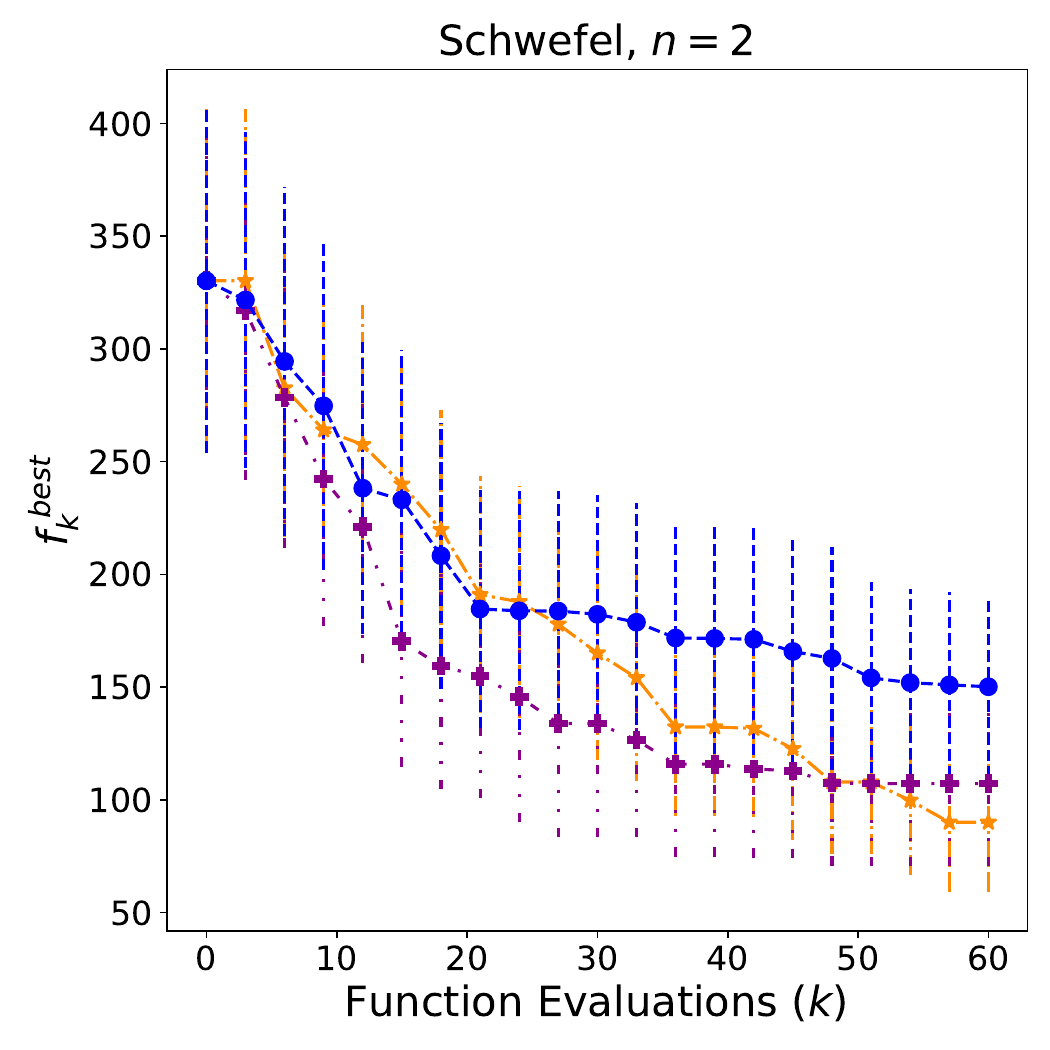}}
	\hfil
	\subfloat[Hartmann, $n = 6$]{\includegraphics[width=0.32\textwidth]{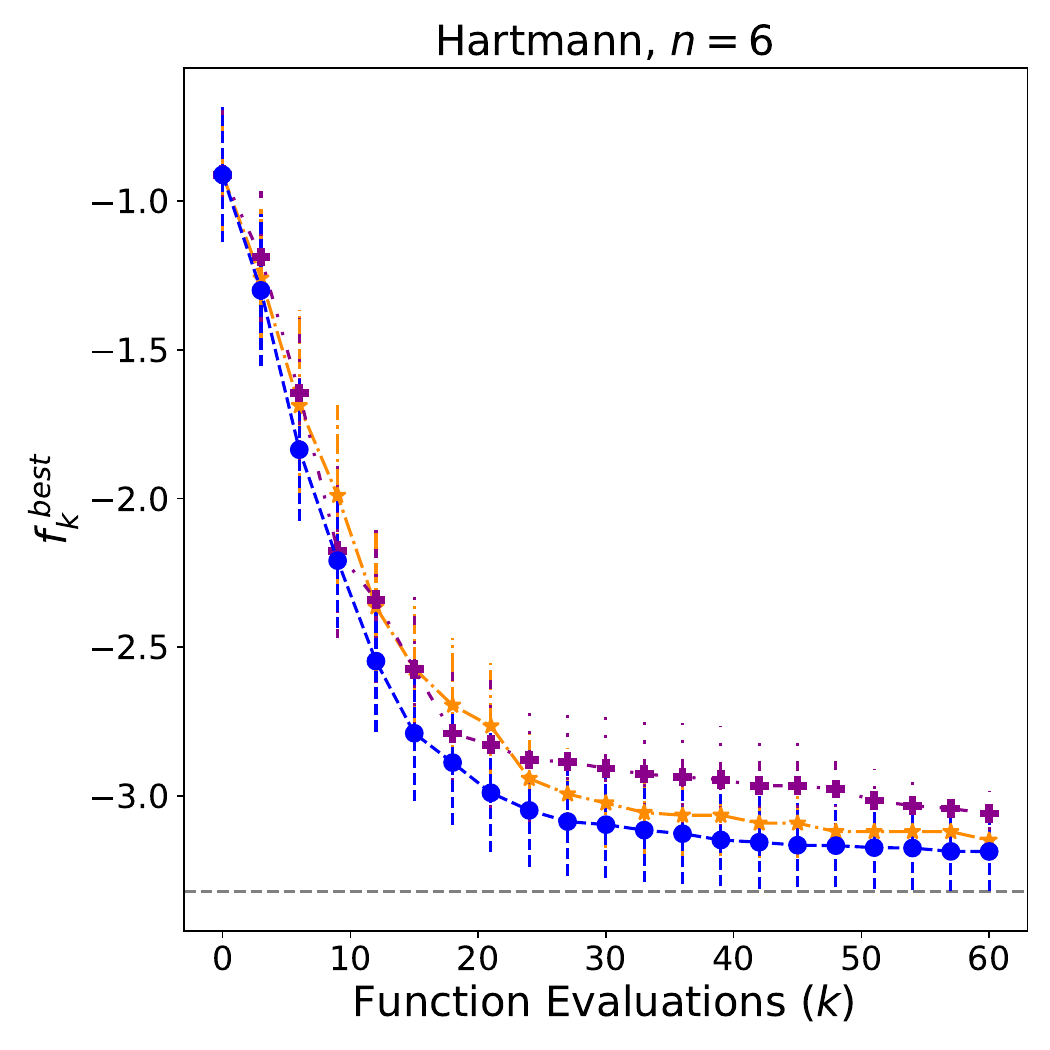}}
	\\
	\subfloat[Ackley\_1, $n = 20$]{\includegraphics[width=0.32\textwidth]{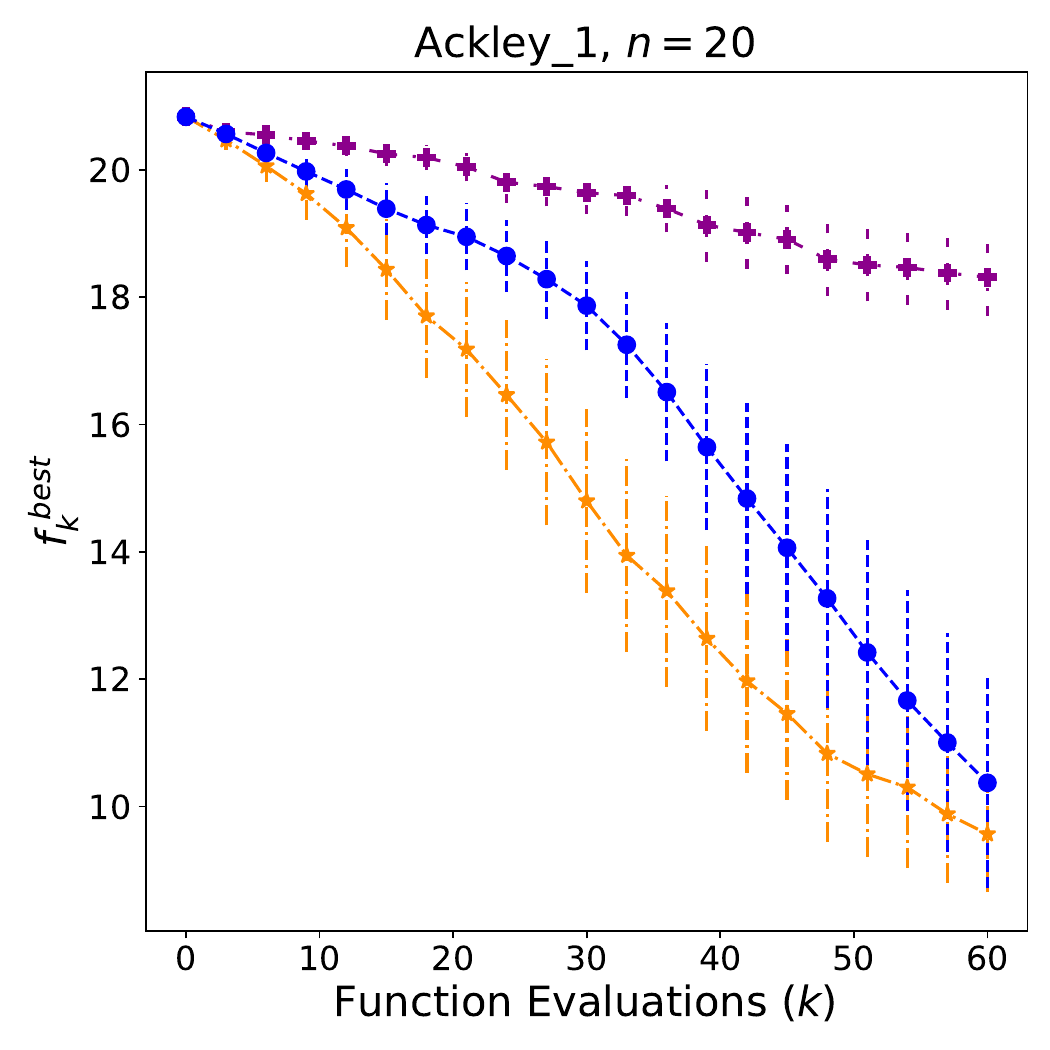}}
	\hfil
	\subfloat[Levy\_8, $n = 20$]{\includegraphics[width=0.32\textwidth]{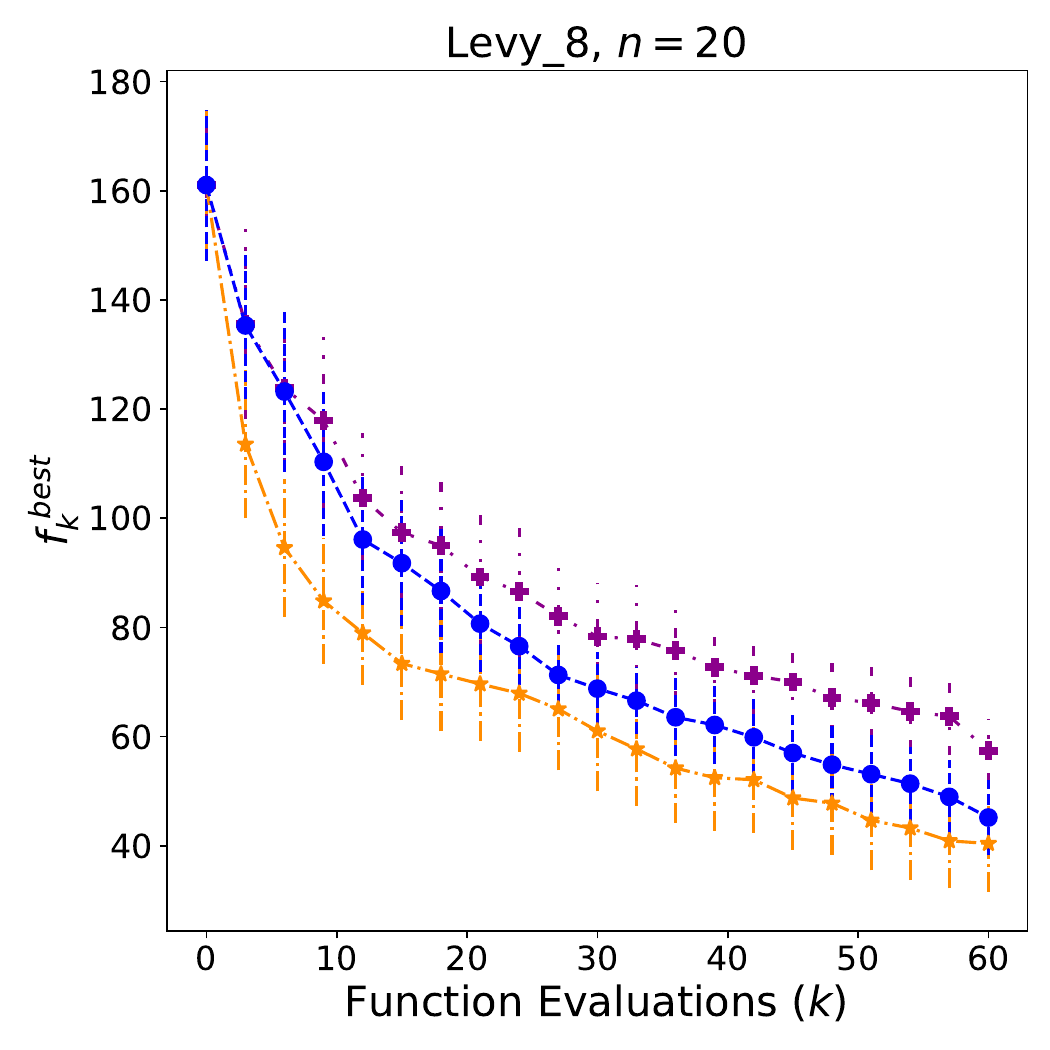}}
	\hfil
	\subfloat[Rastrigin, $n = 50$]{\includegraphics[width=0.32\textwidth]{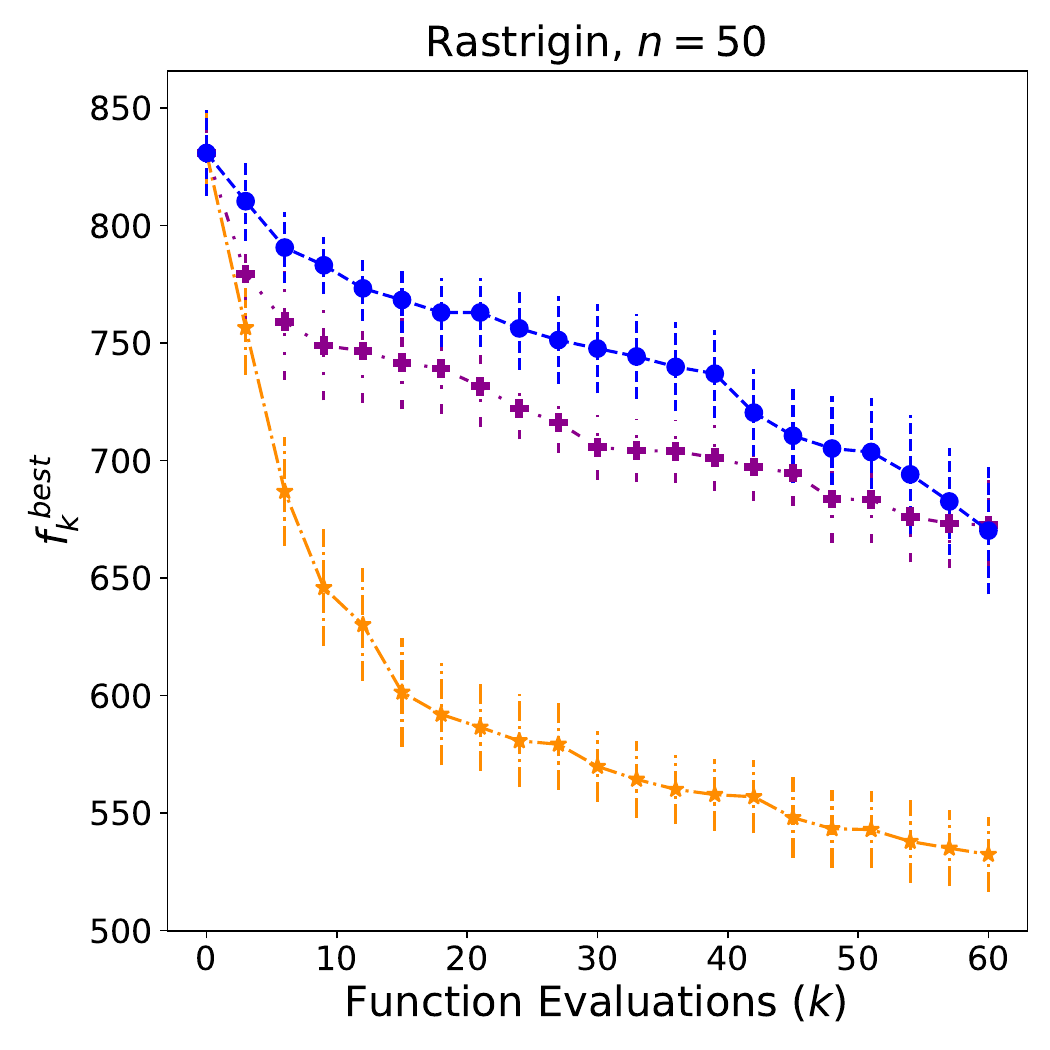}}
	\caption{Plots of the $f^{best}_k$ metric w.r.t.\ the number of function evaluations $k$ for \texttt{NSMA(F)}, \texttt{NSGA-II(F)} and \texttt{q-LCB} on a set of selected functions (see Table \ref{tab::test_functions}). When $f^\star \ne 0$, it is represented by a gray dashed line.}
	\label{fig::bo_biobj}
\end{figure*}

In the first section of computational experiments, we compare the performance of our BOO based acquisition strategy and \texttt{q-LCB} on some selected functions. Note that \texttt{q-LCB} \eqref{eq::qlcb} can be seen as a scalarization of the following BOO problem: $\min_{x \in \Omega} (\mu_k(x), -\sigma_k(x))^\top$, which is similar to the one we use in our strategy; thus, making such a comparison can be useful to appreciate the strengths of reconstructing the Pareto front to choose the next candidate solutions for problem \eqref{eq::bo-prob}. 
In Figure \ref{fig::bo_biobj}, we report the $f^{best}_k$ metric values achieved by \texttt{NSMA(F)}, \texttt{NSGA-II(F)} and \texttt{q-LCB} w.r.t. the number of function evaluations $k$.

On low-dimensional scenarios ($n \le 6$), the BOO based acquisition strategy turned out to be effective in 2 out of 3 problems; from the very first iterations \texttt{NSGA-II(F)} managed to find better solutions than its competitors; moreover, both \texttt{NSMA(F)} and \texttt{NSGA-II(F)} allowed to find better solutions than \texttt{q-LCB} at the end of the Bayes-Opt procedure. In the third problem, i.e., Hartmann, \texttt{NSMA(F)} and \texttt{NSGA-II(F)} were outperformed by \texttt{q-LCB}, but their performance still remained competitive.

When $n \ge 20$, \texttt{NSMA(F)} is the clear winner. In particular, two results of the memetic algorithm have to be highlighted. The first one is the increasing performance gap between \texttt{NSMA(F)} and \texttt{q-LCB} as the value for $n$ increases. This fact highlights an issue typical of scalarization-based methodologies: the choice of the weights (in this context, $\beta$; see Equation \ref{eq::qlcb}) could not be trivial, especially when high-dimensional problems are taken into account. The second result is that \texttt{NSMA(F)} clearly outperformed \texttt{NSGA-II(F)}. As mentioned in \cite[Section 4]{NSMA}, for large values for $n$, \texttt{NSMA} exploited the combined use of the \texttt{NSGA-II} exploration capabilities and the gradient-descent tools to get better and wider Pareto front approximations than the genetic approach. This behavior also occurred in the considered Bayes-Opt context, allowing \texttt{NSMA(F)} to propose better solutions for the objective function $f(\cdot)$ evaluation.

\begin{figure}[!h]
	\centering
	\subfloat[\textit{NR-AUC} -- Rastrigin]{\includegraphics[width=0.32\textwidth]{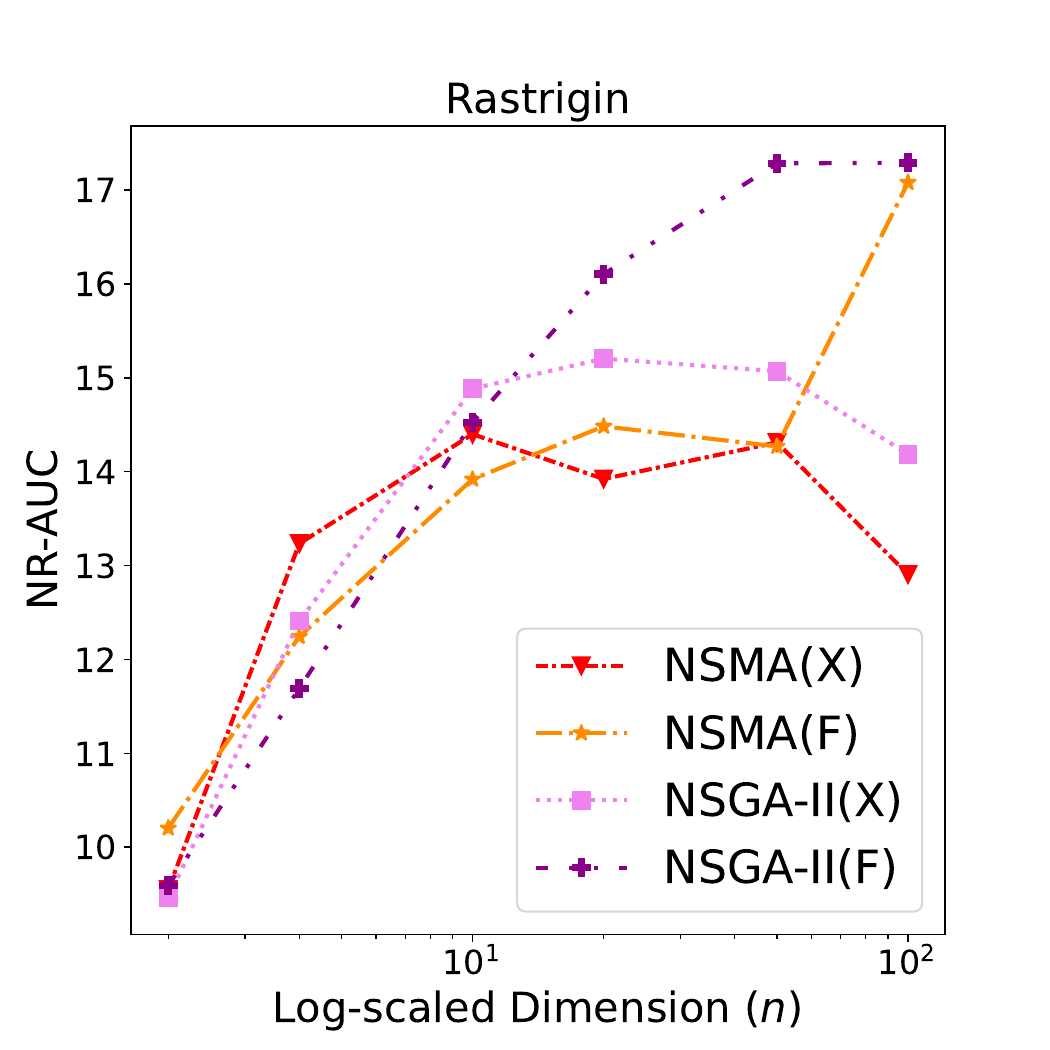}}
	\hfil
	\subfloat[\textit{NR-AUC} -- Ackley\_1]{\includegraphics[width=0.32\textwidth]{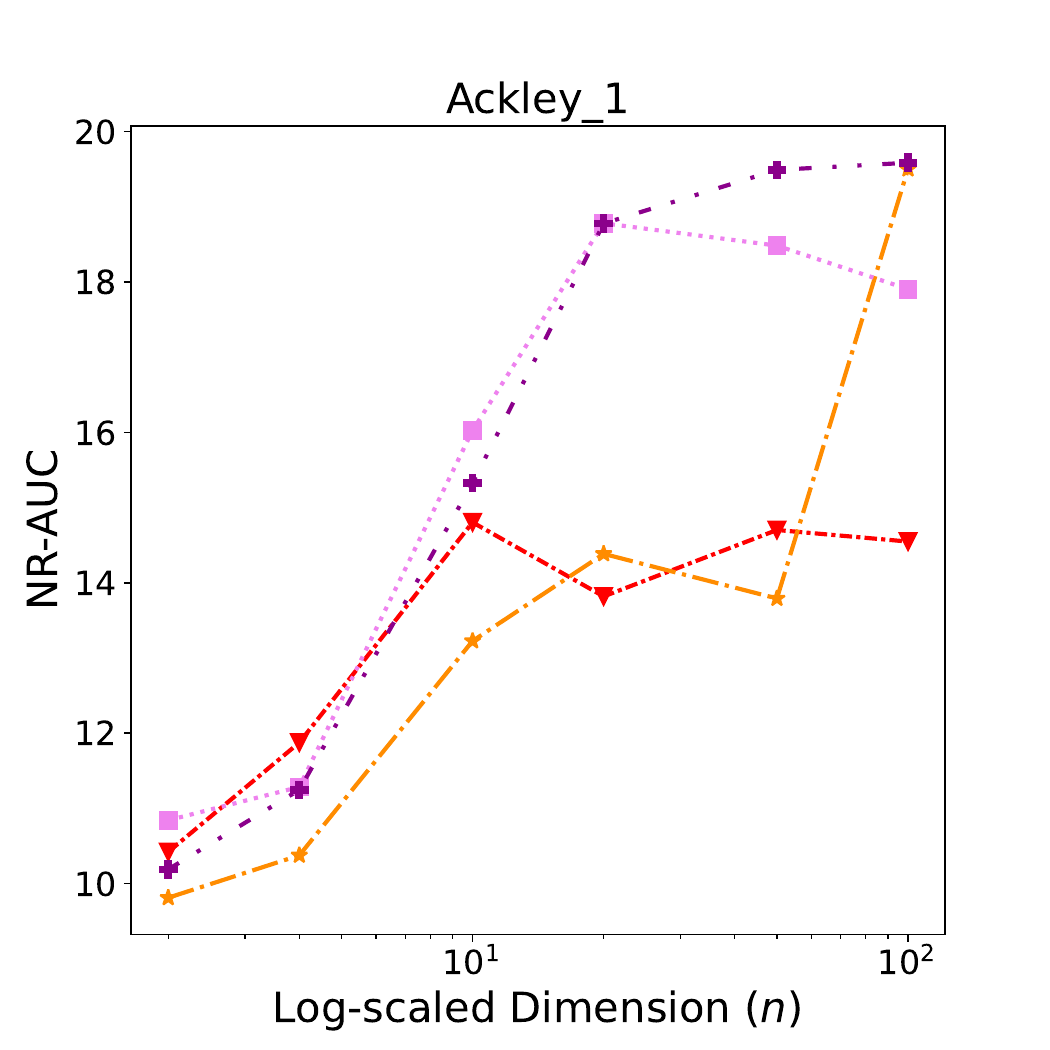}}
	\hfil
	\subfloat[\textit{NR-AUC} -- Levy\_8]{\includegraphics[width=0.32\textwidth]{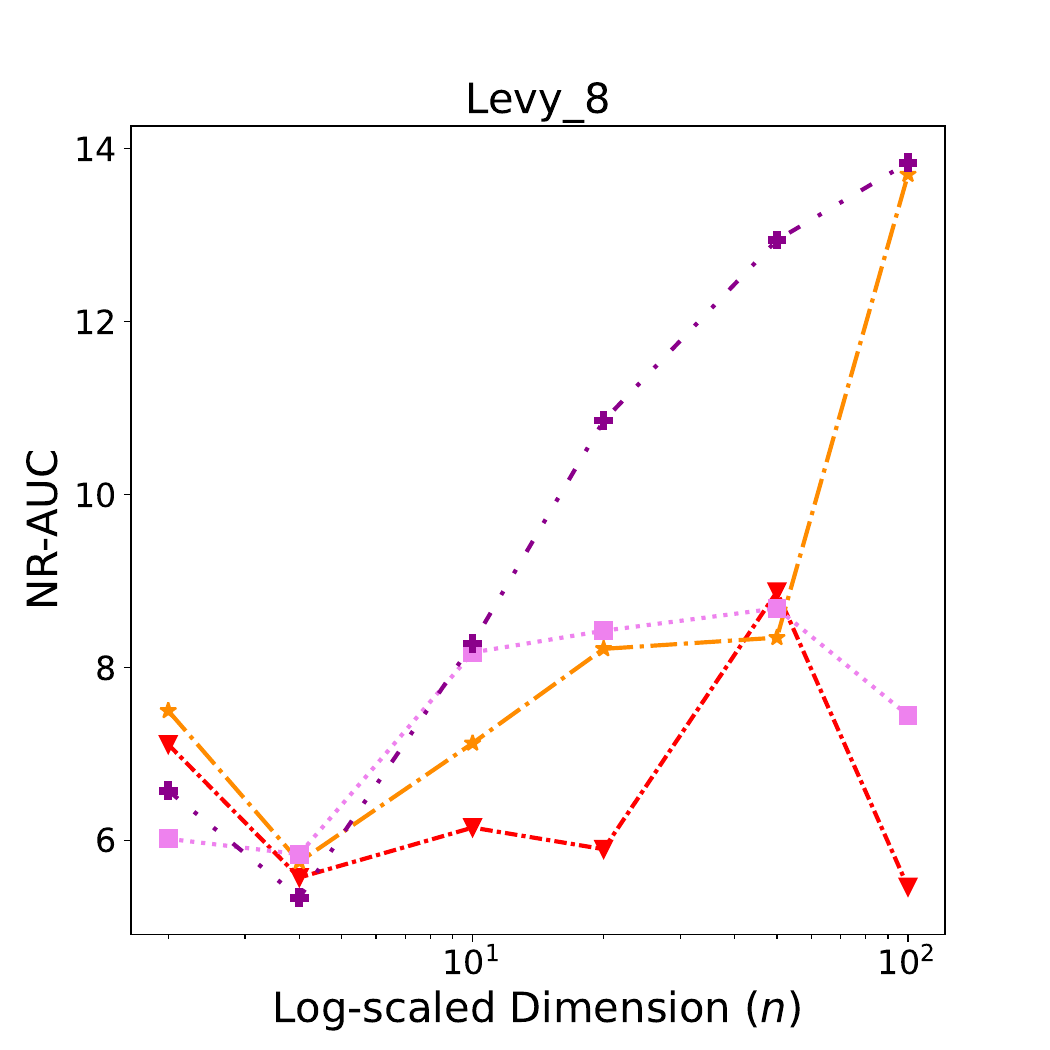}}
	\\
	\subfloat[$f^{best}_L$ -- Rastrigin]{\includegraphics[width=0.32\textwidth]{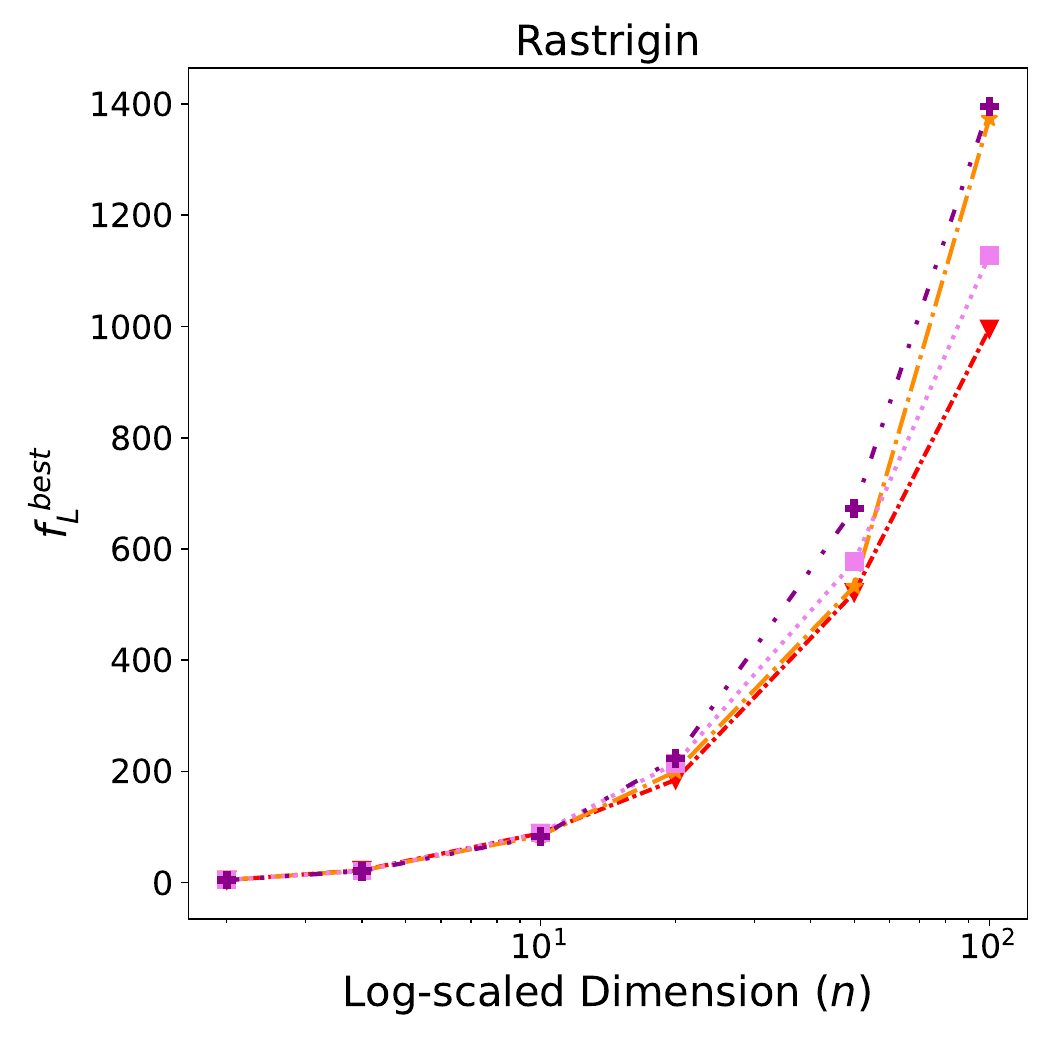}}
	\hfil
	\subfloat[$f^{best}_L$ -- Ackley\_1]{\includegraphics[width=0.32\textwidth]{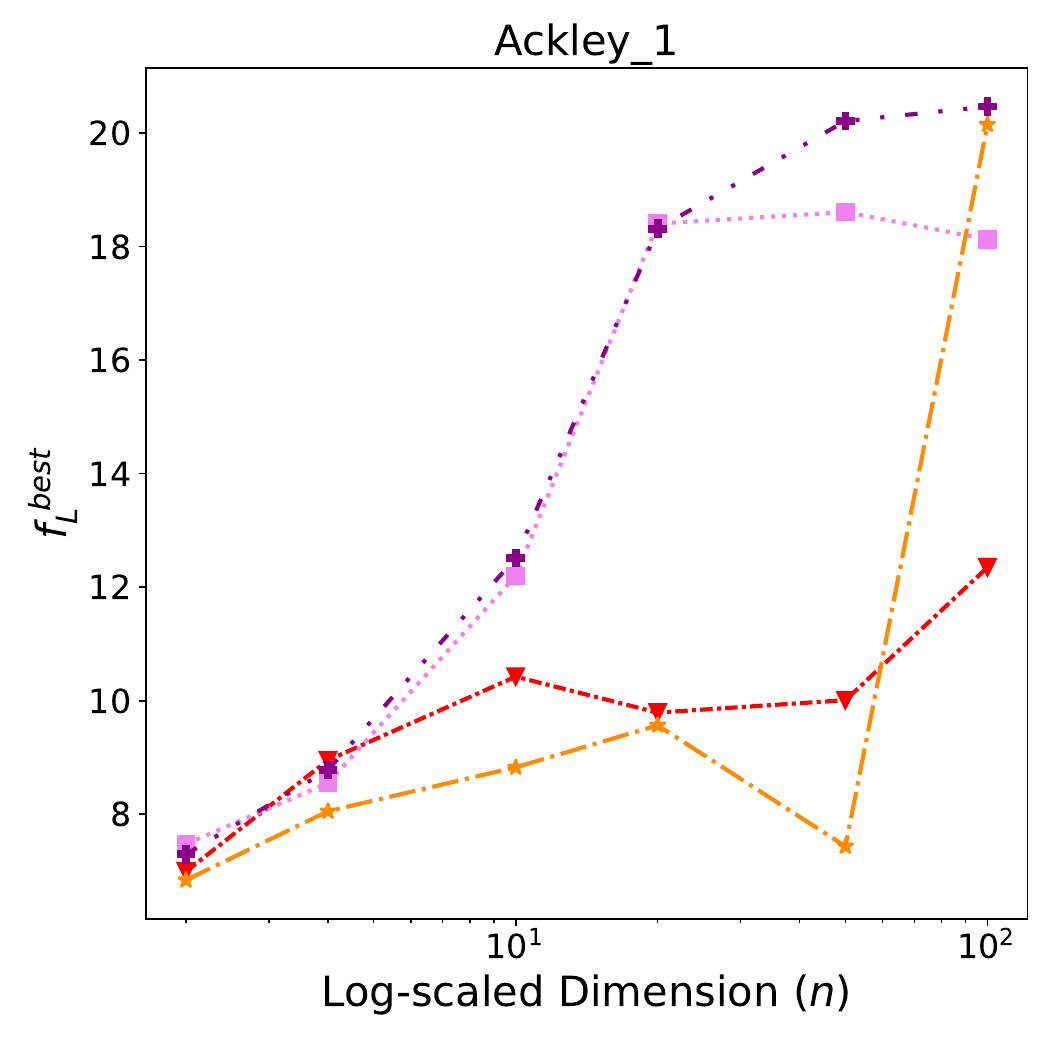}}
	\hfil
	\subfloat[$f^{best}_L$ -- Levy\_8]{\includegraphics[width=0.32\textwidth]{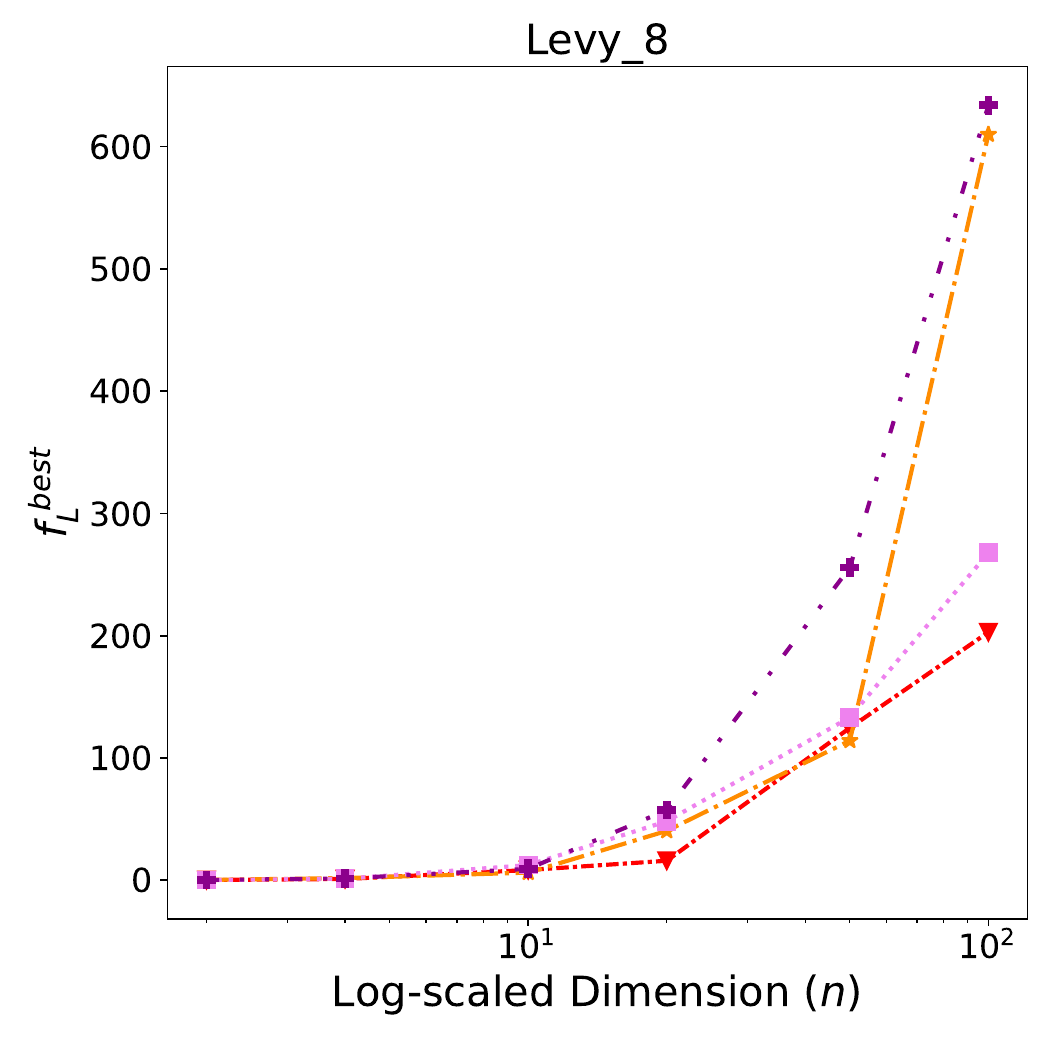}}
	\caption{Plots of the \textit{NR-AUC} and $f^{best}_L$ metrics values achieved by \texttt{NSMA(X)}, \texttt{NSMA(F)}, \texttt{NSGA-II(X)} and \texttt{NSGA-II(F)} on Rastrigin, Ackley\_1 and Levy\_8 (see Table \ref{tab::test_functions}) for values of $n \in \{2, 4, 10, 20, 50, 100\}$. For a better visualization, the x-axis is log-scaled.}
	\label{fig::auc_lastf_dim_plot}
\end{figure}

\subsection{Scalability of the BOO Based Acquisition Strategy Variants}
\label{subsec::clustering_approaches_comparison}

Attested the goodness of our BOO based acquisition strategy, we now focus on the four variants of the latter proposed for the computational experiments, i.e., \texttt{NSMA(X)}, \texttt{NSMA(F)}, \texttt{NSGA-II(X)} and \texttt{NSGA-II(F)}. In particular, we are interested on their scalability, i.e., how well they perform as the value for the dimensionality $n$ increases. In Figure \ref{fig::auc_lastf_dim_plot}, the \textit{NR-AUC} and $f^{best}_L$ metrics values achieved by the four methodologies on a subset of functions tested with $n \in \{2, 4, 10, 20, 50, 100\}$ are shown.

While for low values for $n$ the performance of the four approaches were similar, for high dimensionalities \texttt{NSMA(X)} and \texttt{NSGA-II(X)} got better overall results than \texttt{NSMA(F)} and \texttt{NSGA-II(F)}, respectively. The situation is particularly stressed when $n = 100$, where the gap between the first two methodologies and the second ones is clear. As mentioned in Section \ref{subsec::clustering}, clustering in the variables space could be more effective against the \textit{boundary issue}, allowing to better escape from evaluating points always near to the problem boundaries. The results seem to confirm this hypothesis: exploring feasible space regions different from the ones containing the problem boundaries, \texttt{NSMA(X)} and \texttt{NSGA-II(X)} managed to find better solutions for the original problem \eqref{eq::bo-prob}.

Finally, we again observe the overall superiority of \texttt{NSMA} w.r.t.\ \texttt{NSGA-II}, regardless the employed clustering approach. Thus, after analyzing the results of this section and the previous one, we decided to use only \texttt{NSMA(X)} and \texttt{NSMA(F)} for the rest of the section.

\subsection{Comparisons with the State-of-the-art Acquisition Strategies}
\label{subsec::comparison_state-of-the-art}

\begin{figure*}[!t]
	\centering
	\subfloat[HolderTable\_2, $n = 2$]{\includegraphics[width=0.32\textwidth]{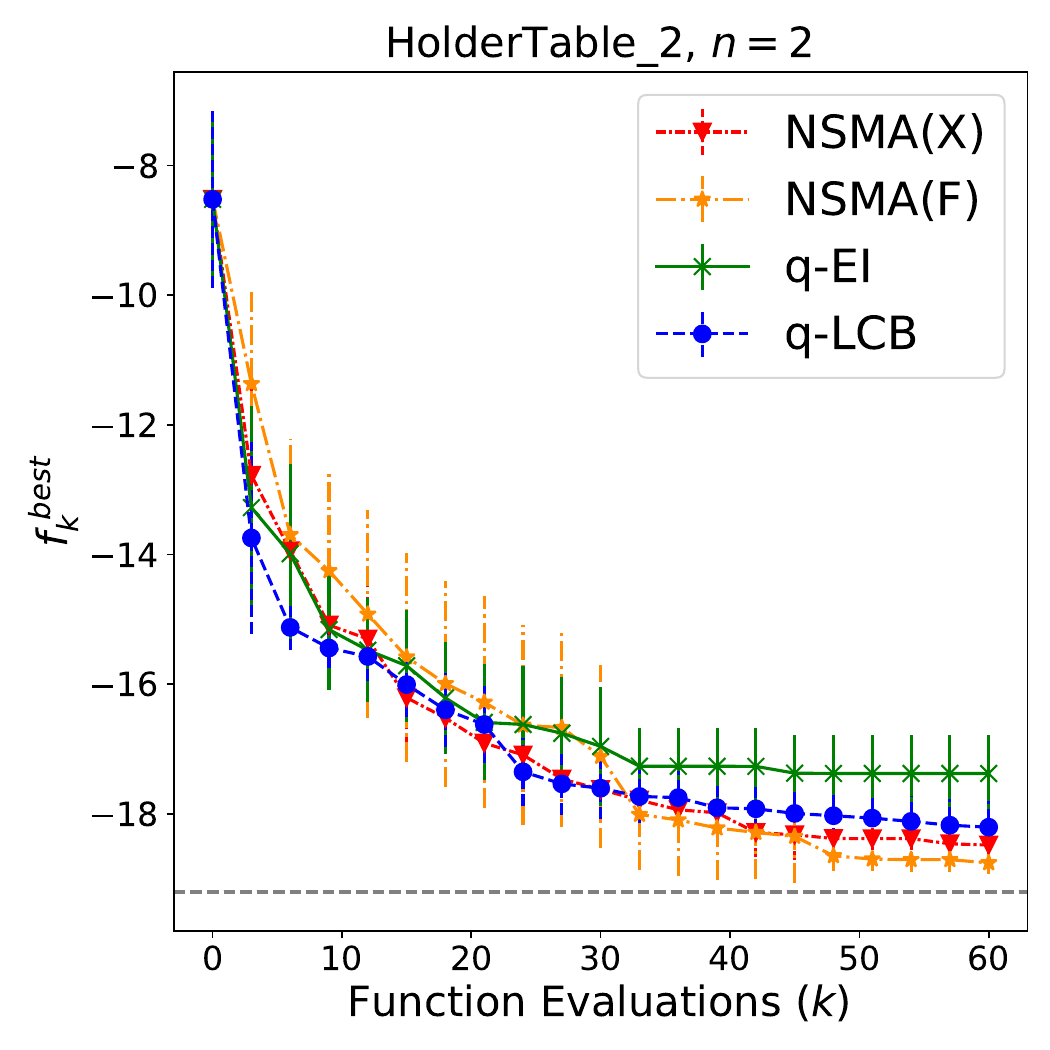}}
	\hfil
	\subfloat[Hartmann, $n = 6$]{\includegraphics[width=0.32\textwidth]{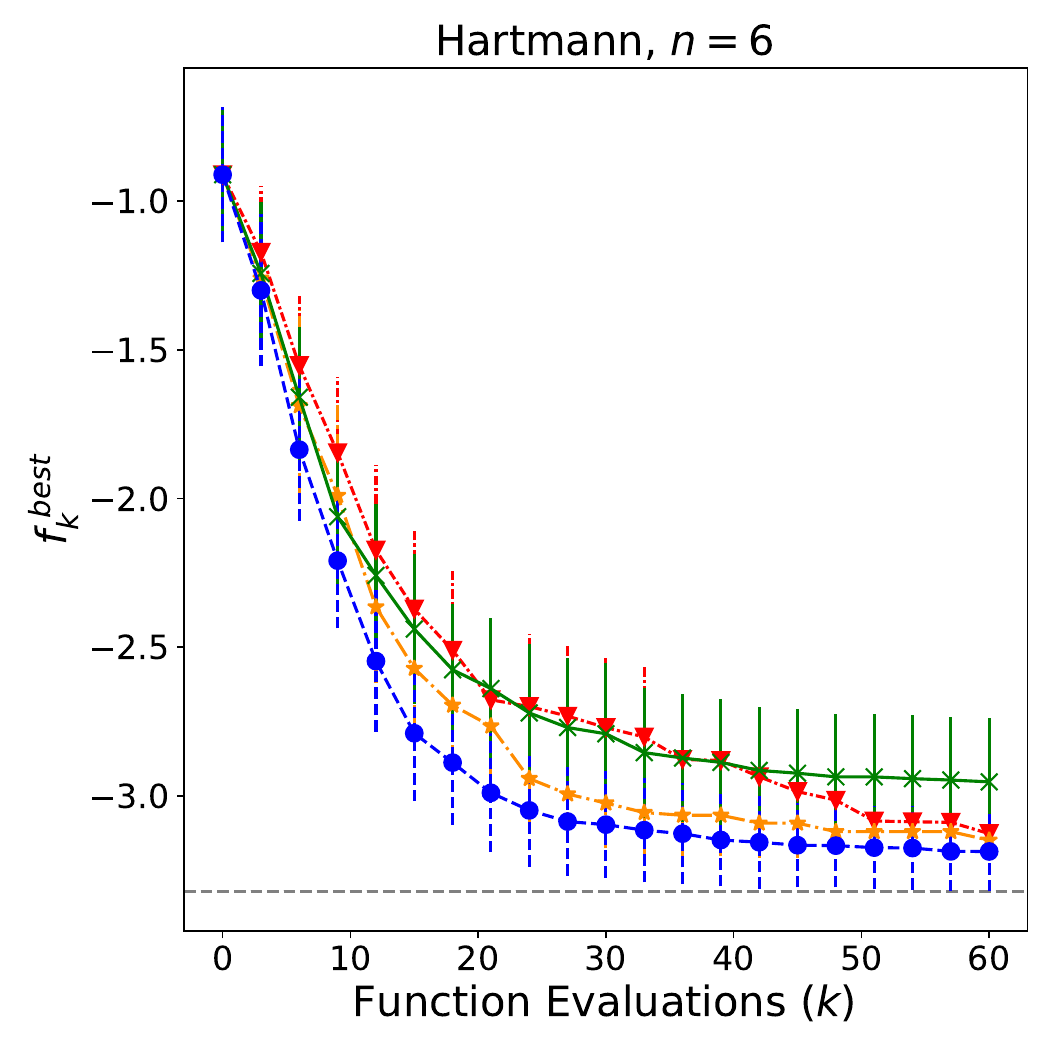}}
	\\
	\subfloat[Rosenbrock, $n = 20$]{\includegraphics[width=0.32\textwidth]{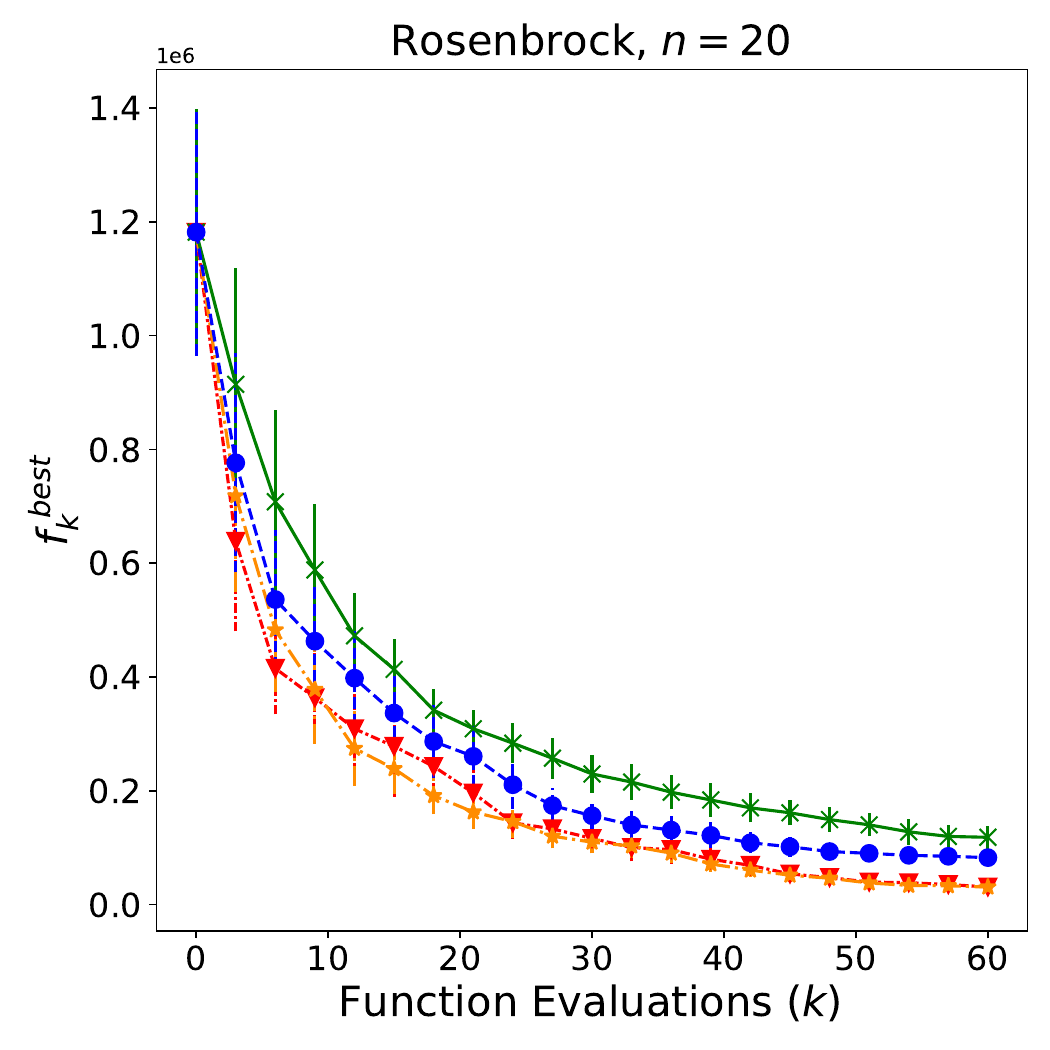}}
	\hfil
	\subfloat[Ackley\_1, $n = 20$]{\includegraphics[width=0.32\textwidth]{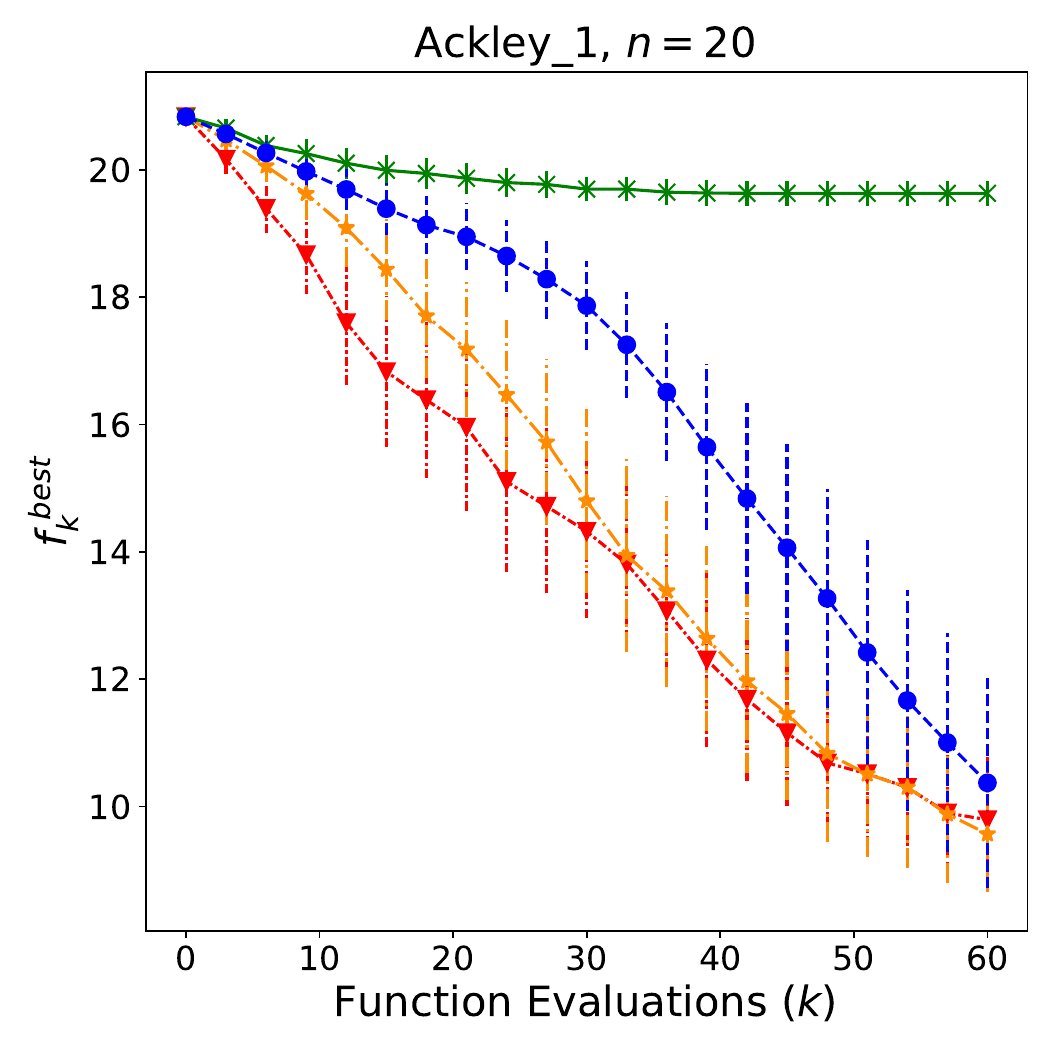}}
	\hfil
	\subfloat[Rastrigin, $n = 50$]{\includegraphics[width=0.32\textwidth]{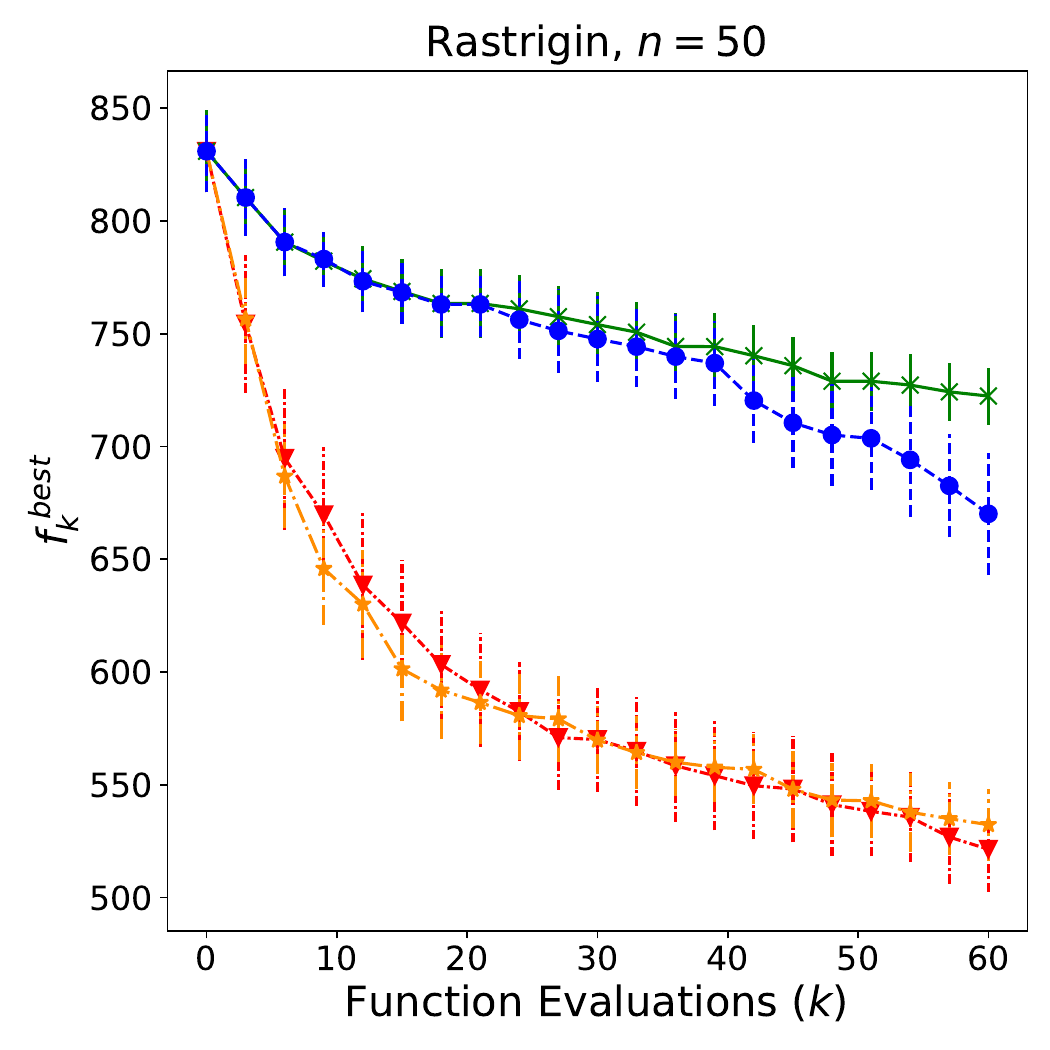}}
	\\
	\subfloat[Schewefel, $n = 100$]{\includegraphics[width=0.32\textwidth]{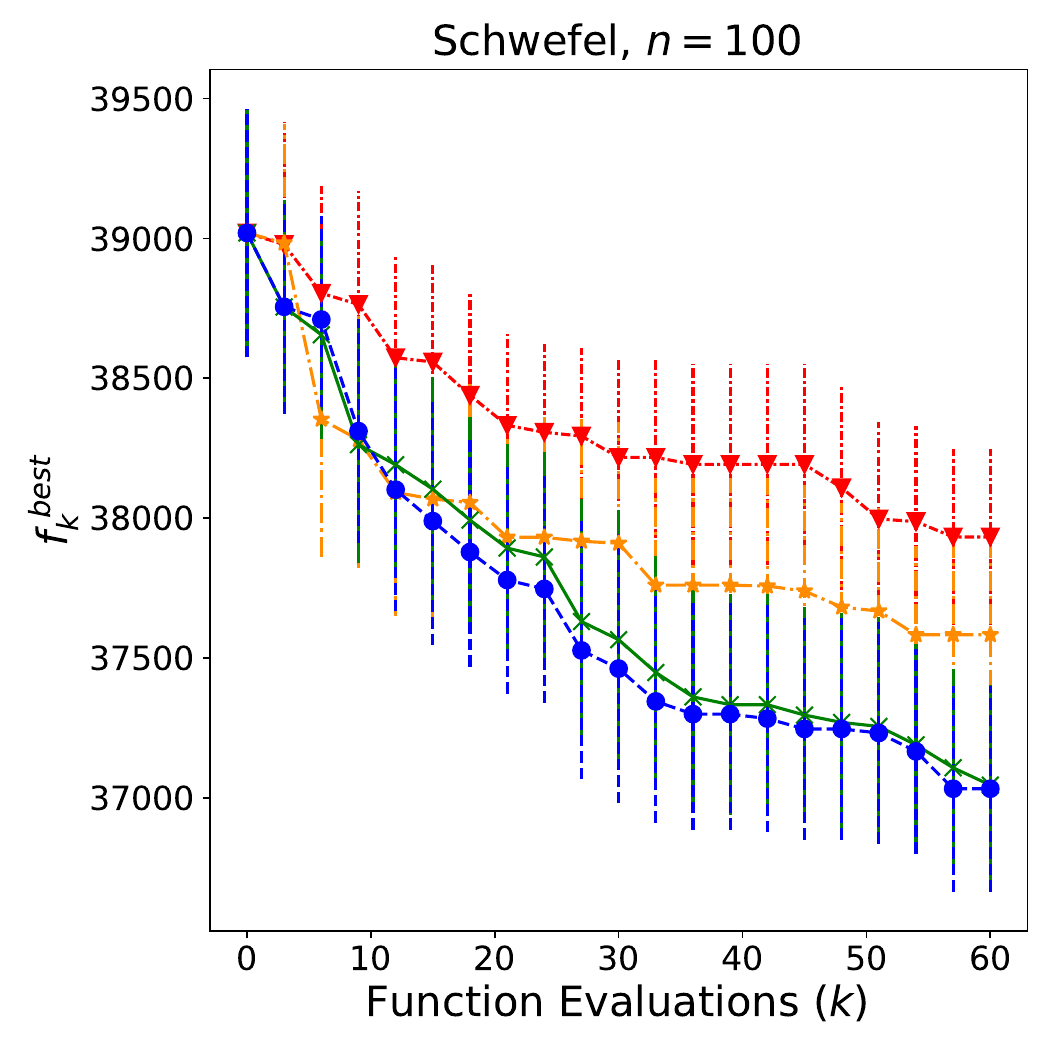}}
	\hfil
	\subfloat[Levy\_8, $n = 100$]{\includegraphics[width=0.32\textwidth]{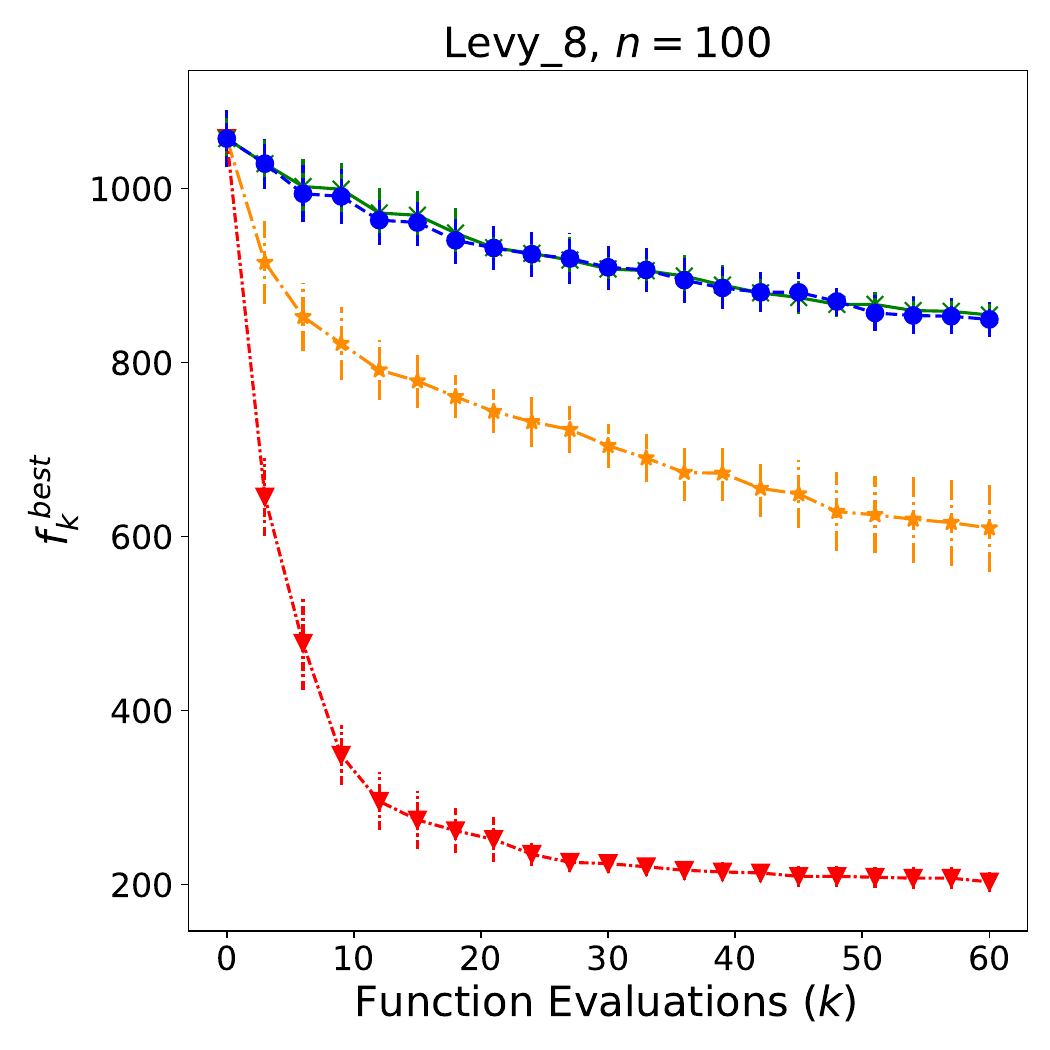}}
	\hfil
	\subfloat[Alpine\_1, $n = 100$]{\includegraphics[width=0.32\textwidth]{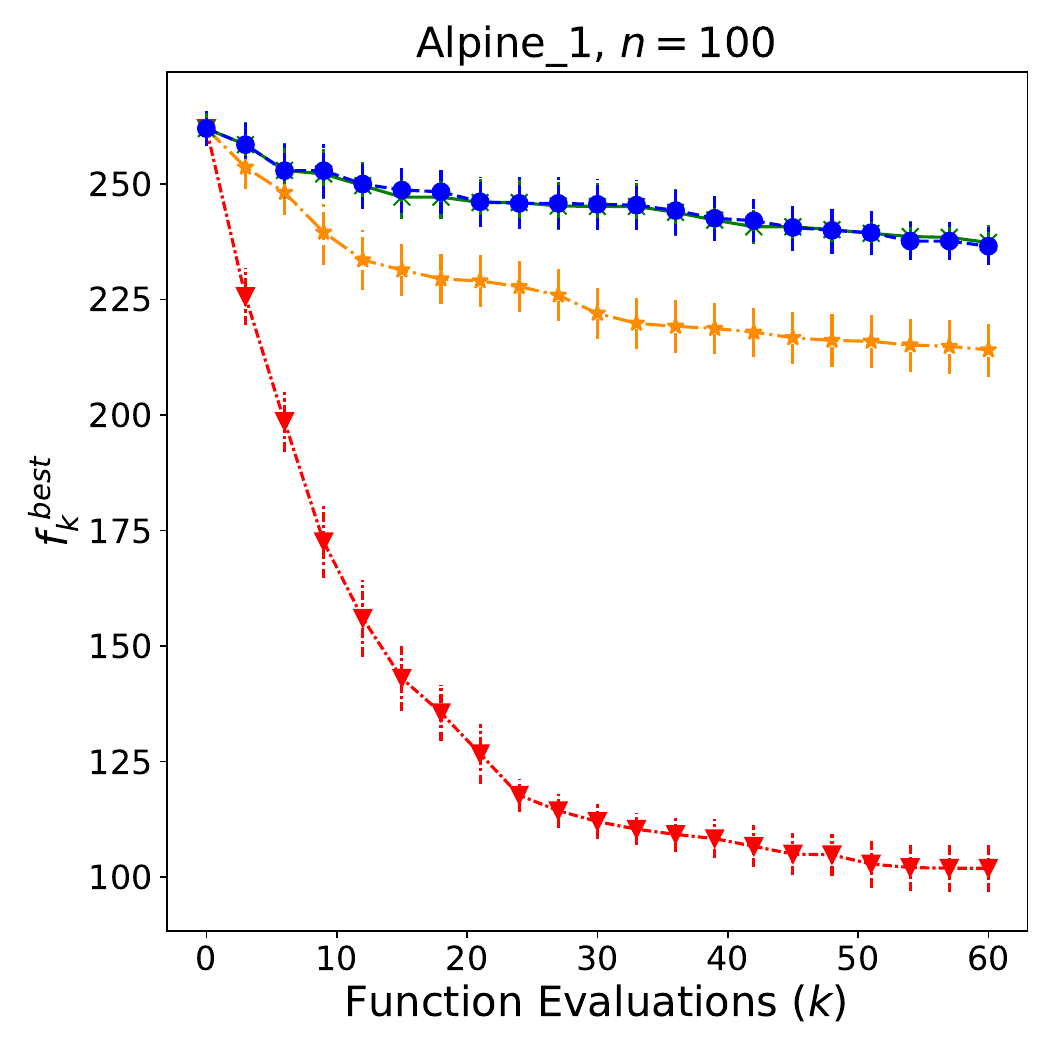}}
	\caption{Plots of the $f^{best}_k$ metric w.r.t.\ the number of function evaluations $k$ for \texttt{NSMA(X)}, \texttt{NSMA(F)}, \texttt{q-EI} and \texttt{q-LCB} on a set of selected functions (see Table \ref{tab::test_functions}). When $f^\star \ne 0$, it is represented by a gray dashed line.}
	\label{fig::nsma_vs_bo}
\end{figure*}

In this section, we investigate the performance of \texttt{NSMA(X)}, \texttt{NSMA(F)}, \texttt{q-EI} and \texttt{q-LCB} more specifically. In Figure \ref{fig::nsma_vs_bo}, we report the plots of the $f^{best}_k$ metric values obtained by the four acquisition methodologies on a set of selected functions, some of which were already addressed in Section \ref{subsec::preliminary_assessment}. 

In the low-dimensional instances ($n \le 6$), the proposed acquisition strategies were competitive: \texttt{NSMA(X)} and \texttt{NSMA(F)} had a similar trend as \texttt{q-LCB}, especially performing well in terms of $f^{best}_L$. Unlike the three mentioned approaches, in these functions \texttt{q-EI} seems to struggle, having competitive performance but never outperforming its competitors. 

As the value for $n$ increased, \texttt{NSMA(X)} and \texttt{NSMA(F)} turned out to be the best algorithms overall. The standard acquisition strategies \texttt{q-EI} and \texttt{q-LCB} still had good results on some scenarios: both approaches had similar behavior as our proposals on the Rosenbrock function with $n=20$; a competitive result was also obtained by \texttt{q-LCB} on the Ackley\_1 function at the same dimensionality. However, for higher values of $n$, the general superiority of \texttt{NSMA(X)} and \texttt{NSMA(F)} become clear. 

\begin{figure*}[!b]
	\centering
	\subfloat[Schwefel, $n=100$]{\includegraphics[width=0.32\textwidth]{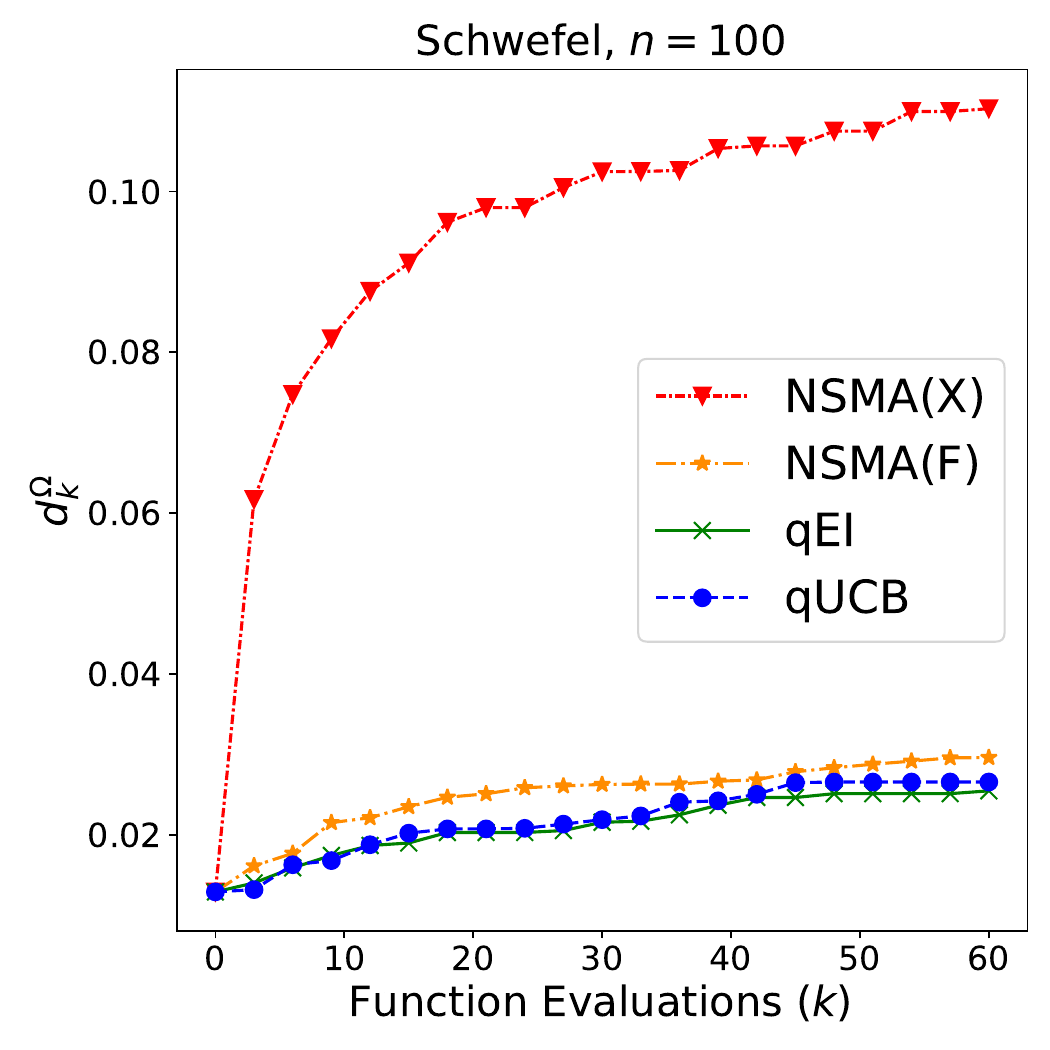}}
	\hfil
	\subfloat[Levy\_8, $n=100$]{\includegraphics[width=0.32\textwidth]{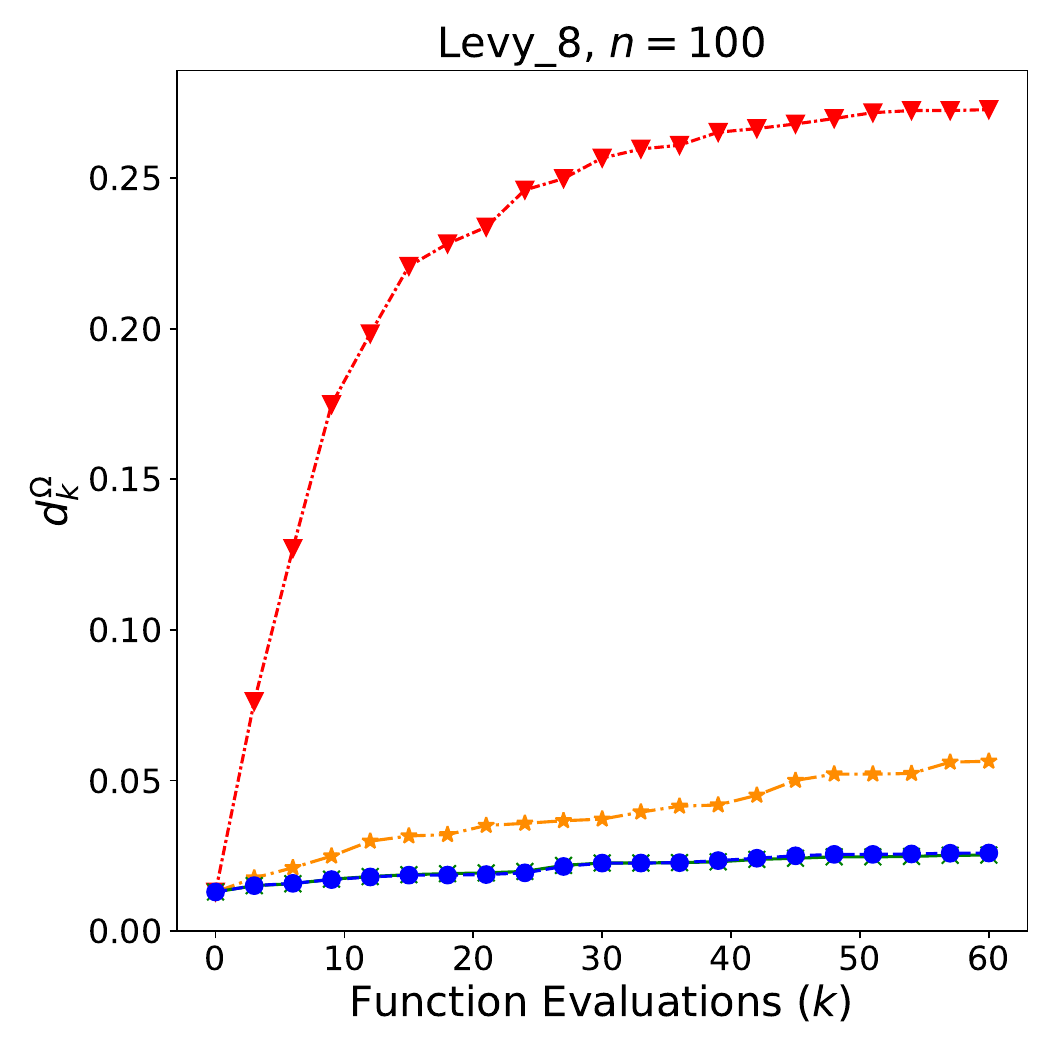}}
	\hfil
	\subfloat[Alpine\_1, $n=100$]{\includegraphics[width=0.32\textwidth]{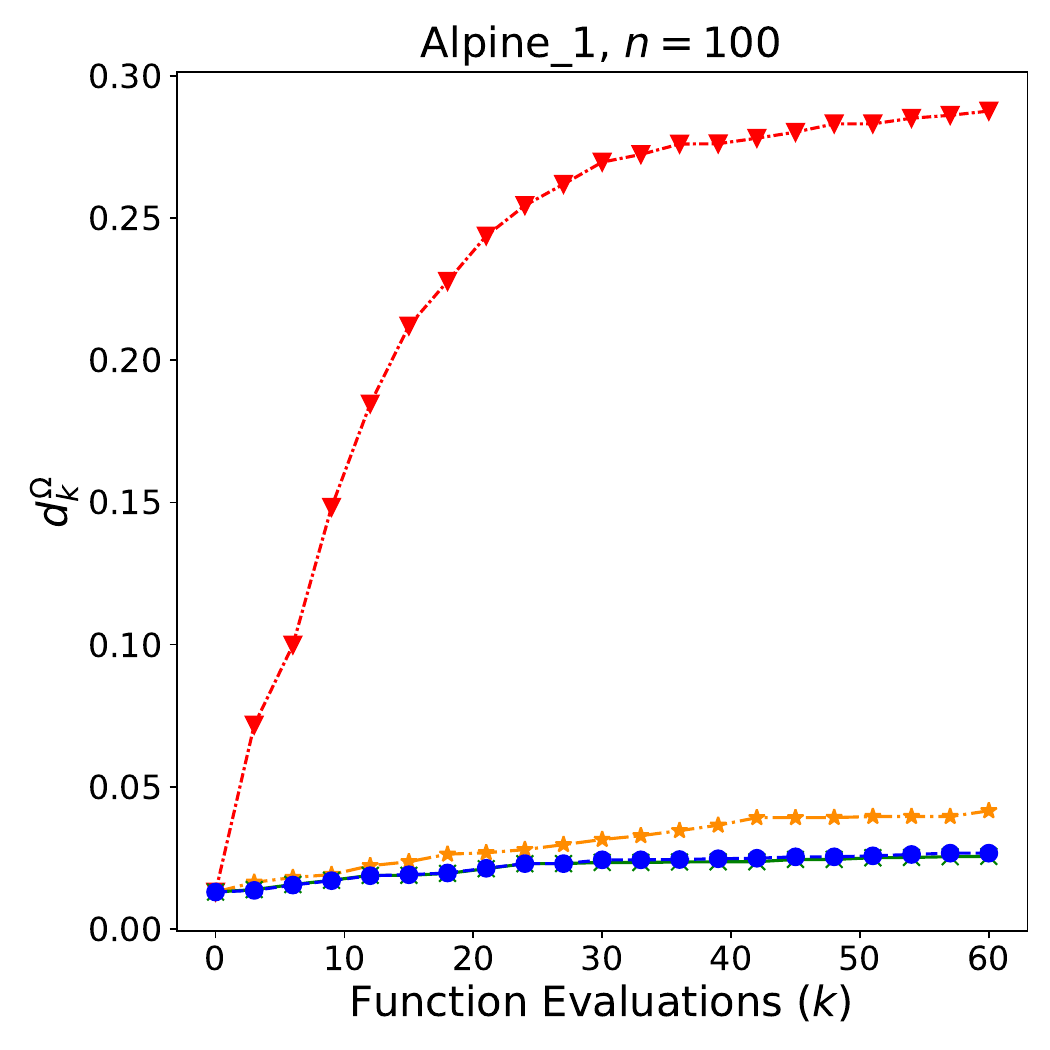}}
	\caption{Plot of $d^\Omega_k$ \eqref{eq::dist-from-bounds} w.r.t.\ the number of function evaluations $k$ for \texttt{NSMA(X)}, \texttt{NSMA(F)}, \texttt{q-EI} and \texttt{q-LCB} on the Schwefel, Levy\_8 and Alpine\_1 functions (see Table \ref{tab::test_functions}) with $n=100$.}
	\label{fig:boundary_issues_100}
\end{figure*}

\begin{table*}[!b]
	\caption{\textit{NR-AUC} and $f^{best}_L$ values achieved by \texttt{NSMA(X)}, \texttt{NSMA(F)}, \texttt{q-EI} and \texttt{q-LCB} on the first set of selected functions shown in Figure \ref{fig::nsma_vs_bo}. The marked-in-bold values and the underlined ones are, respectively, the best and second best for a given metric on a considered function.}
	\label{tab::nsma_vs_bo_1}
	\centering
	\footnotesize	
	\begin{tabular}{c||ccc}
		\toprule
		\multirow{2}{*}{\textbf{Function}}&\multirow{2}{*}{\textbf{Method}}&\multirow{2}{*}{\textbf{NR{-}AUC}}&\multirow{2}{*}{\textbf{$\mathbf{f^{best}_L}$}}\\
		&&&\\
		\midrule
		\midrule
		
		&\texttt{NSMA(X)}&\underline{4.470}&\underline{-18.484}\\
		\cmidrule{2-4}
		HolderTable\_2&\texttt{NSMA(F)}&4.891&\textbf{-18.758}\\
		\cmidrule{2-4}
		($n=2$)&\texttt{q-EI}&5.473&{-}17.379\\
		\cmidrule{2-4}
		&\texttt{q-LCB}&\textbf{4.461}&{-}18.208\\
		
		\midrule
		
		&\texttt{NSMA(X)}&6.292&{-}3.127\\
		\cmidrule{2-4}
		Hartmann&\texttt{NSMA(F)}&\underline{5.092}&\underline{-3.149}\\
		\cmidrule{2-4}
		($n=6$)&\texttt{q-EI}&6.282&{-}2.953\\
		\cmidrule{2-4}
		&\texttt{q-LCB}&\textbf{4.229}&\textbf{-3.187}\\
		
		\midrule
		
		&\texttt{NSMA(X)}&\underline{3.385}&\underline{30865.769}\\
		\cmidrule{2-4}
		Rosenbrock&\texttt{NSMA(F)}&\textbf{3.346}&\textbf{30617.591}\\
		\cmidrule{2-4}
		($n=20$)&\texttt{q-EI}&5.611&118150.206\\
		\cmidrule{2-4}
		&\texttt{q-LCB}&4.388&82420.471\\
		
		\midrule
		
		&\texttt{NSMA(X)}&\textbf{13.817}&\underline{9.790}\\
		\cmidrule{2-4}
		Ackley\_1&\texttt{NSMA(F)}&\underline{14.384}&\textbf{9.566}\\
		\cmidrule{2-4}
		($n=20$)&\texttt{q-EI}&19.078&19.629\\
		\cmidrule{2-4}
		&\texttt{q-LCB}&16.081&10.370\\
		
		\bottomrule
	\end{tabular}
\end{table*}

The results on the instances with $n=100$ are worth mentioning separately. On the Schwefel function, our approaches struggled to reach the performance of \texttt{q-EI} and \texttt{q-LCB}. A possible justification of this fact could be the Schwefel function global optimum position within the feasible set, which is quite near to the boundary. As mentioned in Section \ref{subsec::clustering}, \texttt{NSMA(X)} represents an heuristic solution to the \textit{boundary issue}: it could allow to escape from a massive exploration of the boundaries, which is a typical behavior on high-dimensional instances of acquisition strategies based on the posterior mean ($\mu_k(\cdot)$) and standard deviation ($\sigma_k(\cdot)$) functions. By design, \texttt{NSMA(X)} will be thus driven to find good solutions in the interior. \texttt{NSMA(F)} has a similar, albeit less pronounced, behavior; indeed, it obtained better results than \texttt{NSMA(X)} in the Schwefel function, but, at the same time, it did not perform good enough to reach the performance of \texttt{q-EI} and \texttt{q-LCB}. As a proof of this bias of the \texttt{NSMA}-based acquisition strategies, in Figure \ref{fig:boundary_issues_100} we show the $d^\Omega_k$ \eqref{eq::dist-from-bounds} values achieved by \texttt{NSMA(X)}, \texttt{NSMA(F)}, \texttt{q-EI} and \texttt{q-LCB} on the three $n=100$ instances analyzed in this section. In the Schwefel function, the two standard acquisition strategies took advantage of the blind exploration of the problem boundaries, reaching solutions that, however, are extremely far from the global optimum of the problem. Moreover, the boundary exploration cannot be beneficial with the majority of the benchmark functions considered in the Bayes-Opt literature, which present a global optimum in a more internal position of the feasible set. In these last, most common scenarios, \texttt{NSMA(X)} and \texttt{NSMA(F)} proved to be highly effective: on the Levy\_8 and Alpine\_1 functions with $n=100$, the gap between the two \texttt{NSMA}-based acquisition strategies and the standard ones is clear; in particular, we can observe the remarkable results obtained by \texttt{NSMA(X)} w.r.t.\ all the other competitors.

\begin{table*}
	\caption{\textit{NR-AUC} and $f^{best}_L$ values achieved by \texttt{NSMA(X)}, \texttt{NSMA(F)}, \texttt{q-EI} and \texttt{q-LCB} on the second set of selected functions shown in Figure \ref{fig::nsma_vs_bo}. The marked-in-bold values and the underlined ones are, respectively, the best and second best for a given metric on a considered function.}
	\label{tab::nsma_vs_bo_2}
	\centering
	\footnotesize	
	\begin{tabular}{c||ccc}
		\toprule
		\multirow{2}{*}{\textbf{Function}}&\multirow{2}{*}{\textbf{Method}}&\multirow{2}{*}{\textbf{NR{-}AUC}}&\multirow{2}{*}{\textbf{$\mathbf{f^{best}_L}$}}\\
		&&&\\
		\midrule
		\midrule
		
		&\texttt{NSMA(X)}&\underline{14.311}&\textbf{521.437}\\
		\cmidrule{2-4}
		Rastrigin&\texttt{NSMA(F)}&\textbf{14.269}&\underline{532.328}\\
		\cmidrule{2-4}
		($n=50$)&\texttt{q-EI}&18.206&722.386\\
		\cmidrule{2-4}
		&\texttt{q-LCB}&17.926&670.171\\
		
		\midrule
		
		&\texttt{NSMA(X)}&19.650&37932.468\\
		\cmidrule{2-4}
		Schwefel&\texttt{NSMA(F)}&19.455&37582.778\\
		\cmidrule{2-4}
		($n=100$)&\texttt{q-EI}&\underline{19.337}&\underline{37044.992}\\
		\cmidrule{2-4}
		&\texttt{q-LCB}&\textbf{19.309}&\textbf{37032.126}\\
		
		\midrule
		
		&\texttt{NSMA(X)}&\textbf{5.460}&\textbf{202.875}\\
		\cmidrule{2-4}
		Levy\_8&\texttt{NSMA(F)}&\underline{13.699}&\underline{609.754}\\
		\cmidrule{2-4}
		($n=100$)&\texttt{q-EI}&17.456&854.655\\
		\cmidrule{2-4}
		&\texttt{q-LCB}&17.400&849.608\\
		
		\midrule
		
		&\texttt{NSMA(X)}&\textbf{10.057}&\textbf{101.890}\\
		\cmidrule{2-4}
		Alpine\_1&\texttt{NSMA(F)}&\underline{17.301}&\underline{214.062}\\
		\cmidrule{2-4}
		($n=100$)&\texttt{q-EI}&18.735&237.291\\
		\cmidrule{2-4}
		&\texttt{q-LCB}&18.753&236.497\\
		
		\bottomrule
	\end{tabular}
\end{table*}

Figure \ref{fig::nsma_vs_bo} is supported by Tables \ref{tab::nsma_vs_bo_1}-\ref{tab::nsma_vs_bo_2}, where the values of the \textit{NR-AUC} and $f^{best}_L$ metrics obtained by the four methodologies are reported. In low-dimensional cases, \texttt{NSMA} turned out to be competitive, especially when clustering was performed on the objectives space (\texttt{NSMA(F)}). In the Schwefel function with $n = 100$, \texttt{q-LCB} and \texttt{q-EI} were the first and second best algorithms in terms of both metrics; in all the other high-dimensional instances, the memetic method outperformed the competitors, with \texttt{NSMA(X)} being by far the most effective methodology.

\subsection{Overall Performance Evaluation using Relative Gap Cumulative Distributions}
\label{subsec::perfomance_profiles}

\begin{figure*}[!h]
	\centering
	\subfloat[\textit{NR-AUC} for $q = 3$, $n < 10$]{\includegraphics[width=0.375\textwidth]{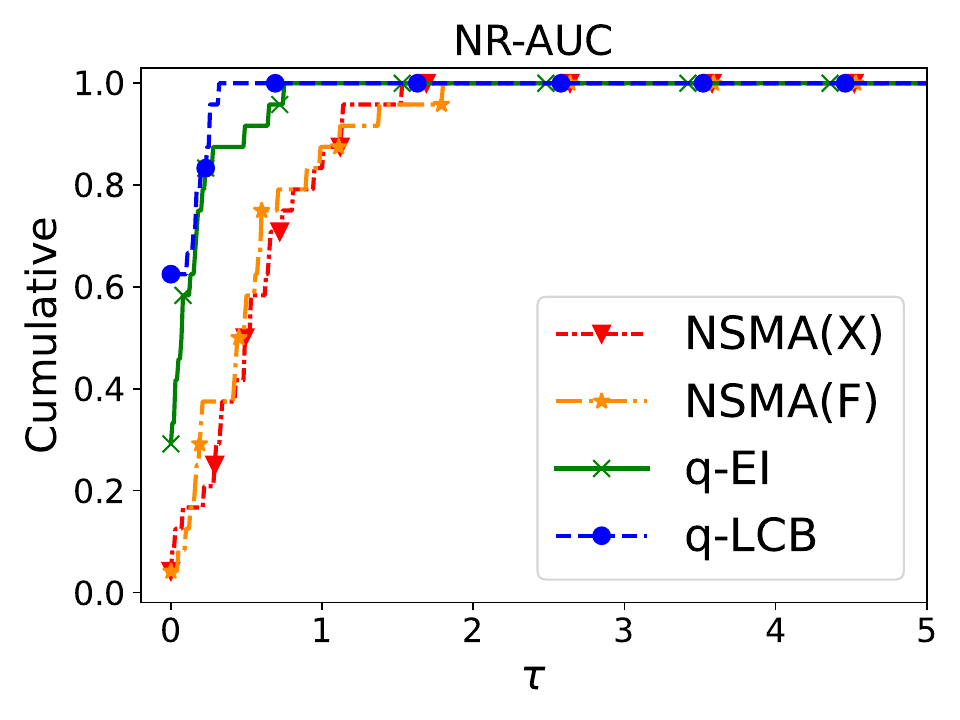}}
	\hfil
	\subfloat[$f^{best}_L$ for $q = 3$, $n < 10$]{\includegraphics[width=0.375\textwidth]{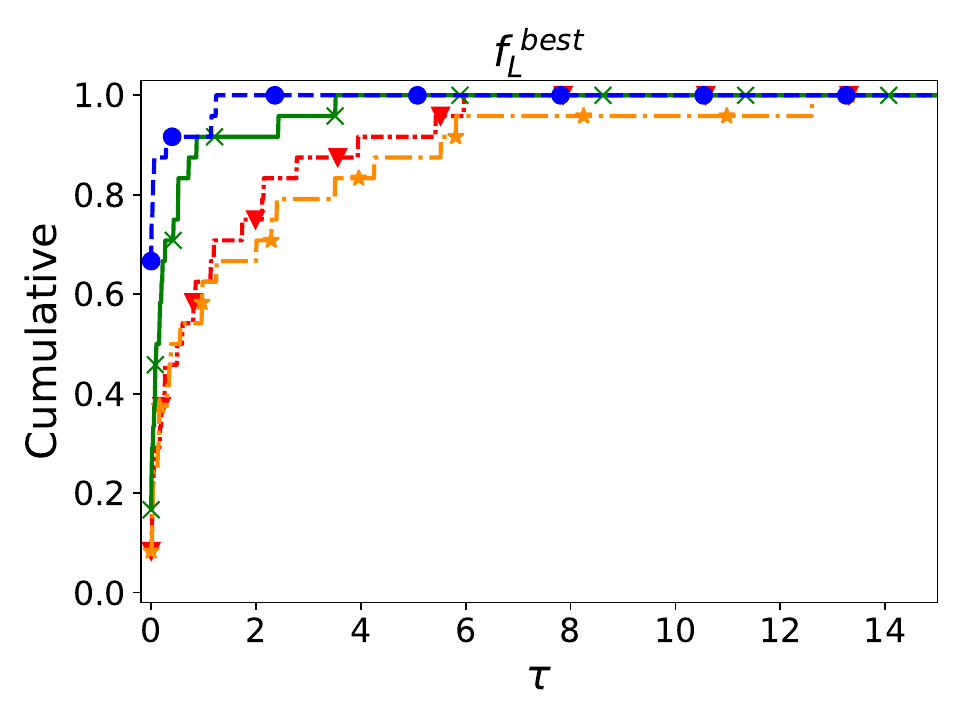}}
	\\
	\subfloat[\textit{NR-AUC} for $q = 3$, $n \ge 10$]{\includegraphics[width=0.375\textwidth]{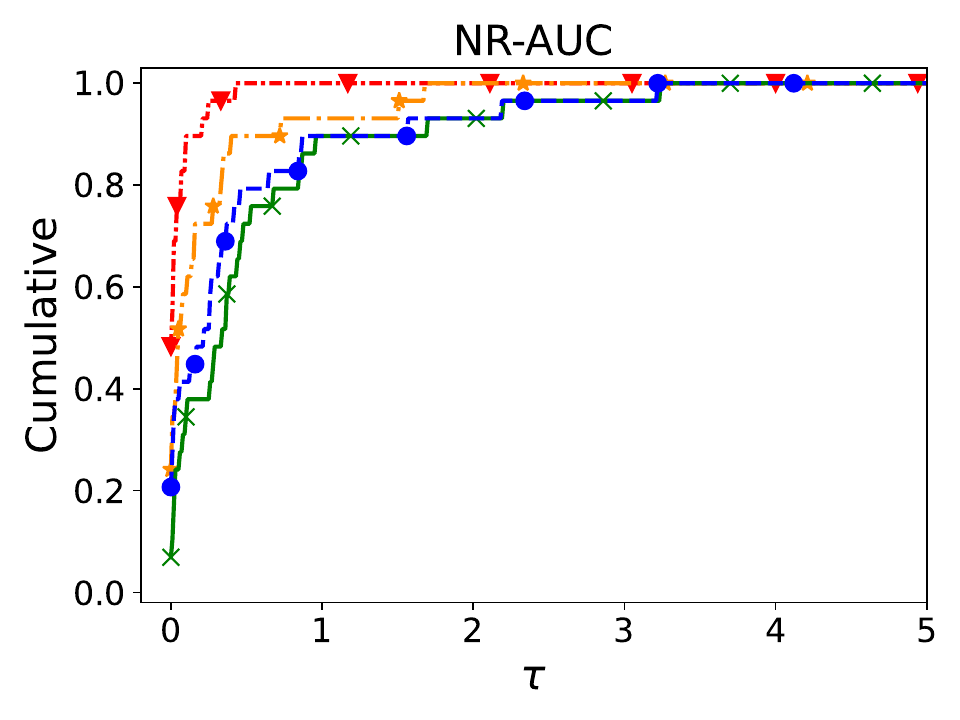}}
	\hfil
	\subfloat[$f^{best}_L$ for $q = 3$, $n \ge 10$]{\includegraphics[width=0.375\textwidth]{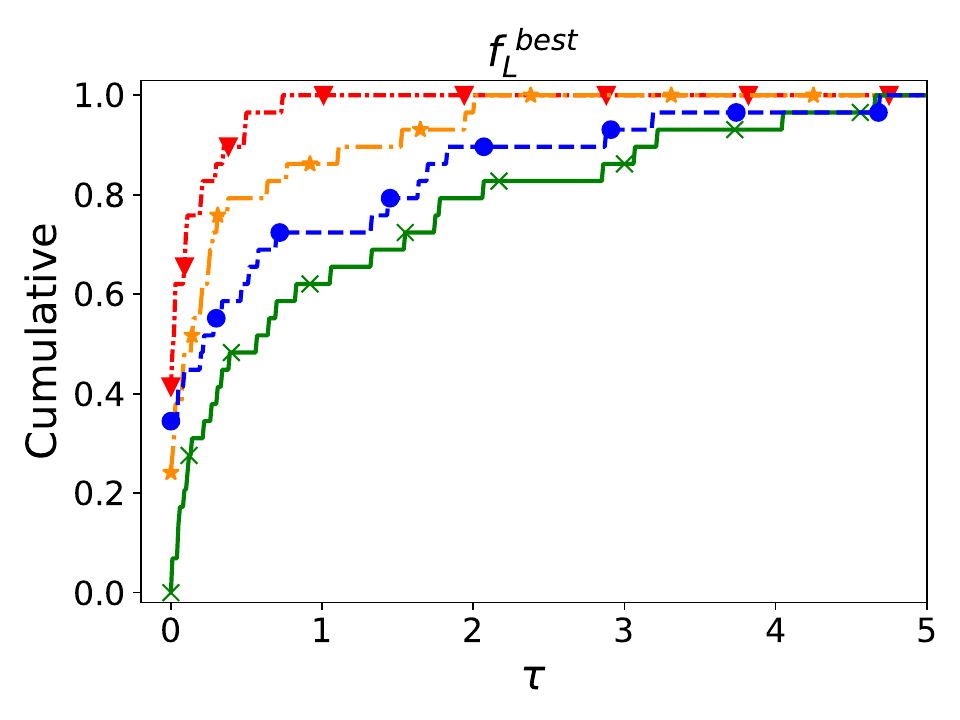}}
	\caption{Plots of the relative gap cumulative distribution for \texttt{NSMA(X)}, \texttt{NSMA(F)}, \texttt{q-EI} and \texttt{q-LCB} considering the low-dimensional ($n < 10$) and high-dimensional ($n \ge 10$) functions listed in Table \ref{tab::test_functions}. Note that the intervals of the x-axis are set for a better visualization of the methodologies results.}
	\label{fig::pp}
\end{figure*}

In the last section of computational experiments, we show the relative gap cumulative distributions for \texttt{NSMA(X)}, \texttt{NSMA(F)}, \texttt{q-EI} and \texttt{q-LCB} on the low-dimensional and high-dimensional functions separately. The distributions are plotted in Figure \ref{fig::pp}.

\begin{figure}
	\centering
	\subfloat[\textit{NR-AUC} for $q=2$, $n < 10$]{\includegraphics[width=0.375\textwidth]{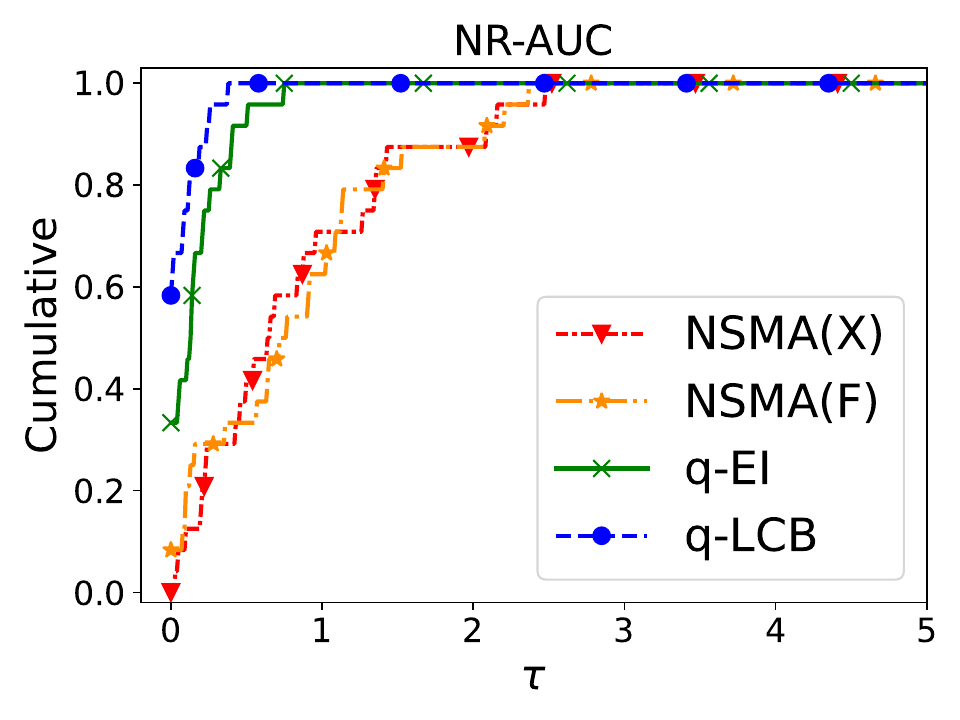}}
	\hfil
	\subfloat[\textit{NR-AUC} for $q=5$, $n < 10$]{\includegraphics[width=0.375\textwidth]{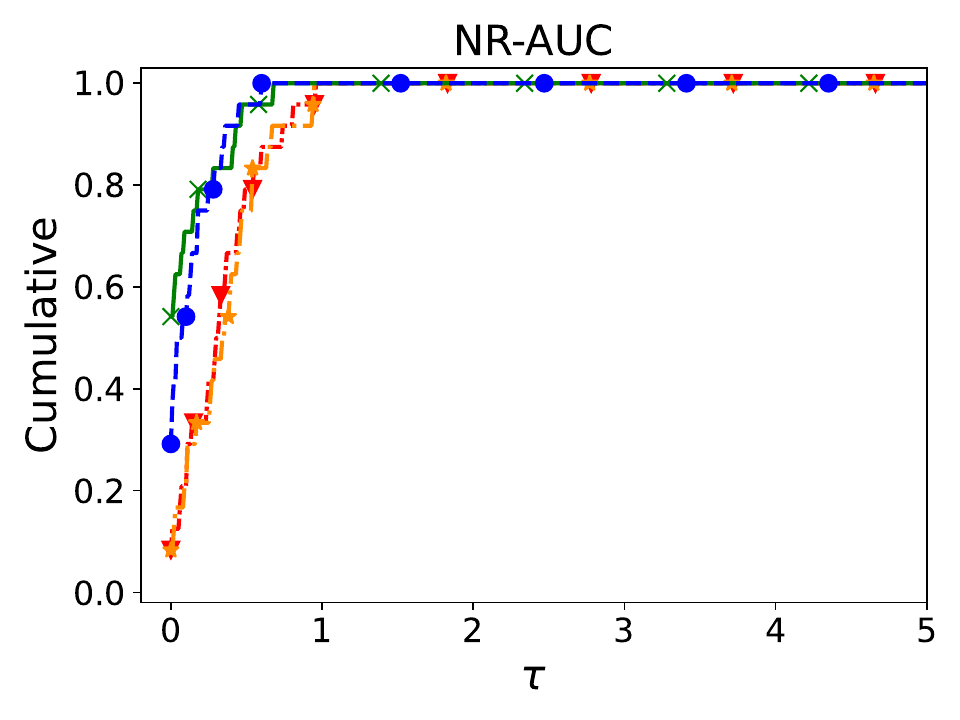}}
	\\
	\subfloat[\textit{NR-AUC} for $q=2$, $n \ge 10$]{\includegraphics[width=0.375\textwidth]{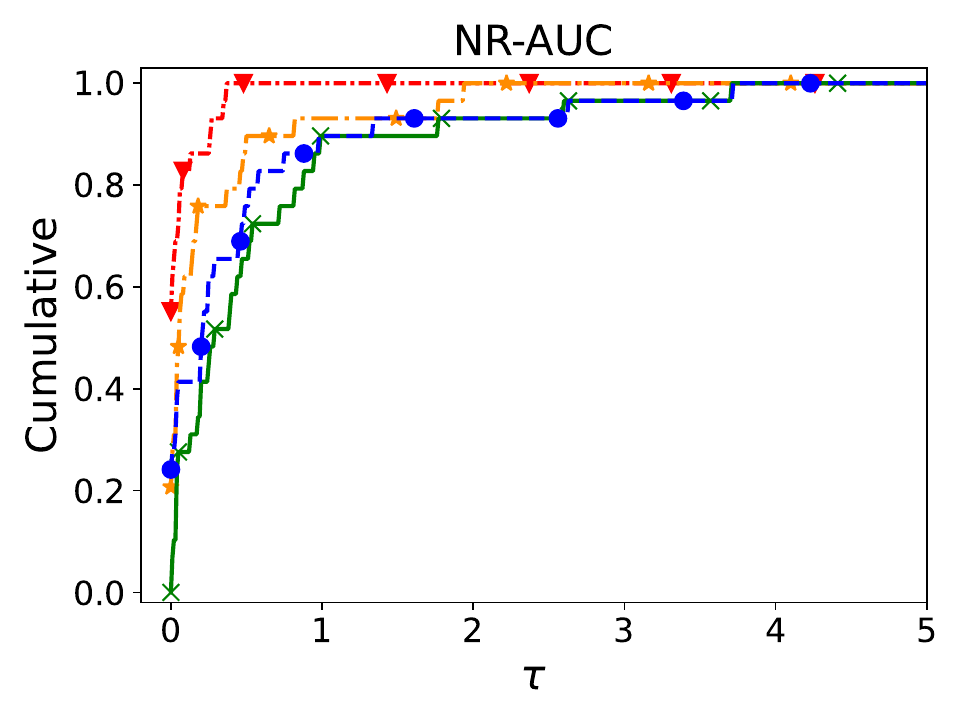}}
	\hfil
	\subfloat[\textit{NR-AUC} for $q=5$, $n \ge 10$]{\includegraphics[width=0.375\textwidth]{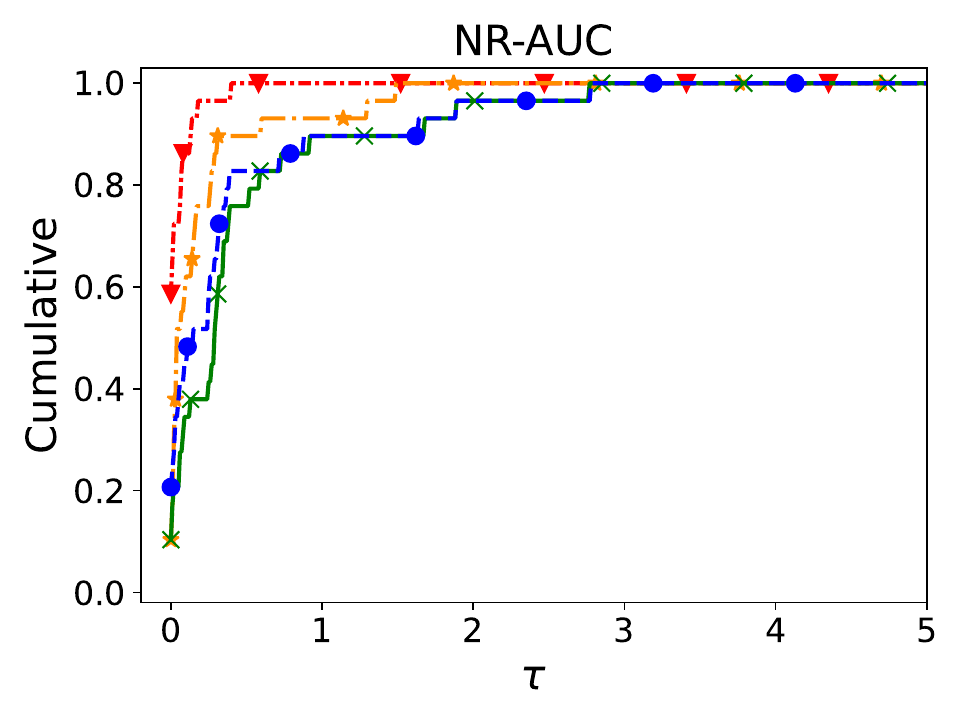}}
	\\
	\subfloat[$f^{best}_L$ for $q=2$, $n < 10$]{\includegraphics[width=0.375\textwidth]{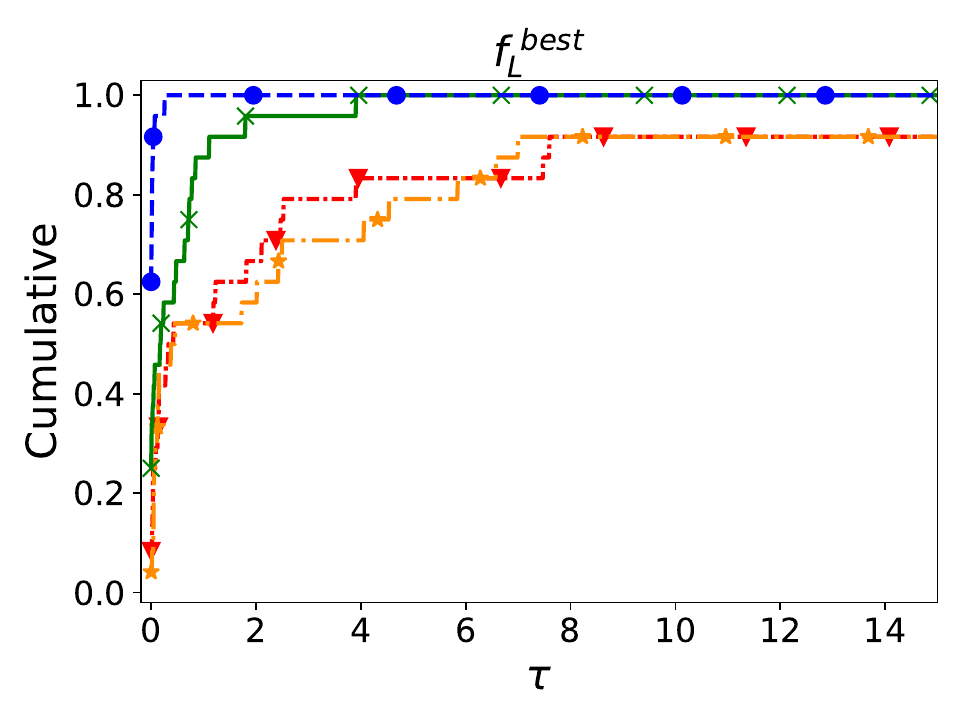}}
	\hfil
	\subfloat[$f^{best}_L$ for $q=5$, $n < 10$]{\includegraphics[width=0.375\textwidth]{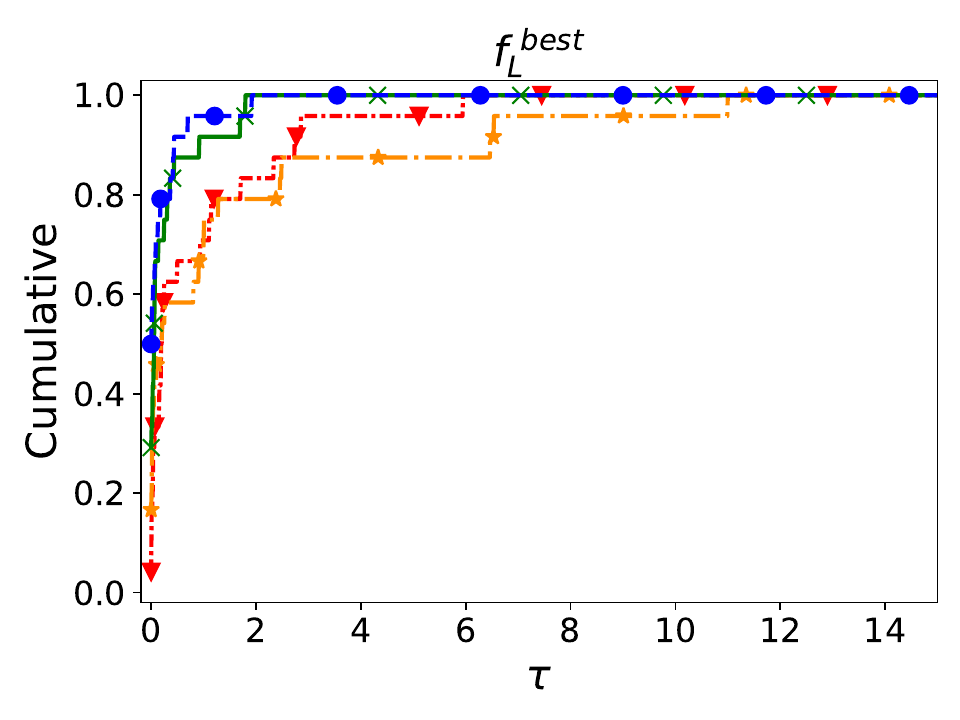}}
	\\
	\subfloat[$f^{best}_L$ for $q=2$, $n \ge 10$]{\includegraphics[width=0.375\textwidth]{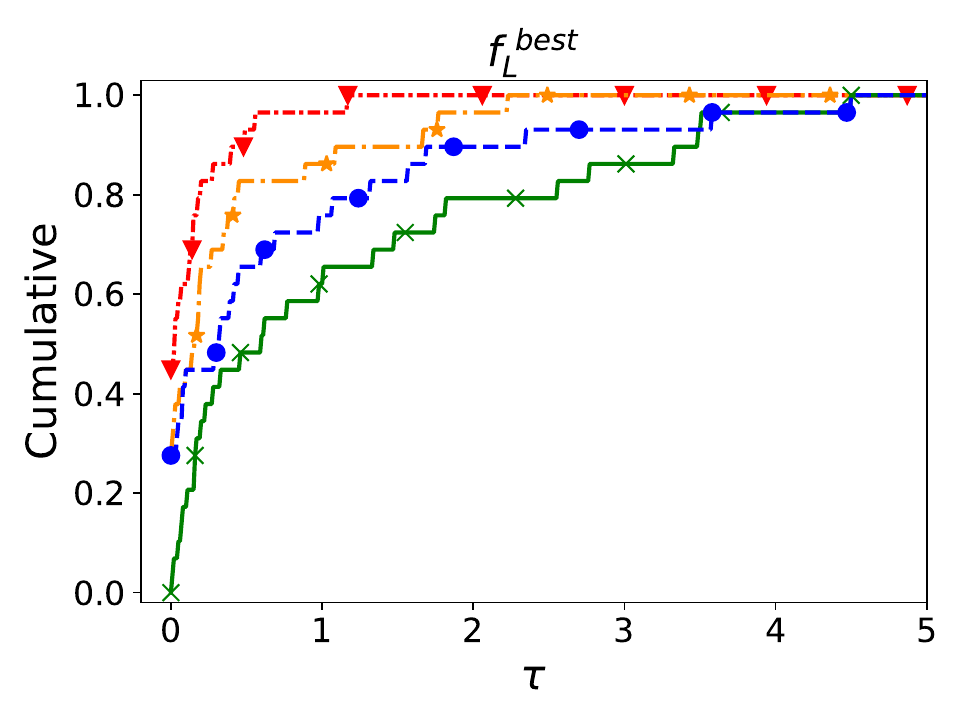}}
	\hfil
	\subfloat[$f^{best}_L$ for $q=5$, $n \ge 10$]{\includegraphics[width=0.375\textwidth]{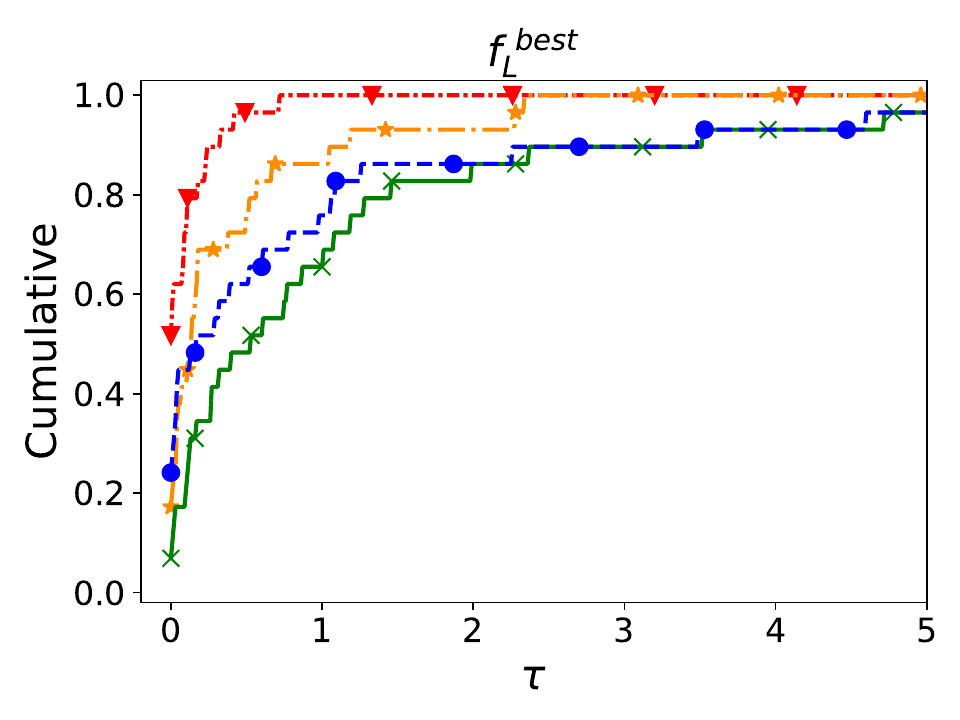}}
	\caption{Plots of the relative gap cumulative distribution for \texttt{NSMA(X)}, \texttt{NSMA(F)}, \texttt{q-EI} and \texttt{q-LCB} considering the low-dimensional ($n < 10$) and high-dimensional ($n \ge 10$) functions listed in Table \ref{tab::test_functions} and different values for the batch size parameter $q$. Note that the intervals of the x-axis are set for a better visualization of the methodologies results.}
	\label{fig::batch_pp}
\end{figure}

The plots summarize some highlights of the previous sections. In the low-dimensional problems, \texttt{q-EI} and \texttt{q-LCB} were the best methodologies in terms of both \textit{NR-AUC} and $f^{best}_L$; the performance obtained by the two standard acquisition strategies were similar. In this scenario, \texttt{NSMA(X)} and \texttt{NSMA(F)} were competitive in terms of robustness, reaching the competitors for values of $\tau \ge 2$ on the \textit{NR-AUC} metric and of $\tau \ge 6$ on the $f^{best}_L$ one. In the high-dimensional cases, the situation is the opposite, with \texttt{NSMA(X)} and \texttt{NSMA(F)} outperforming \texttt{q-EI} and \texttt{q-LCB} w.r.t.\ both the metrics. In this context, \texttt{NSMA(X)} was by far the best algorithm. On the other hand, \texttt{NSMA(F)} had the second best performance, proving that, even without any heuristic mechanism to escape from the problem boundaries (see Section \ref{subsec::clustering}), the only Pareto front reconstruction on high-dimensional instances could be more effective than a scalarization-based approach. 

The performance of the tested acquisition strategies were further investigated as the value for the batch size parameter $q$ varied. In particular, we considered $q \in \{2, 5\}$. Note that in these settings the maximum number of function evaluations remained the same, i.e., $N_M = 60$. Thus, we can easily observe that when $q = 2$ the Gaussian Process regression was performed more times than when $q=5$. For more details on the impact of $q$ on the Gaussian Process regression, the reader is referred to Section \ref{subsec::overview_bo}. For the sake of brevity, we again compare the methodologies by plotting in Figure \ref{fig::batch_pp} their relative gap cumulative distributions.

The general behavior of the acquisition strategies is confirmed with these new plots: \texttt{q-EI} and \texttt{q-LCB} performed better in the low-dimensional cases; \texttt{NSMA(X)} and \texttt{NSMA(F)} were the best algorithms for high values for $n$. Comparing the performance for different $q$, we observe that our proposals got better results as the value of $q$ increased. When $q=5$, the gap in the low dimensional scenario between the BOO based acquisition strategies and the standard ones is less marked than when $q=2$ or $q=3$; moreover, the gap obtained considering the high dimensional instances is sharper. It thus seems that our proposals become more effective as the number of points we can evaluate in parallel increases. In the Bayes-Opt context, it is difficult to have the resources to evaluate many points in parallel; however, already with values of $q \in [2, 5]$, \texttt{NSMA(X)} and \texttt{NSMA(F)} obtained interesting results.

\section{Concluding Remarks}
\label{sec::concluding-remarks}

In this paper, we dealt with Batch Bayesian Optimization (Bayes-Opt) and we proposed a novel acquisition strategy for selecting new points where to evaluate the objective function $f(\cdot)$ of the original problem \eqref{eq::bo-prob}. In Bayes-Opt, the choice of the points to evaluate is one of the most crucial operations: the evaluations are usually expensive, causing an high consumption of resources (e.g., CPU/GPU time, financial costs). Ideally, we would like to obtain a global optimum of the original problem in the minimum number of evaluations. The most employed acquisition strategies consist on selecting new points by minimizing an acquisition function: \texttt{q-EI} \eqref{eq::qei} and \texttt{q-LCB} \eqref{eq::qlcb} are certainly among the most famous ones.

In the proposed acquisition strategy, we consider a Bi-Objective Optimization (BOO) problem \eqref{eq::acq-f-biob} based on the posterior mean ($\mu_k(\cdot)$) and variance ($\sigma_k^2(\cdot)$) functions of the Gaussian Process used in the Bayes-Opt procedure. By means of a recently proposed memetic method for MOO, called \texttt{NSMA} \cite{NSMA}, we aim to reconstruct, as accurately as possible, the Pareto front of \eqref{eq::acq-f-biob}, in order to select then potentially optimal solutions for problem \eqref{eq::bo-prob}. The Pareto front reconstruction leads to the generation of multiple solutions representing different compromises between exploitation of the current information about $f(\cdot)$ (low $\mu_k(\cdot)$ values) and exploration of new regions of the problem feasible set (high $\sigma_k^2(\cdot)$ values). All these exploitation-exploration trade-offs would be difficult to reach with the employment of the standard acquisition strategies, usually based on combinations of the two posterior distribution functions. The selection of points from the Pareto front can then be performed through the two proposed clustering approaches. These latter ones work on different spaces: the first approach employs \texttt{K-MEANS} in the variables space, while the second one operates in the objectives space.

The proposed strategy was tested with thorough computational experiments, comparing its performance with the one of well-known acquisition methodologies from the Bayes-Opt literature. In the low-dimensional case our proposal turned out to be competitive, while in the high-dimensional scenario it achieved remarkable results, outperforming by far the competitors. In particular, the general superiority of \texttt{NSMA} w.r.t. \texttt{NSGA-II}, attested in the original paper \cite{NSMA}, is also confirmed in this context. The clustering approach in the objectives space generally outperformed the tested acquisition strategies from the literature, proving to be also effective on the high dimensional scenarios; this is quite remarkable as this approach does not consider the domain space and it does not have any heuristic mechanism to prevent the boundary issue, typical in the Bayesian literature in such cases. On the other hand, the clustering method in the domain space was undoubtedly the best approach when the dimensionality $n \ge 50$. Its potential effectiveness in these scenarios, analyzed in the theoretical part of the paper, is then confirmed by the computational experiments.

As possible future research directions, the following ones certainly stand out: a variant of the proposed methodology where the batch size parameter $q$ is chosen dynamically at each iteration; an effectiveness study of the presented acquisition strategy in Bayes-Opt procedures different from the standard one.

\section*{Conflict of Interest}
The authors have no competing interests to declare that are relevant to the content of this article.

\section*{Funding Sources}
No funding was received for conducting this study.

\section*{Data Availability Statement}
Data sharing not applicable to this article as no datasets were generated or analyzed during the current study.

\section*{Code availability statement}
The full code of the experiments can be found at \url{https://github.com/FranciC19/biobj_acquistion_function_for_BO}.

\bibliographystyle{abbrv}

\end{document}